\documentstyle[11pt,amssymb,amsfonts]{article}

\begin{document}
\newcommand{\p}{\parallel }
\makeatletter \makeatother
\newtheorem{th}{Theorem}[section]
\newtheorem{lem}{Lemma}[section]
\newtheorem{de}{Definition}[section]
\newtheorem{rem}{Remark}[section]
\newtheorem{cor}{Corollary}[section]
\renewcommand{\theequation}{\thesection.\arabic {equation}}

\title{{\bf The first terms in the expansion of the Bergman kernel in higher degrees: mixed curvature case}}

\author{ Yong Wang$^*$, Aihui Sun }

\date{}
\maketitle

\begin{abstract}
We establish the cancellation of the first $|2j-q|$ terms in the diagonal asymptotic
expansion of the restriction to the $(0, 2j)$-forms of the Bergman kernel associated to the modified spin$^c$
Dirac operator on high tensor powers of a line bundle with mixed curvature twisted by a (non necessarily
holomorphic) complex vector bundle, over a compact symplectic manifold. Moreover, we give a
local formula for the first and the second (non-zero) leading coefficients which generalizes the Puchol-Zhu's results.
\\

\noindent{\bf Keywords:}\quad Modified spin$^c$ Dirac operator; Bergman kernel; Asymptotic expansion
\\

\noindent {\bf MSC(2010)}: 58J35, 32L10\\

\end{abstract}

\section{Introduction}
    \quad The study of the asymptotic expansion of Bergman kernel has attracted much attention
recently. The existence of the diagonal asymptotic expansion of the Bergman
kernel of high tensor powers of a positive line bundle over a compact complex manifold
was first established by Tian [Ti], Ruan [Ru], Catlin [Ca], and Zelditch [Ze].
Tian [Ti], followed by Lu [LuZ] and Wang [Wa], derived explicit formulae for several
terms of the asymptotic expansion on the diagonal, via Tian's method of peak
sections.\\
 \indent Using Bismut-Lebeau's analytic localization techniques, Dai, Liu, and Ma [DLM] established
the full off-diagonal asymptotic expansion of the Bergman kernel of the Spinc Dirac operator associated with high powers of a Hermitian line bundle with positive curvature in the general context of symplectic manifolds. Moreover, they calculated
the second coefficient of the expansion in the case of K${\rm\ddot{a}}$hler manifolds. Later,
Ma and Marinescu [MM08] studied the expansion of generalized Bergman kernel associated
with Bochner Laplacians and developed a method of formal power series to
compute the coefficients. By the same method, Ma and Marinescu [MM06, Theorem 2.1]
computed the second coefficient of the asymptotic expansion of the Bergman kernel
of the Spin$^c$ Dirac operator acting on high tensor powers of line bundles with positive
curvature in the case of symplectic manifolds. \\
\indent Recently, this asymptotic in
the symplectic case found an application in the study of variation of Hodge structures of vector
bundles by Charbonneau and Stern in [CS]. In their setting, the Bergman kernel is the kernel
of a Kodaira-like Laplacian on a twisted bundle $L\otimes E$, where $E$ is a (not necessarily holomorphic)
complex vector bundle. Because of that, the Bergman kernel no longer concentrates in degree 0
(unlike it did in the K${\rm\ddot{a}}$hler case), and the asymptotic of its restriction to the $(0, 2j)$-forms is
related to the degree of non-holomorphicity of $E.$
In [PZ], Puchol and Zhu showed that the leading term in the asymptotic of the restriction to the
$(0, 2j)$-forms of the Bergman kernel is of order $p^{{\rm {dimX}}-2j}$ for K${\rm\ddot{a}}$hler manifolds and they computed the first and the second terms
 in this asymptotic. In [LuW1],
Lu calculated the second coefficient of the asymptotic expansion of the
Bergman kernel of the Hodge-Dolbeault operator associated with high powers of a
Hermitian line bundle with non-degenerate curvature, using the method of formal
power series developed by Ma and Marinescu. In this paper, we extend the Puchol-Zhu's results to the non-degenerate curvature case.\\
 \indent
Let $(X, J)$ be a compact connected complex manifold with complex structure $J$
and ${\rm dim}_{\mathbb{C}}X = n.$
 Let $(L,h^L)$ be a holomorphic Hermitian line bundle on $X$, and let
$\nabla^L$ be the Chern connection of $(L,h^L)$ with the curvature $R^L =(\nabla^L)^2$.\\
\indent {\bf Our basic assumption} is that $\omega=\frac{\sqrt{-1}}{2\pi}R^L$
defines a symplectic form on $X$.\\
\indent The complex structure $J$ induces a splitting $TX\otimes_{\mathbb{R}}{\mathbb{C}}=T^{(1,0)}X\oplus T^{(0,1)}X,$ where
$T^{(1,0)}X$ and $T^{(0,1)}X$ are the eigenbundles of $J$ corresponding to the eigenvalues $\sqrt{-1}$ and $-\sqrt{-1}$.
Since the $J$-invariant bilinear form $\omega(\cdot.J\cdot)$ is nondegenerate
on $TX$, there exist $J$-invariant subbundles denoted $V,V^\bot\subset TX$
 such that
$$\omega(\cdot.J\cdot)|_V<0,~~\omega(\cdot.J\cdot)|_{V^\bot}>0\eqno(1.1)$$
and $V,V^\bot$
are orthogonal with respect to $\omega(\cdot.J\cdot)$. Equivalently, there exist subbundles
$W,W^\bot\subset T^{(1,0)}X$
such that $W\oplus W^\bot=T^{(1,0)}X$, $W,W^\bot$
 orthogonal with respect to $\omega$
and
$$R^L(u,\overline{u})<0, ~~{\rm for}~~ u \in W; ~~R^L(u,\overline{u})>0, ~~{\rm for}~~ u \in W^\bot.\eqno(1.2)$$
Set ${\rm rank}W = q$. Then the curvature $R^L$ is non-degenerate of signature $(q, n- q)$.
Now take the Riemannian metric $g^{TX}$ on $X$ to be
$$g^{TX}:=-\omega(\cdot.J\cdot)|_V\oplus\omega(\cdot.J\cdot)|_{V^\bot}.\eqno(1.3)$$
Since $\omega$ is compatible with the complex structure $J$ , the metric $g^{TX}$ is also compatible
with $J$. Note that $(X, g^{TX})$ is not necessarily K${\rm \ddot{a}}$hler. Denote by $\left<\cdot,\cdot\right>$ the
${\mathbb{C}}$-bilinear form on $TX\otimes_{\mathbb{R}}{\mathbb{C}}$ induced by $g^{TX}$.
Note that $\left<\cdot,\cdot\right>$ vanishes on $T^{(1,0)}X\times T^{(1,0)}X$ and on $T^{(0,1)}X\times T^{(0,1)}X$.\\
\indent Let $\nabla^{TX}$
denote the Levi-Civita connection on $X$ and $R^{TX}$, $r^X$ denote the curvature and the scalar curvature of $(TX,g^{TX})$. The metric
$g^{TX}$ induces a Hermitian metric $h^{T^{(1,0)}X}$ on $T^{(1,0)}X$ and a metric $h^{\wedge^{0,{\bf \cdot}}}$ on
$\wedge^{0,{\bf \cdot}}(T^*X):=\wedge^{{\bf \cdot}}(T^{(0,1)}X).$ Let $\nabla^{T^{(1,0)}X}$ denote the Chern connection on $(T^{(1,0)}X,
h^{T^{(1,0)}X})$ whose curvature is $R^{T^{(1,0)}X}$ and $\nabla^{T^{(1,0)}X}$ induces a Chern connection $\nabla^{{\rm det}(T^{(0,1)}X)}$
on ${\rm det}(T^{(0,1)}X):=\wedge^nT^{(0,1)}X.$ Let $(E,h^E)$ be a Hermitian complex vector bundle with a Hermitian connection $\nabla^E$,
whose curvature is $R^E=(\nabla^E)^2$. Let $L^p=L^{\otimes p}$ be the $p^{th}$ tensor power of $L$ and $\Omega^{0,\bullet}(X,L^p\otimes E)=\Gamma(X.\wedge^{0,\bullet}(T^*X)\otimes L^p\otimes E)$. We still denote by $\left<\cdot,\cdot\right>$ be the fibrewise metric on
$\wedge^{0,\bullet}(T^*X)\otimes L^p\otimes E$ induced by $g^{TX},h^L, $ and $h^E$. Let $d\nu_X$ be the Riemannian volume of $(X,g^{TX})$.
The $L^2$-scalar product on $\Omega^{0,\bullet}(X,L^p\otimes E)$ is given by
$$\left<s_1,s_2\right>=\int_X\left<s_1(x),s_2(x)\right>d\nu_X(x).\eqno(1.4)$$
\indent The Chern connection $\nabla^{T^{(1,0)}X}$
on $T^{(1,0)}X$ induces naturally a Hermitian connection $\nabla^{T^{(0,1)}X}$ on $T^{(0,1)}X$. Set
$$\widetilde{\nabla}^{TX}=\nabla^{T^{(1,0)}X}\oplus \nabla^{T^{(0,1)}X}.\eqno(1.5)$$
Then $\widetilde{\nabla}^{TX}$ is a Hermitian connection on $TX\otimes_{\mathbb{R}}{\mathbb{C}}$ which preserves the decomposition
$TX\otimes_{\mathbb{R}}{\mathbb{C}}=T^{(1,0)}X\oplus T^{(0,1)}X$. We still denote by $\widetilde{\nabla}^{TX}$ the induced connection
on $TX$. Let $T$ be the torsion of the connection $\widetilde{\nabla}^{TX}$, and let $T_{as}$ be the
anti-symmetrization of the tensor $T$ , i.e., for $U,V,W \in TX,$
$$T_{as}(U,V,W) = \left<T(U,V),W\right>+\left<T(V,W),U\right>+\left<T(W,U),V\right> . \eqno(1.6)$$
Denote by $S^B$ the $1$-form with values in the space of antisymmetric elements of
${\rm End}(T X)$ which satisfies for $U,V,W \in TX,$
$$\left<S_B(U)V ,W\right>=-\frac{1}{2}T_{as}(U,V,W).\eqno(1.7)$$
Then the Bismut connection on $TX$ is defined by
$$\nabla^B=\nabla^{TX}+S^B.\eqno(1.8)$$
By [Bi, Prop. 2.5], $\nabla^B$ preserves the metric $g^{TX}$ and the complex structure $J$ of $TX$.\\
 \indent For any $v\in TX\otimes_{\mathbb{R}}{\mathbb{C}}$ with decomposition $v= v^{1,0} + v^{0,1}\in T ^{(1,0)}X\oplus T^{(0,1)}X,$ let
$\overline{v}^*_{1,0}$ be the $\left<\cdot,\cdot\right>$ dual of $v_{1,0}$. Then $c(v)=\sqrt{2}(\overline{v}^*_{1,0}\wedge-i_{v_{0,1}})$
defines the Clifford
action of $v$ on $\wedge(T^{*(0,1)}X)$, where $\wedge$ and $i$ denote the standard exterior and interior
multiplication, respectively. Let $\nabla^{Cl}$ denote the Clifford connection
on $\wedge(T^{*(0,1)}X)$
induced canonically by $\nabla^{TX}$ and $\nabla^{{\rm det}(T^{(1,0)}X)}$. If
$ e^1, \cdots, e^{2n}$ denotes an orthonormal frame of $T^*X$,
then define
$$^c(e^{i_1}\wedge\cdots\wedge e^{i_j})=c(e_{i_1})\cdots c(e_{i_j}), ~~{\rm for}~i_1<\cdots<i_j.\eqno(1.9)$$
In this sense $^cB$ is defined for any $B\in \wedge^*(T^*X)\otimes_{\mathbb{R}}{\mathbb{C}}$
 by extending $\mathbb{C}$-linearity. Take $U\in TX,$ Let
$$\nabla^{B,\wedge^{0,\bullet}}_U=\nabla^{Cl}_U-{\frac{1}{4}} {^c(i_UT_{as})}\eqno(1.10)$$
denote the Hermitian connection on $\wedge(T^{*(0,1)}X)$ induced by $\nabla^{Cl}$ and $T_{as}$. Then
$\nabla^{B,\wedge^{0,\bullet}}$
is the Clifford connection on the spinor bundle $\wedge(T^{*(0,1)}X)$ induced by $\nabla^B$
on $TX$ and $\nabla^{{\rm det}(T ^{(1,0)}X)}$ on ${\rm det}(T ^{(1,0)}X)$.
The connection $\nabla^{B,\wedge^{0,\bullet}}$
preserves
the $Z$-grading of $\wedge^{0,\bullet}(TX)$.
 Let $w_1, \cdots , w_n$ be an orthonormal
frame of $T ^{(1,0)}X$, Then $(\overline{w_1},\cdots,\overline{w_n})$ is a local orthonormal frame of $T ^{(0,1)}X$ whose
dual frame is denoted by $(\overline{w^1},\cdots,\overline{w^n})$ and the vectors
$$e_{2j}=\frac{1}{\sqrt{2}}(w_j+\overline{w_j}),~~e_{2j-1}=\frac{\sqrt{-1}}{\sqrt{2}}(w_j-\overline{w_j})\eqno(1.11)$$
form a local orthonormal frame of $TX$.
Set
$$\nabla^{TX}e_j=\Gamma^{TX}e_j,~~~\nabla^{{\rm det}(T ^{(1,0)}X)}(w_1\wedge\cdots\wedge w_n)=\Gamma^{{\rm det}(T ^{(1,0)}X)}
(w_1\wedge\cdots\wedge w_n).\eqno(1.12)$$
By [MM07, (1.3.5)], we have that $\nabla^{\wedge^{0,\bullet}}$, $\nabla^{B,\wedge^{0,\bullet}}$ is given, with respect to the frame $\{\overline{w}^{j_1}\wedge\cdots\wedge
\overline{w}^{j_k},~~1\leq j_1< \cdots< j_k\leq n\}$ of $\wedge(T^{*(0,1)}X)$, by the local formula respectively
$$d+\frac{1}{4}\left<\Gamma^{TX}e_i,e_j\right>c(e_i)c(e_j)+\frac{1}{2}\Gamma^{{\rm det}(T ^{(1,0)}X)};$$
$$d+\frac{1}{4}\left<\Gamma^{TX}e_i,e_j\right>c(e_i)c(e_j)+\frac{1}{2}\Gamma^{{\rm det}(T ^{(1,0)}X)}-\frac{1}{4}
{^c(i_{\bullet}T_{as})}.\eqno(1.13)$$
Let $\Gamma^{B,\wedge^{0,\bullet}}$ be the connection $1$-form of $\nabla^{B,\wedge^{0,\bullet}}$, i.e.,
$$\Gamma^{B,\wedge^{0,\bullet}}=\frac{1}{4}\left<\Gamma^{TX}e_i,e_j\right>c(e_i)c(e_j)+\frac{1}{2}\Gamma^{{\rm det}(T ^{(1,0)}X)}-\frac{1}{4}
{^c(i_{\bullet}T_{as})}.\eqno(1.14)$$
Denote by $\nabla^{L^p\otimes E}$ the Hermitian connection on $L^p\otimes E$ induced by $\nabla^L$ and $\nabla^E$. Set
$$\nabla^{\wedge^{0,\bullet}\otimes L^p\otimes E}=\nabla^{\wedge^{0,\bullet}}\otimes 1+1\otimes
\nabla^{L^p\otimes E};~\nabla^{B,\wedge^{0,\bullet}\otimes L^p\otimes E}=\nabla^{B,\wedge^{0,\bullet}}\otimes 1+1\otimes
\nabla^{L^p\otimes E}.\eqno(1.15)$$
Then $\nabla^{\wedge^{0,\bullet}\otimes L^p\otimes E}$ and $\nabla^{B,\wedge^{0,\bullet}\otimes L^p\otimes E}$ are Hermitian connections on $\wedge^{0,\bullet}(TX)\otimes L^p\otimes E$.
Define the modified spin$^c$ Dirac operator by
$$D^{B}_p:=\sum_{j=1}^{2n}c(e_j)\nabla^{\wedge^{0,\bullet}\otimes L^p\otimes E}_{e_j}-\frac{1}{4}c(T_{as}),\eqno(1.16)$$
which is a first order elliptic self-adjoint differential operator.\\

\noindent {\bf Definition 1.1} Let
$$P_p:\Omega^{(0,\bullet)}(X,L^p\otimes E)\rightarrow {\rm ker}(D^B_p)\eqno(1.15)$$
be the orthogonal projection onto the kernel ${\rm ker}(D^B_p)$ of $D^B_p$. The operator $Pp$ is called the
Bergman projection. It has a smooth kernel with respect to $d\nu_X(y)$, denoted by $Pp(x, y)$, which
is called the {\bf Bergman kernel}.\\

\indent We recall\\

\noindent {\bf Theorem 1.2} ([MM07, Thm. 8.2.4]) {\it There exist smooth sections ${\bf b}_r$ of
${\rm End}(\wedge^{0,{\rm even}}(T^*X)\otimes E)$ such that for any $k\in \mathbb{N}$ and for $p \rightarrow +\infty$:
$$ p^{-n}P_p(x,x)=\sum_{r=0}^k{\bf b}_r(x)p^{-r}+O(p^{-k-1}),\eqno(1.17)$$
that is for every $k, l\in \mathbb{N}$, there exists a constant $C_{k,l} > 0$ such that for any $p\in\mathbb{N}$,
$$\left|p^{-n}P_p(x,x)-\sum_{r=0}^k{\bf b}_r(x)p^{-r}\right|_{L^l(X)}\leq C_{k,l}p^{-k-1}.\eqno(1.18)$$
Here $| \cdot |_{L^l(X)}$ is the $L^l$-norm for the variable $x \in X$.}\\

Let ${\mathcal{R}}= (R^E)^{0,2} \in\Omega^{(0,2)}(X,{\rm End}(E))$
be the $(0, 2)$-part of $R^E$ (which is zero if $E$ is holomorphic). Let
$${\mathcal{R}}={\mathcal{R}}^\top+{\mathcal{R}}^0+{\mathcal{R}}^\bot,\eqno(1.19)$$
where ${\mathcal{R}}^\top\in \Gamma(X,\wedge^2({\overline{W}}^*)\otimes {\rm End}(E))$, ${\mathcal{R}}^0\in \Gamma(X,{\overline{W}}^*\otimes
{{\overline{W}}^*}^\bot\otimes {\rm End}(E))$, ${\mathcal{R}}^\bot\in \Gamma(X,\wedge^2({{\overline{W}}^*}^\bot)\otimes {\rm End}(E))$.
For $1\leq j\leq n$, let
$$I_j:\wedge^{0,\bullet}(T^*X)\otimes E\rightarrow \wedge^{0,j}(T^*X)\otimes E\eqno(1.20)$$
be the natural orthogonal projection. Let $I_{{\rm det}(\overline{W}^*)\otimes E}$ denote the projection on ${\rm det}(\overline{W}^*)\otimes E$.
For $j,k,l\in\mathbb{N}$, we define $B^{k,j}_l$ for $k\leq j$ by
$$B^{k,j}_l:=\frac{1}{(2k+l)\times\cdots\times (2j+l)};~~B^{k,0}_l=1.\eqno(1.21)$$
Then we have\\

\noindent {\bf Theorem 1.3} {\it For any $k\in\mathbb{N}$, $k\geq 2j$, we have when $p \rightarrow +\infty$:
$$p^{-n}I_{2j}P_p(x,x)I_{2j}=\sum_{r=2j}^kI_{2j}{\bf b}_{|r-q|}(x)I_{2j}p^{-|r-q|}+O(p^{-|k-q|-1}),\eqno(1.22)$$
and moreover, when $2j\geq q$
$$I_{2j}{\bf b}_{2j-q}(x)I_{2j}=\frac{1}{(4\pi)^{2j-q}}(B^{1,j-\frac{q}{2}}_0)^2I_{2j}({{\mathcal{R}}_x^\bot})^{j-\frac{q}{2}}
I_{{\rm det}(\overline{W}^*)\otimes E}({{\mathcal{R}}_x^{\bot,*}})^{j-\frac{q}{2}}I_{2j},\eqno(1.23)$$
where ${{\mathcal{R}}_x^{\bot,*}} $ is the dual of ${{\mathcal{R}}_x^\bot}$ acting on $(\wedge^{0,\bullet}(T^*X)\otimes E)_x$.
When $2j<q$
$$I_{2j}{\bf b}_{q-2j}(x)I_{2j}=\frac{1}{(4\pi)^{q-2j}}(B^{1,\frac{q}{2}-j}_0)^2I_{2j}({{\mathcal{R}}_x^{\top,*}})^{\frac{q}{2}-j}
I_{{\rm det}(\overline{W}^*)\otimes E}({{\mathcal{R}}_x^\top})^{\frac{q}{2}-j}I_{2j}.\eqno(1.24)$$
}\\

Let ${\bf J} :TX\rightarrow TX$ be
the almost complex structure defined by
$$\omega(U,V)=g^{TX}({\bf J}U,V)~~{\rm for}~~~U,V\in TX.\eqno(1.25)$$
Then $J$ commutes with ${\bf J}$. Let $w_1, \cdots , w_n$ be an orthonormal frame of $(T^{(1,0)}X,h^{T^{(1,0)}X})$ such that the subbundle $W$ is spanned by $w_1,\cdots , w_q$, and let $w^1,\cdots,w^q$ be
the dual frame. Then
$${\bf J}w_j=-\sqrt{-1}w_j,~~{\rm for}~~j\leq q;~~{\bf J}w_j=\sqrt{-1}w_j,~~{\rm for}~~j\geq q+1.\eqno(1.26)$$
Let $T^{(1,0)}_{{\bf J}}X$ and  $T^{(0,1)}_{{\bf J}}X$
 be the eigenbundles of ${\bf J}$ corresponding to the eigenvalues $\sqrt{-1}$ and $-\sqrt{-1}$
respectively. Set
$$u_l = \overline{w}_l ~~{\rm if}~~ l \leq q ~~{\rm and}~~ u_l= w_l ~~{\rm otherwise}.\eqno(1.27) $$
Then $u_1, \cdots , u_n$ forms an orthonormal frame of the subbundle $T^{(1,0)}_{{\bf J}}X$
We denote
by $u^1,\cdots, u^n$ its dual frame. Then
$$\omega=\sqrt{-1}\sum_{l=1}^nu^l\wedge \overline{u}^l.\eqno(1.28)$$
For $m,k,l\in\mathbb{N}$, we define $C_m(k)$ for $k\leq m$ by
$$~~C_m(k):=\frac{1}{(4\pi)^{m}}\frac{1}{2^kk!}\frac{1}{\prod_{s=k+1}^m(2s+1)},
\eqno(1.29)$$
with the convention that $\prod_{s\in\emptyset}=1$.
Let $\triangle^{\wedge^{0,\bullet}\otimes E}$ be the Laplacian induced by $\nabla^{\wedge^{0,\bullet}\otimes E}.$ and
 $\widetilde{J}=-2\pi\sqrt{-1}{\bf J}$. Denote by $R^{\wedge^{0,\bullet},B}$ and $R^{\rm det}$ the curvature associated to the connections
 $\nabla^{\wedge^{0,\bullet},B}$ and $\nabla^{\rm det}$. If $e_1,\cdots,e_{2n}$ denotes an orthonormal frame of $TX$, then set
 $$|A_1|^2=\sum_{i<j<k}|A_1(e_i,e_j,e_k)|^2 ~~{\rm for}~A_1\in \wedge^3(T^*X).\eqno(1.30)$$
 Define
$$A(Y,U,V):=\left<(\nabla^{TX}_{Y}\widetilde{J})U,V\right>;~~(\nabla^X\nabla^X\widetilde{J})_{(U,V)}=\nabla^X_U\nabla^X_V(\widetilde{J})
-\nabla^X_{\nabla^X_UV}(\widetilde{J});\eqno(1.31)$$
$$
\Psi(X_0,Y,U,V)=-\frac{1}{9}\left<(\nabla^{TX}_{X_0}\widetilde{J})Y,(\nabla^{TX}_{U}\widetilde{J})V\right>;\eqno(1.32)$$
$$\widehat{\Psi}(X_0,Y,U,V)=\frac{1}{2}\left<(\nabla^{TX}\nabla^{TX})_{(X_0,Y)}U,V\right>$$
$$+\frac{1}{6}\left[\left
<R^{TX}(X_0,V)Y,\widetilde{J}U\right>-\left
<R^{TX}(X_0,U)Y,\widetilde{J}V\right>\right];\eqno(1.33)$$
$$\Phi(U,V)=\frac{1}{2}\left<(\nabla^X\nabla^X\widetilde{J})_{(U,V)}{ u_l},\overline{u}_l\right>.\eqno(1.34)$$
Write
$$A_{jl\overline{r}}=A(u_j,u_l,\overline{u_r});~~A_{e_rle_r}=A(e_r,u_l,e_r);~~$$
$${\Psi}_{jl\overline{r}\overline{s}}=\Psi(u_j,u_l,\overline{u_r},\overline{u_s});~~
\widehat{\Psi}_{jl\overline{r}\overline{s}}=\widehat{\Psi}(u_j,u_l,\overline{u_r},\overline{u_s}).\eqno(1.35)$$\\

\noindent {\bf Theorem 1.4} {\it When $n=4$, $j=2$ and $q=2$, we have}
$$I_4{\bf b}_3(x)I_4=I+II+III+IV+V,\eqno(1.36)$$
{\it where}
$${\rm I}=\frac{1}{1152\pi^3}I_4
(\nabla^{\wedge^{0,\bullet}\otimes E}_{\overline{u_l}}{{\mathcal{R}}}_{.})^\bot(x)I_{{\rm det}(\overline{W}^*)\otimes E}
\left[(\nabla^{\wedge^{0,\bullet}\otimes E}_{\overline{u_l}}{{\mathcal{R}}}_{.})^\bot(x)\right]^*I_4,\eqno(1.37)$$

$${\rm II}=\frac{1}{768\pi^3}I_4
(\triangle^{\wedge^{0,\bullet}\otimes E}{{\mathcal{R}}}_{.})^\bot(x)I_{{\rm det}(\overline{W}^*)\otimes E}({\mathcal{R}}_x^\bot)^*I_4$$
$$
+\frac{1}{768\pi^3}I_4{\mathcal{R}}_x^\bot I_{{\rm det}(\overline{W}^*)\otimes E}(\triangle^{\wedge^{0,\bullet}\otimes E}{{\mathcal{R}}}_{.})^{\bot,*}(x)I_4.
\eqno(1.38)$$\\
{\it III is determined by (4.11), (4.34), (4.37), (4.76), (4.79)-(4.83). IV is determined by (4.198), (4.212), (4.214), (4.215), (4.222), (4.223), (4.225), (4.226),
(4.229), (4.230), (4.232). V is determined by
 (4.86), (4.151)-(4.153),
(4.159), (4.162)-(4.172), (4.182), (4.183), (4.189), (4.194)-(4.197). }\\

\indent This paper is organized as follows: in Section 2 we use a
local trivialization to rescale $(D^B_p)^2$, and then give the Taylor expansion of the rescaled operator. In
Section 3, we use this expansion to give a formula for the coefficients ${\bf b}_r$. Finally, in Section 4, we prove Theorem 1.4. In this whole paper, when an index variables appears twice in a single term, it
means that we are summing over all its possible values.

\section{Rescaling $(D^B_p)^2$ and Taylor expansion }

 \quad In this Section, we will give a formula for the square of $(D^B_p)^2$,
then rescale the operator $(D^B_p)^2$ to get an operator ${\mathcal{L}}_t$, and give the Taylor
expansion of the rescaled operator. We also study more precisely the limit operator ${\mathcal{L}}_0$.\\
\indent Let the Laplacian
$$\triangle^{B,\wedge^{0,\bullet}\otimes L^p\otimes E}=-\sum_{j=1}^{2n}\left[(\nabla^{B,\wedge^{0,\bullet}\otimes L^p\otimes E}_{e_j})^2
-\nabla^{B,\wedge^{0,\bullet}\otimes L^p\otimes E}_{\nabla^{TX}_{e_j}e_j}\right].\eqno(2.1)$$
Then by [MM07, Thm 1.3.7], we have
$$(D^B_p)^2= \triangle^{B,\wedge^{0,\bullet}\otimes L^p\otimes E}+\frac{1}{2}pR^L(e_i,e_j)c(e_i)c(e_j)+\frac{r^X}{4}
+^c(R^E+\frac{1}{2}R^{\rm det})-\frac{1}{4}c(dT_{as})-\frac{1}{8}|T_{as}|^2.\eqno(2.2)$$
 For any skew-adjoint endomorphism $A$ of $TX$, it holds that
 $$\frac{1}{4}\left<Ae_i,e_j\right>c(e_i)c(e_j)=-\frac{1}{2}\left<Aw_j,\overline{w}_j\right>
 +\left<Aw_l,\overline{w}_m\right>\overline{w}^m\wedge i_{\overline{w}_l}$$
 $$+\frac{1}{2}
 \left<Aw_l,{w_m}\right>i_{\overline{w}_l}i_{\overline{w}_m}
 +\frac{1}{2}\left<A\overline{w}_l,\overline{w}_m\right>\overline{w}^l\wedge\overline{w}^m\wedge.\eqno(2.3)$$
 On $\Omega^{0,\bullet}(X)$, we define the operator $\omega_{d,x}$ by
 $$\omega_{d,x}=-2\pi\sum_{l=1}^qi_{\overline{w}_l}\wedge \overline{w}^l-2\pi\sum_{l=q+1}^n\overline{w}^l\wedge i_{\overline{w}_l}.\eqno(2.4)$$
By (2.2)-(2.4) and (1.27), (1.28) and the fundamental assumption, we get
\\

 \noindent {\bf Proposition 2.1} {\it It holds that}
 $$(D^B_p)^2= \triangle^{B,\wedge^{0,\bullet}\otimes L^p\otimes E}-(R^E+\frac{1}{2}R^{\rm det})(w_l,\overline{w}_l)
 +\frac{r^X}{4}-2\pi pn-2p\omega_d-\frac{1}{4}c(dT_{as})$$
 $$-\frac{1}{8}|T_{as}|^2
 +2(R^E+\frac{1}{2}R^{\rm det})(w_l,\overline{w}_m)\overline{w}^m\wedge i_{\overline{w}_l}
 +R^E(w_l,{w_m})i_{\overline{w}_l}i_{\overline{w}_m}
 +R^E(\overline{w}_l,\overline{w}_m)\overline{w}^l\wedge\overline{w}^m\wedge.\eqno(2.5)$$\\

 \indent Fix $x_0 \in X$ and ${w_j}$ an orthonormal basis of $T^{(1,0)}_{x_0}X$, with
dual basis ${w^j}$, and we construct an orthonormal basis $\{e_l\}$ of $T_{x_0}X$ from $\{w_j\}$ as in (1.11).
For $\varepsilon > 0$, we denote by $B^X(x_0,\varepsilon)$ and $B^{T_{x_0}X}(0,\varepsilon)$ the open balls in $X$ and $T_{x_0}X$ with center
$x_0$ and $0$ and radius $\varepsilon$ respectively. If ${\rm exp}^X_{x_0}$
is the Riemannian exponential of $X$, then for $\varepsilon$
small enough, $Z\in B^{T_{x_0}X}(0,\varepsilon)\rightarrow {\rm exp}^X_{x_0}(Z)\in B^X(x_0,\varepsilon)$
is a diffeomorphism, which gives local
coordinates by identifying $T_{x_0}X$ with ${\mathbb{R}}^{2n}$ via the orthonormal basis $\{e_l\}$:
$$(Z_1,\cdots,Z_{2n})\in {\mathbb{R}}^{2n}\mapsto \sum_{l}Z_le_l\in T_{x_0}X.\eqno(2.6)$$
From now on, we will always identify $B^{T_{x_0}X}(0,\varepsilon)$ and $B^X(x_0,\varepsilon)$. Note that in this identification,
the radial vector field ${\mathcal{R}}=\sum_lZ_le_l$ becomes ${\mathcal{R}}=Z$,
 so $Z$ can be viewed as a point or as a tangent
vector.\\
\indent For $Z \in B^{T_{x_0}X}(0,\varepsilon)$, we identify $(L_Z, h^L_Z),~(E_Z,h^E_Z)$ and $(\wedge^{0,\bullet}(T^*X)_Z,h^{\wedge^{0,\bullet}}_{Z})$ with $(L_{x_0}, h^L_{x_0}),~(E_{x_0},h^E_{x_0})$ and $(\wedge^{0,\bullet}(T^*X)_{x_0},h^{\wedge^{0,\bullet}}_{x_0})$
 by parallel transport with respect to the connection $\nabla^L,~\nabla^E$ and $\nabla^{\wedge^{0,\bullet},B}$.
Let $S_L$ be a unit vector of $L_{x_0}$. It gives an isometry $L^p_{x_0}\simeq {\mathbb{C}}$,
 which yields to an isometry ${\mathbb{E}}_p:=\wedge^{0,\bullet}(T^*X)\otimes L^p\otimes E\simeq \wedge^{0,\bullet}(T^*X)\otimes E$.
Let $d\nu_{TX}$ be the Riemannian volume form of $(T_{x_0}X, g^{T_{x_0}X})$, and $\kappa(Z)$ be the smooth positive
function defined for $|Z|\leq \varepsilon$ by
$$ d\nu_X(Z) = \kappa(Z)d\nu_{TX}(Z),\eqno(2.7)$$
with $\kappa(0) = 1$.\\

\noindent{\bf Definition 2.2} We denote by $\nabla_U$ the ordinary differentiation operator in the direction $U$ on
$T_{x_0}X$. For $s\in\Gamma({\mathbb{R}}^{2n},{\mathbb{E}}_{x_0})$, and for $t=\frac{1}{\sqrt{p}}$, set
 $$(S_ts)(Z)=s(Z/t),$$
 $$\nabla^B_t=tS^{-1}_t\kappa^{\frac{1}{2}}\nabla^{Cl_0,B}\kappa^{-\frac{1}{2}}S_t,$$
 $$\nabla_{0,.}=\nabla_.+\frac{1}{2}R^L_{x_0}(Z,\cdot),\eqno(2.8)$$
 $${\mathcal{L}}_t=t^2S^{-1}_t\kappa^{\frac{1}{2}}(D^B_p)^2\kappa^{-\frac{1}{2}}S_t,$$
 $${\mathcal{L}}_0=-(\nabla_{0,e_j})^2-2n\pi-2\omega_{d,x_0}.$$

 Let $||\cdot||_{L^2}$ be the $L^2$-norm induced by $h^{{\mathbb{E}}_{x_0}}$ and $d\nu_{TX}$, we have\\

\noindent{\bf Theorem 2.3} {\it There exist second order formally self-adjoint (with respect to $||\cdot||_{L^2}$) differential
operators ${\mathcal{O}}_r$ with polynomial coefficients such that for all $m\in {\mathbb{N}}$,
$${\mathcal{L}}_t={\mathcal{L}}_0+\sum_{r=1}^mt^r{\mathcal{O}}_r+O(t^{m+1})\eqno(2.9)$$
Furthermore, each ${\mathcal{O}}_r$ can be decomposed as
$${\mathcal{O}}_r={\mathcal{O}}^0_r+{\mathcal{O}}^{+2}_r+{\mathcal{O}}^{-2}_r,\eqno(2.10)$$
where ${\mathcal{O}}^k_r$ changes the degree of the form which it acts on by $k$.}\\

 \noindent{\bf Proof.} By (5.22) in [LuW1], we have
 $$\frac{1}{4}{^c(dT_{as})}=\frac{1}{8}dT_{as}(w_i,\overline{w}_i,w_j,\overline{w}_j)
 -\frac{1}{2}dT_{as}(w_i,\overline{w}_j,w_k,\overline{w}_k)\overline{w}^j\wedge i_{\overline{w}_i}$$
 $$
 +\frac{1}{4}{dT_{as}}
 (w_i,w_j,\overline{w}_k,\overline{w}_l)\overline{w}^k\wedge\overline{w}^l\wedge i_{\overline{w}_i}i_{\overline{w}_j}.\eqno(2.11)$$
 So $\frac{1}{4}{^c(dT_{as})}$ preserves the degree of $\wedge^{0,\bullet}(T^*X)$. By $\nabla^B(J)=0$, we know that $\nabla^{Cl,B}$
 also preserves the degree of $\wedge^{0,\bullet}(T^*X)$. By Proposition 2.1, similar to the proof of Theorem 1.8 in [PZ], we can prove this theorem. $\Box$\\

\indent Recall Theorem 3.5.1 in [LuW2], we may get\\

\noindent{\bf Theorem 2.4} {\it It holds that}
$${\mathcal{O}}_1=-\frac{2}{3}(\partial_sR^L)_{x_0}({\mathcal{R}},e_i)Z_s\nabla_{0,e_i}-\frac{1}{3}(\partial_sR^L)_{x_0}({\mathcal{R}},e_s)
-\pi\sqrt{-1}\left<(\nabla^B_{{\mathcal{R}}}{\bf J})e_i,e_l\right>c(e_i)c(e_l),
\eqno(2.12)$$
\begin{eqnarray*}
{\mathcal{O}}^0_2&=&
\frac{1}{3}\left<R^{TX}_{x_0}({\mathcal{R}},e_i){\mathcal{R}},e_s\right>\nabla_{0,e_i}\nabla_{0,e_s}
+\left[\frac{2}{3}\left<R^{TX}_{x_0}({\mathcal{R}},e_s)e_s,e_i\right>\right.\\
&&\left.-\left(\frac{1}{2}\sum_{|\alpha|=2}(\partial_\alpha R^L)_{x_0}\frac{Z^\alpha}{\alpha!}
+R^E_{x_0}+R^{\wedge^{0,\bullet},B}_{x_0}\right)({\mathcal{R}},e_i)\right]\nabla_{0,e_i}\\
&&-\frac{1}{4}\nabla_{e_s}\left(\sum_{|\alpha|=2}(\partial_\alpha R^L)_{x_0}
({\mathcal{R}},e_s)\frac{Z^\alpha}{\alpha!}\right)\\
&&-\frac{1}{9}\sum_i\left[\sum_s(\partial_sR^L)_{x_0}({\mathcal{R}},e_i)Z_s\right]^2\\
&&+\left[-(R^E+\frac{1}{2}R^{\rm det})(w_l,\overline{w}_l)
 +\frac{r^X}{4}-\frac{1}{4}c(dT_{as})\right.\\
 &&\left.-\frac{1}{8}|T_{as}|^2
 +2(R^E+\frac{1}{2}R^{\rm det})(w_l,\overline{w}_m)\overline{w}^m\wedge i_{\overline{w}_l}
 \right](x_0)\\
 &&+\frac{1}{12}\left[\sum_s(\nabla_{0,e_s})^2,\left<R^{TX}_{x_0}({\mathcal{R}},e_l){\mathcal{R}},e_l\right>\right]\\
 &&-\frac{\pi}{2}\sqrt{-1}\left<(\nabla^B\nabla^B{\bf J})_{({{\mathcal{R}}},{{\mathcal{R}}})}e_i,e_l\right>c(e_i)c(e_l)
 ~~~~~~~~~~~~~~~~~~~~~~~~~~~~~~~~~~~~~~~~(2.13)
\end{eqnarray*}
$${\mathcal{O}}^{+2}_2={\mathcal{R}}_{x_0},~~~{\mathcal{O}}^{-2}_2=({\mathcal{R}}_{x_0})^*.\eqno(2.14)$$\\

\indent In the following, we will study the limit operator ${\mathcal{L}}_0$.
We introduce the complex coordinates $\xi=(\xi_1,\cdots,\xi_n)$ on ${\mathbb{C}}^n\simeq {\mathbb{R}}^{2n}$. We get $Z=\xi+\overline{\xi}$,
$w_j=\sqrt{2}\frac{\partial}{\partial \xi_j}$ and $\overline{w}_j=\sqrt{2}\frac{\partial}{\partial \overline{\xi}_j}$. We will identify $\xi$ to
$\sum_j\xi_j\frac{\partial}{\partial \xi_j}$ and $\overline{\xi}$ to $\sum_j\overline{\xi}_j\frac{\partial}{\partial \overline{\xi}_j}$ when we consider
$\xi$ and $\overline{\xi}$ as vector fields. Set
$$z=(\overline{\xi}_1,\cdots,\overline{\xi}_q,\xi_{q+1},\cdots,\xi_n),~~\overline{z}=({\xi}_1,\cdots,{\xi}_q,\overline{\xi}_{q+1},\cdots,
\overline{\xi}_n).\eqno(2.15)$$
Then $${\bf J}\partial_{z_l}=\sqrt{-1}\partial_{z_l},~~{\bf J}\partial_{\overline{z}_l}=-\sqrt{-1}\partial_{\overline{z}_l},~~
{\rm for} ~~l=1,\cdots,n.\eqno(2.16)$$
We will identify $z$ to
$\sum_jz_j\frac{\partial}{\partial z_j}$ and $\overline{z}$ to $\sum_j\overline{z}_j\frac{\partial}{\partial \overline{z}_j}$ when we consider
$z$ and $\overline{z}$ as vector fields. Then $Z=\xi+\overline{\xi}=z+\overline{z}$. Set $u_j=\sqrt{2}\partial_{z_j}$ and
$$f_{2j}=\frac{1}{\sqrt{2}}(u_j+\overline{u_j}),~~f_{2j-1}=\frac{\sqrt{-1}}{\sqrt{2}}(u_j-\overline{u_j}).\eqno(2.17)$$
Then $\{u_1,\cdots,u_n\}$ forms an orthonormal basis of $T^{(1,0)}_{{\bf J},x_0}X$ and ${f_1,\cdots , f_{2n}}$ is an
orthonormal basis of $T_{x_0}X$.
Set
$$b_i=-2\nabla_{0,\frac{\partial}{\partial z_i}},~~b^+_i=2\nabla_{0,\frac{\partial}{\partial \overline{z}_i}},$$
$$b=(b_1,\cdots,b_n),~~{\mathcal{L}}=-\sum_i(\nabla_{0,e_i})^2-2\pi n.\eqno(2.18)$$
By definition, $\nabla_0=\nabla+\frac{1}{2}R^L_{x_0}(Z,\cdot)$, so we get
$$b_i=-2\frac{\partial}{\partial z_i}+\pi\overline{z}_i,~~b^+_i=2\frac{\partial}{\partial \overline{z}_i}+\pi{z_i},\eqno(2.19)$$
and for any polynomial $g(z,\overline{z})$ in $z$ and $\overline{z}$,
$$[b_i,b^+_j]=-4\pi \delta_{ij},~~~[b_i,b_j]=[b^+_i,b^+_j]=0,$$
$$[g(z,\overline{z}),b_j]=2\frac{\partial}{\partial z_j}g(z,\overline{z}),~~~[g(z,\overline{z}),b^+_j]=-2\frac{\partial}{\partial \overline{z}_j}g(z,\overline{z}).\eqno(2.20)$$
and
$${\mathcal{L}}=\sum^n_{l=1}b_lb^+_l,~~~{\mathcal{L}}_0={\mathcal{L}}-2\omega_{d,x_0}.\eqno(2.21)$$
Then $b^+_l=(b_l)^*$, and ${\mathcal{L}},~{\mathcal{L}}_0$ are self-adjoint with respect to the $L^2$ norm.
We recall\\

\noindent{\bf Theorem 2.5} ([MM07, Thm. 8.2.3]) {\it The spectrum of the restriction of ${\mathcal{L }}$ to $L^2({\mathbb{R}}^{2n})$
 is ${\rm Sp}({\mathcal{L }}|_{L^2({\mathbb{R}}^{2n})})= 4\pi{\mathbb{N}}$ and
an orthogonal basis of the eigenspace for the eigenvalue $4\pi k$ is }
$$b^\alpha\left(z^\beta{\rm exp}\left(-\frac{\pi}{2}|z|^2\right)\right),~~~{\rm with}~~\alpha,\beta\in {\mathbb{N}}^n~~{\rm and} ~~\sum_i\alpha_i=k.\eqno(2.22)$$\\

\indent Especially, an orthonormal basis of ${\rm ker}({\mathcal{L }}|_{L^2({\mathbb{R}}^{2n})})$ is
$$\left(\frac{\pi^{|\beta|}}{\beta!}\right)^{\frac{1}{2}}z^\beta{\rm exp}\left(-\frac{\pi}{2}|z|^2\right),\eqno(2.23)$$
and thus if ${\mathcal{P}}(Z,Z')$ is the smooth kernel of ${\mathcal{P}}$ the orthogonal projection from $(L^2({\mathbb{R}}^{2n}),||\cdot||_0)$
onto ${\rm ker}({\mathcal{L}})$ (where $||\cdot ||_0$ is the $L^2$-norm associated to $g^{TX}_{x_0}$)
with respect to $d\nu_{TX}(Z')$, we have
$${\mathcal{P}}(Z,Z¡ä)={\rm exp}\left(-\frac{\pi}{2}(|z|^2+|z'|^2-2zz')\right).\eqno(2.24)$$
Now let $P^N$ be the orthogonal projection from $(L^2({\mathbb{R}}^{2n}, {\mathbb{E}}_{x_0}), ||\cdot||_{L^2})$ onto $N := {\rm ker}({\mathcal{L }}_0)$, and
$P^N(Z,Z')$ be its smooth kernel with respect to $d\nu_{TX}(Z')$. From (2.21), we have:
$$P^N(Z,Z')={\mathcal{P}}(Z,Z')I_{{\rm det}(\overline{W}^*)\otimes E}.\eqno(2.25)$$\\

\section{The first coefficient in the asymptotic expansion}

\quad In this section, we will compute the first coefficient in the asymptotic expansion. By Theorem 2.5 and (2.21), we get that
for every $\lambda\in S^1$ the unit circle in ${\mathbb{C}}$, $(\lambda-{\mathcal{L}}_0)^{-1}$ exists. Let $f(\lambda,t)$ be a formal power series on $t$ with values in ${\rm End}(L^2({\mathbb{R}}^{2n},{\mathbb{E}}_{x_0}))$:
$$f(\lambda,t)=\sum_{r=0}^{+\infty}t^rf_r(\lambda)~~{\rm with}~f_r(\lambda)\in {\rm End}(L^2({\mathbb{R}}^{2n},{\mathbb{E}}_{x_0})).\eqno(3.1)$$
\indent Consider the equation of formal power series on $t$ for $\lambda\in S^1$:
$$\left(\lambda-{\mathcal{L}}_0-\sum_{r=1}^{+\infty}t^r{\mathcal{O}}_r\right)f(\lambda,t)={\rm Id}_{L^2({\mathbb{R}}^{2n},{\mathbb{E}}_{x_0})}.\eqno(3.2)$$
By induction, we get
$$f_r(\lambda)=\left(\sum_{\begin{array}{lcr}
  r_1+\cdots+r_k=r  \\
  ~~~~~~ r_j\geq 1
\end{array}}
(\lambda-{\mathcal{L}}_0)^{-1}{\mathcal{O}}_{r_1}\cdots(\lambda-{\mathcal{L}}_0)^{-1}
{\mathcal{O}}_{r_k}\right)(\lambda-{\mathcal{L}}_0)^{-1}.\eqno(3.3)$$
We define ${\mathcal{F}}_r$ by
$${\mathcal{F}}_r=\frac{1}{2\pi\sqrt{-1}}\int_{S^1}f_r(\lambda)d\lambda,\eqno(3.4)$$
and we denote by ${\mathcal{F}}_r(Z,Z')$ its smooth kernel with respect to $d\nu_{TX}(Z')$.\\

{\noindent{\bf Theorem 3.1}([MM07, Thm 8.2.4) {\it The following equation holds}
$${\bf b}_r(x_0)={\mathcal{F}}_{2r}(0,0).\eqno(3.5)$$\\

\noindent{\bf Proof of Theorem 1.3.} Let $T_{\bf r}(\lambda)=(\lambda-{\mathcal{L}}_0)^{-1}{\mathcal{O}}_{r_1}\cdots
(\lambda-{\mathcal{L}}_0)^{-1}{\mathcal{O}}_{r_k}(\lambda-{\mathcal{L}}_0)^{-1}$ be the term in the sum (3.3) corresponding to
${\bf r}=(r_1,\cdots,r_k)$. Let $N^{\bot}$ be the orthogonal of $N$ in $L^2({\mathbb{R}}^{2n},{\mathbb{E}}_{x_0})$, and $P^{N^\bot}$
 be the associated orthogonal projector. In $T_r(\lambda)$, each term $(\lambda-{\mathcal{L}}_0)^{-1}$ can be decomposed as
 $$(\lambda-{\mathcal{L}}_0)^{-1}=(\lambda-{\mathcal{L}}_0)^{-1}P^{N^\bot}+\frac{1}{\lambda}P^N.\eqno(3.6)$$
 Set
 $$L^{N^{\bot}}(\lambda)=(\lambda-{\mathcal{L}}_0)^{-1}P^{N^\bot},~~L^{N}(\lambda)=\frac{1}{\lambda}P^N.\eqno(3.7)$$
 Then $(\lambda-{\mathcal{L}}_0)^{-1},~L^{N^{\bot}}(\lambda),~L^{N}(\lambda)$ preserve the degree.
 For $\eta=(\eta_1,\cdots,\eta_{k+1})\in \{N,N^{\bot}\}^{k+1}$, let
 $$T^\eta_{\bf r}(\lambda)=L^{\eta_1}(\lambda){\mathcal{O}}_{r_1}\cdots L^{\eta_k}(\lambda){\mathcal{O}}_{r_k}L^{\eta_{k+1}}(\lambda).\eqno(3.8)$$
 So we have
 $$T_{\bf r}(\lambda)=\sum_{\eta=(\eta_1,\cdots,\eta_{k+1})}T^\eta_{\bf r}(\lambda),\eqno(3.9)$$
 $${\mathcal{F}}_{2r}
 =\frac{1}{2\pi \sqrt{-1}}\sum_{\begin{array}{lcr}
  |{\bf r}|=2r, \eta\\
   r_j\geq1
\end{array}}
\int_{S^1}T^\eta_{\bf r}(\lambda)d\lambda.\eqno(3.10)$$
Since $L^{N^{\bot}}(\lambda)$ is a holomorphic function of $\lambda$, so in (3.10), every non-zero term that appears at least
one $L^{N}(\lambda)$, that is there exists $i_0$ such that $\eta_{i_0}=N$.
Now fix $k$ and $j$ in ${\mathbb{N}}$. Let $s\in L^2({\mathbb{R}}^{2n},{\mathbb{E}}_{x_0})$ be a form of degree $2j$,
${\bf r}\in ({\mathbb{N}}/0)^k$ such that
$\sum_ir_i=2r$ and $\eta=(\eta_1,\cdots,\eta_{k+1})\in \{N,N^{\bot}\}^{k+1}$ such that there is a $i_0$ satisfying $\eta_{i_0}= N$. We
want to find a necessary condition for $I_{2j}T^\eta_{\bf r}(\lambda)I_{2j}$ to be non-zero.\\

\noindent{\bf Case I) all $r_l\geq 2$.}
Since $L^{\eta_{i_0}}=\frac{1}{\lambda}P^N$, and $N$ is concentrated in degree $q$,
we must have that $${\rm deg}({\mathcal{O}}_{r_{i_0}}L^{\eta_{i_0+1}}(\lambda){\mathcal{O}}_{r_{i_0}+1}\cdots L^{\eta_k}(\lambda){\mathcal{O}}_{r_k}L^{\eta_{k+1}}(\lambda)I_{2j}s)$$
has the degree $q$ component. But each $L^{\eta_l}(\lambda)$ preserves the degree and ${\mathcal{O}}_{r_i}$ for $r_i\geq 2$ rises or lowers the degree at most by $2$, so
$$q\geq 2j-2(k-i_0+1),~~ 2j+2(k-i_0+1)\geq q.\eqno(3.11)$$
Similarly, $L^{\eta_{1}}(\lambda){\mathcal{O}}_{r_{1}}\cdots L^{\eta_k}(\lambda){\mathcal{O}}_{r_k}L^{\eta_{k+1}}(\lambda)I_{2j}s$ must have a non-zero component in degree $2j$, then
$$2j\leq 2(i_0-1)+q,~~q-2(i_0-1)\leq 2j.\eqno(3.12)$$
By (3.11) and (3.12), we have
$$4j\leq 2q+2k,~~2q-4j\leq 2k.\eqno(3.13)$$
By $r_l\geq 2$, then $2k\leq 2r$. So
$$|2j-q|\leq r.\eqno(3.14)$$\\

\noindent{\bf Case II) at least one $r_j=1$.} We assume that $l\geq 1$ and there are $l_1$ terms ${\mathcal{O}}_{r_\alpha}={\mathcal{O}}_{1}$ for $k\geq \alpha\geq i_0$ and there are $l-l_1$ terms ${\mathcal{O}}_{r_\beta}={\mathcal{O}}_{1}$ for $1\leq \beta\leq i_0-1$.  Suppose that $I_{2j}T^\eta_{\bf r}(\lambda)I_{2j}\neq 0$. Since $L^{\eta_{i_0}}=\frac{1}{\gamma}
P^N$, $N$ is concentrated in degree $q$, we must have
$${\rm deg}({\mathcal{O}}_{r_{i_0}}L^{\eta_{i_0+1}}(\lambda){\mathcal{O}}_{r_{i_0}+1}\cdots L^{\eta_k}(\lambda){\mathcal{O}}_{r_k}L^{\eta_{k+1}}(\lambda)I_{2j}s)=q\eqno(3.15)$$
We know that $L^{\eta_i}(\lambda)$ and
${\mathcal{O}}_{1}$ preserve the degree, and ${\mathcal{O}}_{r_i}$ for $r_i\geq 2$ rises or lowers the degree at most by $2$, so
$$q\geq 2j-2(k-i_0+1-l_1),~q\leq 2j+2(k-i_0+1-l_1).\eqno(3.16)$$
\indent Similarly, $L^{\eta_1}(\lambda){\mathcal{O}}_{r_1}\cdots L^{\eta_k}(\lambda){\mathcal{O}}_{r_k}L^{\eta_{k+1}}(\lambda)I_{2j}s$
must have a non-zero component in degree $2j$ and ${\mathcal{O}}_{r_i}$ for $r_i\geq 2$ rises or lowers the degree at most by $2$, so we have
$$2j\leq q+2[i_0-1-(l-l_1)],~~2j\geq q-2[i_0-1-(l-l_1)]\eqno(3.17)$$
By (3.16) and (3.17), we get $4j\leq 2q+2k-2l$ and $2q\leq 4j+2k-2l.$
Finally, since $2r=\sum_{i=1}^kr_i\geq 2(k-l)+l=2k-l$, we have $2k-l\leq 2r$. So
$$4j\leq 2q+2r-l,~~2q\leq 4j+2r-l.\eqno(3.18)$$
By $1\leq l$, then $$|2j-q|+\frac{1}{2}\leq |2j-q|+\frac{l}{2} \leq r.\eqno(3.19)$$
\indent Consequently, if $r<|2j-q|$, we have $I_{2j}T^\eta_{\bf r}(\lambda)I_{2j}=0$, and if $r=|2j-q|$, only terms in case I) contribute to $I_{2j}T^\eta_{\bf r}(\lambda)I_{2j}$. \\
\indent For the second part of this theorem, let us assume that we are in the limit case where $r = |2j-q|$. When $2j\geq q$,
by (3.13)
and $2k\leq 2r=4j-2q$, we get $k=2j-q$ and $r=k$. By $r_i\geq 2$, then $r_i=2$ for any $i$. By (3.11) and (3.12), then
$$i_0=\frac{k}{2}+1=j+1-\frac{q}{2},\eqno(3.20)$$
and similar to (2.18) in [PZ], we have
$$I_{2j}{\mathcal{F}}_{2(2j-q)}I_{2j}=I_{2j}({\mathcal{L}}_0^{-1}{\mathcal{O}}_2^{+2})^{j-\frac{q}{2}}P^N
({\mathcal{O}}_2^{-2}{\mathcal{L}}_0^{-1})^{j-\frac{q}{2}}I_{2j}.\eqno(3.21)$$
Let $A=I_{2j}({\mathcal{L}}_0^{-1}{\mathcal{O}}_2^{+2})^{j-\frac{q}{2}}P^N$, then
$$I_{2j}{\mathcal{F}}_{2(2j-q)}I_{2j}=AA^*.\eqno(3.22)$$
By (2.4),(2.21) and (2.25), similar to (2.20) in [PZ], we have
$$A=\frac{1}{(4\pi)^{j-\frac{q}{2}}}B_0^{1,j-\frac{q}{2}}I_{2j}({\mathcal{R}}_{x_0}^\bot)^{j-\frac{q}{2}}P^N.\eqno(3.23)$$
By (2.21) in [PZ] and (3.22),(3.23), we get (1.23). Using the same method, we have (1.24) when $2j<q$.
~~$\Box$

\section{ The second coefficient in the asymptotic expansion}

\quad In this section, we prove Theorem 1.4. Using (3.5), we know that
$$I_{2j}{\bf b}_{2j+1-q}I_{2j}(0, 0) = I_{2j}{\mathcal{F}}_{4j+2-2q}I_{2j }(0, 0).\eqno(4.1)$$
In Section 4.1, we decompose this term into 5 terms, and then in Sections 4.2, 4.3 and 4.4, we
handle them separately.
For the rest of the section we fix an integer $j \in [0, n]$. For every smoothing operator $F$ acting on
$L^2({\mathbb{R}}^{2n}, {\mathbb{E}}_{x_0})$ that appears in this section, we will denote by $F(Z,Z')$ its smooth kernel with
respect to $d\nu_{TX}(Z')$.\\

\noindent{\bf 4.1. Decomposition of the problem.}\\

For $r=2j+1-q$, using the same discussions in Section 3 in [PZ], we see that in $I_{2j}{\mathcal{F}}_{4j+2-2q}I_{2j }(0, 0)$
there are $3$ types of terms $T^{\eta}_{\bf b}(\lambda)$ from case I) with non-zero integral, in which:\\
$~~~~~~~~~~$$\bullet$ for $k=2j-q$:\\
$~~~~~~~~~~~~~$-there are $2j - 2-q$ ${\mathcal{O}}_{r_i}$ equal to ${\mathcal{O}}_{2}$ and $2$ equal to ${\mathcal{O}}_{3}$: we will denote by I the sum
of these terms,\\
$~~~~~~~~~~~~~$-there are $2j - 1-q$ ${\mathcal{O}}_{r_i}$ equal to ${\mathcal{O}}_{2}$ and $1$ equal to ${\mathcal{O}}_{4}$: we will denote by II the sum
of these terms,\\
$~~~~~~~~~~$$\bullet$ for $k = 2j + 1-q$:\\
$~~~~~~~~~~~~~$-all the ${\mathcal{O}}_{r_i}$ are equal to ${\mathcal{O}}_{2}$: we will denote by III the sum of these terms.\\
\indent
For $r=2j+1-q$, by (3.18), we have
$$4j\leq 2q+2k-2l\leq 2q+2r-l=4j+2-l,\eqno(4.2)$$
then $1\leq l\leq 2,$ and $l=1$ or $2$. \\
\indent For $l=1$, then $4j\leq 2q+2k-2\leq 4j+1$, so $k=2j+1-q$. By $4j+2-2q=2r=\sum_{i=1}^kr_i$ and $l=1$, so there are
$2j - 1-q$ ${\mathcal{O}}_{r_i}$ equal to ${\mathcal{O}}_{2}$ and $1$ equal to ${\mathcal{O}}_{3}$ and $1$ equal to ${\mathcal{O}}_{1}$;
we will denote by IV the sum of these terms.\\
\indent For $l=2$, then $4j\leq 2q+2k-4\leq 4j$, so $k=2j+2-q$. By $4j+2-2q=2r=\sum_{s=1}^kr_s$ and $l=2$, so there are
$2j-q$ ${\mathcal{O}}_{r_i}$ equal to ${\mathcal{O}}_{2}$ and $2$ equal to ${\mathcal{O}}_{1}$;
we will denote by V the sum of these terms.\\
\indent We have a decomposition
$$I_{2j}{\mathcal{F}}_{4j+2-2q}I_{2j }(0, 0)={\rm I+II+III+IV+V}.\eqno(4.3)$$
 \indent Similar to Remark 3.1 in [PZ], we know that only ${\mathcal{O}}^{\pm 2}_{2}$, ${\mathcal{O}}^{\pm 2}_{3}$ and ${\mathcal{O}}^{\pm 2}_{4}$
in I and II, and not some ${\mathcal{O}}^{0}_{r_i}$. So by (0.13),(0.14) and (0.15) in [PZ], we get cases I, II and we only need change $j$ to
$j-\frac{q}{2}$ and change ${\mathcal{R}}$ to ${\mathcal{R}}^\bot$. Then we have

$${\rm I=I}_a+{\rm I}^*_a+{\rm I}_b,\eqno(4.4)$$
if ${j-\frac{q}{2}}=0,1$, then ${\rm I}_a=0,$ and if ${j-\frac{q}{2}}\geq 2$,
\begin{eqnarray*}
{\rm I}_a&=&\frac{C_{j-\frac{q}{2}}({j-\frac{q}{2}})}{2\pi}I_{2j}\sum_{\alpha=0}^{{j-\frac{q}{2}}-2}
\sum_{s=0}^\alpha\left\{(C_{j-\frac{q}{2}}({j-\frac{q}{2}})-C_{j-\frac{q}{2}}(\alpha+1))\right.\\
&&\cdot{{\mathcal{R}}^\bot}^{{j-\frac{q}{2}}-(\alpha+2)}_x
(\nabla^{\wedge^{0,\bullet}\otimes E}_{\overline{u_l}}
{{\mathcal{R}}}_{.})^\bot(x)
{{\mathcal{R}}^\bot}^{\alpha-s}_x
(\nabla^{\wedge^{0,\bullet}\otimes E}_{{u_l}}{{\mathcal{R}}})^\bot(x){{\mathcal{R}}^\bot}^{s}_x\\
&&+C_{j-\frac{q}{2}}(s)\left[\prod_{s_1=\alpha+2}^{j-\frac{q}{2}}
(1+\frac{1}{2s_1})-1\right]{{\mathcal{R}}^\bot}^{{j-\frac{q}{2}}-(\alpha+2)}_x
(\nabla^{\wedge^{0,\bullet}\otimes E}_{{u_l}}{{\mathcal{R}}}_{.})^\bot(x)\\
&&\left.{{\mathcal{R}}^\bot}^{\alpha-s}_x
(\nabla^{\wedge^{0,\bullet}\otimes E}_{\overline{u_l}}{{\mathcal{R}}})^\bot(x){{\mathcal{R}}^\bot}^{s}_x\right\}
I_{{\rm det}(\overline{W}^*)\otimes E}({{\mathcal{R}}^\bot}^{j-\frac{q}{2}}_x)^*I_{2j},
~~~~~~~~~~~~~~~~~~~~~~~(4.5)
\end{eqnarray*}
if ${j-\frac{q}{2}}=0$, then ${\rm I}_b=0$, and if ${j-\frac{q}{2}}\geq 1$,
$${\rm I}_b=\frac{1}{2\pi}{\rm I}_{2j}\left[\sum_{k=0}^{{j-\frac{q}{2}}-1}(C_{j-\frac{q}{2}}({j-\frac{q}{2}})-C_{j-\frac{q}{2}}(k)){({\mathcal{R}}_x^\bot)}^{{j-\frac{q}{2}}-k-1}
(\nabla^{\wedge^{0,\bullet}\otimes E}_{\overline{u_l}}{{\mathcal{R}}}_{.})^\bot(x){{\mathcal{R}}_x^\bot}^{k}\right]$$
$$\times I_{{\rm det}(\overline{W}^*)\otimes E}\left[\sum_{k=0}^{{j-\frac{q}{2}}-1}(C_{j-\frac{q}{2}}({j-\frac{q}{2}})-C_{j-\frac{q}{2}}(k)){{\mathcal{R}}_x^\bot}^{{j-\frac{q}{2}}-k-1}
(\nabla^{\wedge^{0,\bullet}\otimes E}_{\overline{u_l}}{{\mathcal{R}}}_{.})^\bot(x){{\mathcal{R}}_x^\bot}^{k}\right]^*I_{2j},\eqno(4.6)$$
$${\rm II}={\rm II}_a+{\rm II}^*_a,\eqno(4.7)$$
if ${j-\frac{q}{2}}=0$, then ${\rm II}_a=0$, and if ${j-\frac{q}{2}}\geq1$,}
$${\rm II}_a=\frac{C_{j-\frac{q}{2}}({j-\frac{q}{2}})}{4\pi}I_{2j}\sum_{k=0}^{{j-\frac{q}{2}}-1}\left\{(C_{j-\frac{q}{2}}({j-\frac{q}{2}})
-C_{j-\frac{q}{2}}(k))\right.$$
$$\left.\cdot({\mathcal{R}}_x^\bot)^{{j-\frac{q}{2}}-(k+1)}_x
(\triangle^{\wedge^{0,\bullet}\otimes E}{{\mathcal{R}}}_{.})^\bot(x)({\mathcal{R}}_x^\bot)^{k}\right\}I_{{\rm det}(\overline{W}^*)\otimes E}(({\mathcal{R}}_x^\bot)^{{j-\frac{q}{2}}})^*I_{2j}.\eqno(4.8)$$\\

\noindent
{\bf 4.2 The term involving only ${\mathcal{O}}_{2}$}\\

\noindent{\bf Lemma 4.1} {\it Any term $T^{\eta}_{\bf r}$
appearing in the term III (with non-vanishing integral) has three types, ${\rm III}_0$: terms with $\eta_{{j-\frac{q}{2}}+1}=\eta_{{j-\frac{q}{2}}+2}=N$ and other $\eta_i=N^{\bot}$,
${\rm III}_a$: terms with $\eta_{{j-\frac{q}{2}}+1}=N$ and other $\eta_i=N^{\bot}$, ${\rm III}_b$: terms with $\eta_{{j-\frac{q}{2}}+2}=N$ and other $\eta_i=N^{\bot}$.
 We have:}
$${\rm III}_0=\sum_{l=1}^{j-\frac{q}{2}}({\mathcal{L}}_0^{-1}{\mathcal{O}}^{+2}_{2})^{l-1} ({\mathcal{L}}_0^{-2}{\mathcal{O}}^{+2}_{2})
({\mathcal{L}}_0^{-1}{\mathcal{O}}^{+2}_{2})^{{j-\frac{q}{2}}-l}P^N{\mathcal{O}}^{0}_{2}P^N({\mathcal{O}}^{-2}_{2}{\mathcal{L}}_0^{-1})^{j-\frac{q}{2}}I_{2j}$$
$$
+\sum_{l=1}^{j-\frac{q}{2}}({\mathcal{L}}_0^{-1}{\mathcal{O}}^{+2}_{2})^{{j-\frac{q}{2}}} P^N{\mathcal{O}}^{0}_{2}P^N
({\mathcal{O}}^{-2}_{2}{\mathcal{L}}_0^{-1})^{l-1}
({\mathcal{O}}^{-2}_{2}{\mathcal{L}}_0^{-2})
({\mathcal{O}}^{-2}_{2}{\mathcal{L}}_0^{-1})^{{j-\frac{q}{2}}-l}
I_{2j},\eqno(4.9)$$
$${\rm III}_a=\sum_{k=0}^{j-\frac{q}{2}}I_{2j}({\mathcal{L}}_0^{-1}{\mathcal{O}}^{+2}_{2})^{{j-\frac{q}{2}}-k}
({\mathcal{L}}_0^{-1}{\mathcal{O}}^{0}_{2})
({\mathcal{L}}_0^{-1}{\mathcal{O}}^{+2}_{2})^{k}P^N
({\mathcal{O}}^{-2}_{2}{\mathcal{L}}_0^{-1})^{{j-\frac{q}{2}}}I_{2j},\eqno(4.10)$$
$${\rm III}_b=({\rm III}_a)^*,~~~{\rm III}={\rm III}_0+{\rm III}_a+{\rm III}_b.\eqno(4.11)$$

\noindent{\bf Proof.} Fix a term $T^{\eta}_{\bf r}$ appearing in the term III with non-vanishing integral. Using again
the same reasoning as in Section 2.2 in [PZ], we see that there exists at most two indices $i_0$ such that
 $\eta_{i_0}=N$, and that they are in $\{{j-\frac{q}{2}}+1, {j-\frac{q}{2}} +2\}$. Now, the only possible term with  $\eta_{{j-\frac{q}{2}}+1}=\eta_{{j-\frac{q}{2}}+2}=N$ is:
$$({\mathcal{L}}_0^{-1}{\mathcal{O}}^{+2}_{2})^{{j-\frac{q}{2}}}P^N{\mathcal{O}}^{0}_{2}P^N({\mathcal{O}}^{-2}_{2}{\mathcal{L}}_0^{-1})^{j-\frac{q}{2}}.$$
By the Cauchy integral formula, we get the term ${\rm III}_0$. Similar to the discussions in Lemma 3.2 in [PZ], we can get ${\rm III}_a$ and ${\rm III}_b$.
~~$\Box$\\

\indent Nextly, we compute $P^N{\mathcal{O}}^{0}_{2}P^N$. We note that it is zero in the K${\rm \ddot{a}}$hler case.
Let
\begin{eqnarray*}
{\mathcal{O}}'_2&=&\frac{1}{3}\left<R^{TX}_{x_0}({\mathcal{R}},e_i){\mathcal{R}},e_s\right>\nabla_{0,e_i}\nabla_{0,e_s}
+\left[\frac{2}{3}\left<R^{TX}_{x_0}({\mathcal{R}},e_s)e_s,e_i\right>\right.\\
&&\left.-\left(\frac{1}{2}\sum_{|\alpha|=2}(\partial_\alpha R^L)_{x_0}\frac{Z^\alpha}{\alpha!}
+R^E_{x_0}\right)({\mathcal{R}},e_i)\right]\nabla_{0,e_i}\\
&&-\frac{1}{4}\nabla_{e_s}\left(\sum_{|\alpha|=2}(\partial_\alpha R^L)_{x_0}
({\mathcal{R}},e_s)\frac{Z^\alpha}{\alpha!}\right)\\
&&-\frac{1}{9}\sum_i\left[\sum_s(\partial_sR^L)_{x_0}({\mathcal{R}},e_i)Z_s\right]^2\\
 &&+\frac{1}{12}\left[\sum_s(\nabla_{0,e_s})^2,\left<R^{TX}_{x_0}({\mathcal{R}},e_l){\mathcal{R}},e_l\right>\right].
 ~~~~~~~~~~~~~~~~~~~~~~~~~~~~~~~~~~~~~~~~~ (4.12)
\end{eqnarray*}
Then by (2.13)

$${\mathcal{O}}^0_2={\mathcal{O}}'_2
-R^{\wedge^{0,\bullet},B}_{x_0}({\mathcal{R}},e_i)\nabla_{0,e_i}
+\left[-(R^E+\frac{1}{2}R^{\rm det})(w_s,\overline{w}_s)
 +\frac{r^X}{4}-\frac{1}{4}c(dT_{as})
 -\frac{1}{8}|T_{as}|^2\right.$$
 $$\left.+2(R^E+\frac{1}{2}R^{\rm det})(w_l,\overline{w}_m)\overline{w}^m\wedge i_{\overline{w}_l}
 \right](x_0)-\frac{\pi}{2}\sqrt{-1}\left<(\nabla^B\nabla^B{\bf J})_{({{\mathcal{R}}},{{\mathcal{R}}})}e_i,e_l\right>c(e_i)c(e_l).
\eqno(4.13)$$

\noindent{\bf Lemma 4.2} (Lemma 3.6.3 in [LuW2]) {\it For the operator ${\mathcal{O}}'_2$, we have}
\begin{eqnarray*}
{\mathcal{O}}'_2{\mathcal{P}}&=&\left\{\frac{1}{3}
b_ib_s\left<R^{TX}_{x_0}({\mathcal{R}},\partial_{\overline{z}_i}){\mathcal{R}},\partial_{\overline{z}_s}\right>
+\frac{1}{2}b_i\sum_{|\alpha|=2}(\partial_\alpha R^L)_{x_0}({\mathcal{R}},\partial_{\overline{z}_i})
\frac{Z^\alpha}{\alpha!}\right.\\
&&+\frac{4}{3}b_s\left[\left<R^{TX}_{x_0}(\partial_{{z_i}},\partial_{\overline{z}_i}){\mathcal{R}},\partial_{\overline{z}_s}\right>
-\left<R^{TX}_{x_0}({\mathcal{R}},\partial_{z_i})\partial_{\overline{z}_i},\partial_{\overline{z}_s}\right>\right]
+R^E({\mathcal{R}},\partial_{\overline{z}_s})b_s\\
&&+\left<(\nabla^X\nabla^X\widetilde{{J}})_{({\mathcal{R}},{\mathcal{R}})}
\partial_{z_s},\partial_{\overline{z}_s}\right>_{x_0}
+4\left<R^{TX}_{x_0}(\partial_{z_i},\partial_{z_s})\partial_{\overline{z}_i},\partial_{\overline{z}_s}\right>\\
&&\left.
-\frac{1}{3}\left[{\mathcal{L}}_0,\left<R^{TX}_{x_0}({\mathcal{R}},\partial_{z_s}){\mathcal{R}},\partial_{\overline{z}_s}\right>\right]
+\frac{1}{9}|(\nabla_{\mathcal{R}}\widetilde{J}){\mathcal{R}}|^2\right\}{\mathcal{P}}.~~~~~~~~~~~~~~~~~~~~~~~(4.14)
\end{eqnarray*}\\

\indent
By (2.25), we have
$$P^N{\mathcal{O}}^0_2P^N=I_{{\rm det}(\overline{W}^*)\otimes E}{\mathcal{P}}{\mathcal{O}}^0_2{\mathcal{P}}I_{{\rm det}(\overline{W}^*)\otimes E},\eqno(4.15)$$
so we only need compute ${\mathcal{P}}{\mathcal{O}}^0_2{\mathcal{P}}$. By [MM12, (4.3)], we have for $\phi\in T^*X$, then
$$\phi(e_i)\nabla_{0,e_i}=\phi(\partial_{z_j})b^+_j-\phi(\partial_{\overline{z}_j})b_j.\eqno(4.16)$$
Recall
$$b^+_l{\mathcal{P}}=0~~{\rm and}~~~(b_l{\mathcal{P}})(Z,Z')=2\pi(\overline{z}_l-\overline{z}'_l){\mathcal{P}}(Z,Z').\eqno(4.17)$$
By Theorem 2.5 and (2.23) and (4.17), we have
$${\mathcal{P}}b^\beta z_\alpha{\mathcal{P}}=0,~~{\mathcal{P}}b^\beta \overline{z}_\alpha{\mathcal{P}}=0.\eqno(4.18)$$
By (4.16),(4.17) and (4.18), we get
$${\mathcal{P}}R^{\wedge^{0,\bullet},B}_{x_0}
({\mathcal{R}},e_i)\nabla_{0,e_i}{\mathcal{P}}=
-2R^{\wedge^{0,\bullet},B}_{x_0}(\partial_{z_s},\partial_{\overline{z}_s}){\mathcal{P}}.\eqno(4.19)$$
By ${\mathcal{P}}^2={\mathcal{P}}$, then
$${\mathcal{P}}\left[-(R^E+\frac{1}{2}R^{\rm det})(w_s,\overline{w}_s)
 +\frac{r^X}{4}
 -\frac{1}{8}|T_{as}|^2\right](x_0){\mathcal{P}}$$
 $$=\left[-(R^E+\frac{1}{2}R^{\rm det})(w_s,\overline{w}_s)
 +\frac{r^X}{4}-
 -\frac{1}{8}|T_{as}|^2\right](x_0){\mathcal{P}}.\eqno(4.20)$$
By (2.11) and direct computations, then
$$I_{{\rm det}(\overline{W}^*)\otimes E}{\mathcal{P}}\left[2(R^E+\frac{1}{2}R^{\rm det})(w_l,\overline{w}_m)\overline{w}^m\wedge i_{\overline{w}_l}
 \right](x_0){\mathcal{P}}I_{{\rm det}(\overline{W}^*)\otimes E}$$
 $$
=2\sum_{m=1}^q(R^E+\frac{1}{2}R^{\rm det})(w_m,\overline{w}_m)
 (x_0){\mathcal{P}}I_{{\rm det}(\overline{W}^*)\otimes E};\eqno(4.21)$$

$$I_{{\rm det}(\overline{W}^*)\otimes E}{\mathcal{P}}\left[-\frac{1}{4}c(dT_{as})
 \right](x_0){\mathcal{P}}I_{{\rm det}(\overline{W}^*)\otimes E}=[-\frac{1}{8}dT_{as}(w_i,\overline{w}_i,w_s,\overline{w}_s)$$
$$+\sum_{i=1}^qdT_{as}(w_i,\overline{w}_i,w_s,\overline{w}_s)
-2\sum_{i,l=1}^qdT_{as}(w_i,w_l,\overline{w}_i,\overline{w}_l)]{\mathcal{P}}I_{{\rm det}(\overline{W}^*)\otimes E};\eqno(4.22)$$
Similar to (4.19), we have
$${\mathcal{P}}R^E_{x_0}({\mathcal{R}},\partial_{\overline{z}_i})b_i{\mathcal{P}}=2R^E_{x_0}
(\partial_{{z}_i},\partial_{\overline{z}_i}){\mathcal{P}}.\eqno(4.23)$$
By (2.22),(2.23) and (2.24), we have
$$({\mathcal{P}}z^\alpha\overline{z}^\beta{\mathcal{P}})(0,0)=\frac{\alpha!}{\pi^{|\alpha|}}\delta_{\alpha\beta}.\eqno(4.24)$$
By (4.24), then
$${\mathcal{P}}\left[\left<(\nabla^X\nabla^X\widetilde{{ J}})_{({\mathcal{R}},{\mathcal{R}})}
\partial_{z_i},\partial_{\overline{z}_i}\right>_{x_0}
\right]{\mathcal{P}}(0,0)~~~~~~~~~~~~~~~~$$
$$=\frac{1}{4\pi}\left[\left<(\nabla^X\nabla^X\widetilde{{ J}})_{(u_s,\overline{u}_s)}u_i,\overline{u_i}\right>
+\left<(\nabla^X\nabla^X\widetilde{{ J}})_{(\overline{u}_s,{u_s})}{u_i},\overline{u_i}\right>
\right].\eqno(4.25)$$
Similar to (4.20), then
$$4\left({\mathcal{P}}\left<R^{TX}_{x_0}(\partial_{z_i},\partial_{z_s})\partial_{\overline{z}_i},\partial_{\overline{z}_s}\right>
{\mathcal{P}}\right)(0,0)=\left<R^{TX}_{x_0}({u_i},{u_s}){\overline{u}_i},{\overline{u}_s}\right>.\eqno(4.26)$$
By the relation $[AB,C]=[A,C]B+A[B,C]$ and (4.18), it holds that
$$-\frac{1}{3}\left({\mathcal{P}}\left[{\mathcal{L}}_0,\left<R^{TX}_{x_0}({\mathcal{R}},\partial_{z_s}){\mathcal{R}},\partial_{\overline{z}_s}\right>\right]
{\mathcal{P}}\right)(0,0)=0.\eqno(4.27)$$
By (4.24),(1.32) and (1.35) and $\Psi(X_0,Y,U,V)=\Psi(U,V,X_0,Y)$, we have
\begin{eqnarray*}
&&\frac{1}{9}\left({\mathcal{P}}|(\nabla_{\mathcal{R}}\widetilde{J}){\mathcal{R}}|^2{\mathcal{P}}\right)(0,0)
\\
&=&-\frac{1}{9}\left({\mathcal{P}}\left<(\nabla_{\mathcal{R}}\widetilde{J}){\mathcal{R}},(\nabla_{\mathcal{R}}\widetilde{J}){\mathcal{R}}\right>
{\mathcal{P}}\right)(0,0)
\\
&=&
\frac{1}{8}\sum_{s\neq l=1}^n\left({\mathcal{P}}z_sz_l\overline{z}_s\overline{z}_l{\mathcal{P}}\right)(0,0)
[\Psi_{sl\overline{s}\overline{l}}+\Psi_{sl\overline{l}\overline{s}}+\Psi_{ls\overline{s}\overline{l}}+\Psi_{ls\overline{l}\overline{s}}+
\Psi_{l\overline{l}s\overline{s}}+\Psi_{l\overline{l}\overline{s}s}\\
&&+\Psi_{\overline{l}ls\overline{s}}+
\Psi_{\overline{l}l\overline{s}s}+\Psi_{s\overline{l}\overline{s}l}+\Psi_{s\overline{l}l\overline{s}}+
\Psi_{\overline{l}s\overline{s}l}+\Psi_{\overline{l}sl\overline{s}}]\\
&&+\frac{1}{4}\sum_{s=1}^n\left({\mathcal{P}}z^2_s\overline{z}^2_s{\mathcal{P}}\right)(0,0)
[2\Psi_{ss\overline{s}\overline{s}}+2\Psi_{s\overline{s}\overline{s}s}+\Psi_{s\overline{s}s\overline{s}}+
\Psi_{\overline{s}s\overline{s}s}]\\
&=&\frac{1}{4\pi^2}(2\Psi_{sl\overline{s}\overline{l}}+2\Psi_{sl\overline{l}\overline{s}}+2\Psi_{s\overline{l}\overline{s}l}
+2\Psi_{\overline{l}ls\overline{s}}+\Psi_{l\overline{l}s\overline{s}}+\Psi_{\overline{l}l\overline{s}s}
+\Psi_{s\overline{l}l\overline{s}}+\Psi_{\overline{l}s\overline{s}l}),~~~~~~~(4.28)
\end{eqnarray*}
By (4.18), it holds that
$$\frac{4}{3}{\mathcal{P}}b_s\left[\left<R^{TX}_{x_0}(\partial_{{z_i}},\partial_{\overline{z}_i}){\mathcal{R}},\partial_{\overline{z}_s}\right>
-\left<R^{TX}_{x_0}({\mathcal{R}},\partial_{z_i})\partial_{\overline{z}_i},\partial_{\overline{z}_s}\right>\right]{\mathcal{P}}=0,
\eqno(4.29)$$
Similar to (4.18), we have
$${\mathcal{P}}b_lz_sz_{j'}\overline{z}_{j''}{\mathcal{P}}={\mathcal{P}}b_lz_s\overline{z}_{j'}\overline{z}_{j''}{\mathcal{P}}=0.\eqno(4.30)$$
So by (4.30) and $R^L$ is $(1,1)$-form and $Z^\alpha$ for $|\alpha|=2$ being the combinators of $z_sz_{j'},~z_s\overline{z}_{j'},~\overline{z}_{s}\overline{z}_{j'}$, then
$$\frac{1}{2}{\mathcal{P}}b_i\sum_{|\alpha|=2}(\partial_\alpha R^L)_{x_0}({\mathcal{R}},\partial_{\overline{z}_i})
\frac{Z^\alpha}{\alpha!}{\mathcal{P}}=0.\eqno(4.31)$$
Similar to (4.31), we can get
$$\frac{1}{3}{\mathcal{P}}
b_ib_s\left<R^{TX}_{x_0}({\mathcal{R}},\partial_{\overline{z}_i}){\mathcal{R}},\partial_{\overline{z}_s}\right>{\mathcal{P}}=0.\eqno(4.32)$$
By (2.3), similar to (4.21) and (4.25), we get
$$-\frac{\pi}{2}\sqrt{-1}I_{{\rm det}(\overline{W}^*)\otimes E}[{\mathcal{P}}\left<(\nabla^B\nabla^B{\bf J})_{({{\mathcal{R}}},{{\mathcal{R}}})}e_i,e_l\right>c(e_i)c(e_l){\mathcal{P}}](0,0)I_{{\rm det}(\overline{W}^*)\otimes E}$$
$$=-\sqrt{-1}\left[-\frac{1}{2}\left<(\nabla^B\nabla^B{\bf J})_{(u_i,\overline{u}_i)}w_s,\overline{w}_s\right>
-\frac{1}{2}\left<(\nabla^B\nabla^B{\bf J})_{(\overline{u}_i,u_i)}w_s,\overline{w}_s\right>\right.$$
$$\left.
+\sum_{m=1}^q\left<(\nabla^B\nabla^B{\bf J})_{(u_i,\overline{u}_i)}w_m,\overline{w}_m\right>
+\sum_{m=1}^q\left<(\nabla^B\nabla^B{\bf J})_{(\overline{u}_l,u_l)}w_m,\overline{w}_m\right>\right]
I_{{\rm det}(\overline{W}^*)\otimes E}.\eqno(4.33)$$

By Theorem 2.4, (4.12)-(4.15),(4.19)-(4.23),(4.25)-(4.27),(4.31)-(4.33), we get\\

\noindent{\bf Lemma 4.3} {\it The following equation holds}
\begin{eqnarray*}
&&({{P}^N}{\mathcal{O}}^0_2{{P}^N})(0,0)~~~~~~~~~~~~~~~~~~~~~~~~~~~~~~~~~~~~~~~~~~~~~~~~~~~~~~~~~~~~~~~~~~~~~~~~~~(4.34)\\
&=&I_{{\rm det}(\overline{W}^*)\otimes E}\left\{
R^{\wedge^{0,\bullet},B}(u_s,\overline{u_s})-(R^E+\frac{1}{2}R^{\rm det})(w_s,\overline{w_s})
+\frac{r^X}{4}-\frac{1}{8}|T_{as}|^2\right.\\
&&
+2\sum_{m=1}^q(R^E+\frac{1}{2}R^{\rm det})(w_m,\overline{w}_m)
-\frac{1}{8}dT_{as}(w_i,\overline{w}_i,w_s,\overline{w}_s)\\
&&+\sum_{i=1}^qdT_{as}(w_i,\overline{w}_i,w_s,\overline{w}_s)
-2\sum_{i,l=1}^qdT_{as}(w_i,w_l,\overline{w}_i,\overline{w}_l)+R^E(u_s,\overline{u_s})\\
&&
+\frac{1}{4\pi}\left[\left<(\nabla^X\nabla^X\widetilde{{ J}})_{(u_s,\overline{u}_s)}u_i,\overline{u}_i\right>
+\left<(\nabla^X\nabla^X\widetilde{{J}})_{(\overline{u}_s,{u_s})}{u_i},\overline{u}_i\right>\right]\\
&&+
\left<R^{TX}_{x_0}({u_i},{u_s}){\overline{u}_i},{\overline{u}_s}\right>
+\frac{1}{4\pi^2}(2\Psi_{sl\overline{s}\overline{l}}+2\Psi_{sl\overline{l}\overline{s}}+2\Psi_{s\overline{l}\overline{s}l}\\
&&
+2\Psi_{\overline{l}ls\overline{s}}+\Psi_{l\overline{l}s\overline{s}}+\Psi_{\overline{l}l\overline{s}s}
+\Psi_{s\overline{l}l\overline{s}}+\Psi_{\overline{l}s\overline{s}l})\\
&&
-\sqrt{-1}\left[-\frac{1}{2}\left<(\nabla^B\nabla^B{\bf J})_{(u_i,\overline{u}_i)}w_s,\overline{w}_s\right>
-\frac{1}{2}\left<(\nabla^B\nabla^B{\bf J})_{(\overline{u}_i,u_i)}w_s,\overline{w}_s\right>\right.\\
&&\left.\left.
+\sum_{m=1}^q\left<(\nabla^B\nabla^B{\bf J})_{(u_i,\overline{u}_i)}w_m,\overline{w}_m\right>
+\sum_{m=1}^q\left<(\nabla^B\nabla^B{\bf J})_{(\overline{u}_l,u_l)}w_m,\overline{w}_m\right>\right]\right\}(x_0)
I_{{\rm det}(\overline{W}^*)\otimes E}.
\end{eqnarray*}

\indent By the definition of ${\mathcal{L}}_0$, we have
$$({\mathcal{L}}_0^{-1}{\mathcal{O}}^{+2}_{2})^{l-1} ({\mathcal{L}}_0^{-2}{\mathcal{O}}^{+2}_{2})
({\mathcal{L}}_0^{-1}{\mathcal{O}}^{+2}_{2})^{{j-\frac{q}{2}}-l}{\mathcal{P}}$$
$$=\frac{1}{(4\pi)^{{j-\frac{q}{2}}+1}
2^{{j-\frac{q}{2}}+1}(j-\frac{q}{2})!({j-\frac{q}{2}}-l+1)}({\mathcal{R}}_{x_0}^\bot)^{j-\frac{q}{2}}
{\mathcal{P}},\eqno(4.35)$$
and
$${\mathcal{P}}({\mathcal{O}}^{-2}_{2}{\mathcal{L}}_0^{-1})^{j-\frac{q}{2}}I_{2j}=
[({\mathcal{L}}_0^{-1}{\mathcal{O}}^{+2}_{2})^{{j-\frac{q}{2}}}{\mathcal{P}}]^*
=\frac{1}{(4\pi)^{j-\frac{q}{2}}2^{j-\frac{q}{2}}({j-\frac{q}{2}})!}{\mathcal{P}}({\mathcal{R}}_{x_0}^{\bot,*})^{j-\frac{q}{2}}.\eqno(4.36)$$
By (4.9),(4.35) and (4.36), we have
$${\rm III}_0(0,0)=\left(\sum_{l=1}^{j-\frac{q}{2}}\frac{1}{l}\right)\frac{1}{(4\pi)^{2j-q+1}2^{2j-q}(({j-\frac{q}{2}})!)^2}
({\mathcal{R}}_{x_0}^\bot)^{j-\frac{q}{2}}({{P}^N}{\mathcal{O}}^0_2{{P}^N})(0,0)({\mathcal{R}}_{x_0}^{\bot,*})^{j-\frac{q}{2}}I_{2j}.\eqno(4.37)$$
By (4.34) and (4.37), we get ${\rm III}_0(0,0)$ .\\

\indent Nextly {\bf we compute ${\rm III}_a$}. We know that $\wedge^{(0,*)}(TX)=\wedge^*(\overline{W}^*)\otimes \wedge ^*({\overline{W}}^{\bot,*})$ has double degree.
Let $\widetilde{R}^{\wedge^{0,\bullet},B}$ and $\frac{1}{4}\widetilde{dT_{as}}$ be the preserving double degree parts
of ${R}^{\wedge^{0,\bullet},B}$ and $\frac{1}{4}{dT_{as}}$ respectively. Then
$${\mathcal{O}}^0_2={\mathcal{O}}^{0,1}_2+{\mathcal{O}}^{0,2}_2,\eqno(4.38)$$
where\begin{eqnarray*}
{\mathcal{O}}^{0,1}_2&=&
\frac{1}{3}\left<R^{TX}_{x_0}({\mathcal{R}},e_i){\mathcal{R}},e_s\right>\nabla_{0,e_i}\nabla_{0,e_s}
+\left[\frac{2}{3}\left<R^{TX}_{x_0}({\mathcal{R}},e_s)e_s,e_i\right>\right.\\
&&\left.-\left(\frac{1}{2}\sum_{|\alpha|=2}(\partial_\alpha R^L)_{x_0}\frac{Z^\alpha}{\alpha!}
+R^E_{x_0}+\widetilde{R}^{\wedge^{0,\bullet},B}_{x_0}\right)({\mathcal{R}},e_i)\right]\nabla_{0,e_i}\\
&&-\frac{1}{4}\nabla_{e_s}\left(\sum_{|\alpha|=2}(\partial_\alpha R^L)_{x_0}
({\mathcal{R}},e_s)\frac{Z^\alpha}{\alpha!}\right)\\
&&-\frac{1}{9}\sum_i\left[\sum_s(\partial_sR^L)_{x_0}({\mathcal{R}},e_i)Z_s\right]^2\\
&&+\left[-(R^E+\frac{1}{2}R^{\rm det})(w_l,\overline{w}_l)
 +\frac{r^X}{4}-\frac{1}{4}\widetilde{c(dT_{as})}\right.\\
 &&\left.-\frac{1}{8}|T_{as}|^2
 +2(\sum_{l,m\leq q}+\sum_{l,m\geq q+1})(R^E+\frac{1}{2}R^{\rm det})(w_l,\overline{w}_m)\overline{w}^m\wedge i_{\overline{w}_l}
 \right]\\
 &&+\frac{1}{12}\left[\sum_s(\nabla_{0,e_s})^2,\left<R^{TX}_{x_0}({\mathcal{R}},e_l){\mathcal{R}},e_l\right>\right]\\
 &&-\frac{\pi}{2}\sqrt{-1}\left[-\frac{1}{2}\left<(\nabla^B\nabla^B{\bf J})_{({{\mathcal{R}}},{{\mathcal{R}}})}w_l,\overline{w}_l\right>\right.\\
 &&\left.+
 (\sum_{l,m\leq q}+\sum_{l,m\geq q+1})\left<(\nabla^B\nabla^B{\bf J})_{({{\mathcal{R}}},{{\mathcal{R}}})}w_l,\overline{w}_m\right>\overline{w}^m
 \wedge i_{\overline{w}_l}\right];~~~~~~~~~~~~~~~~~~~~~(4.39)
\end{eqnarray*}

\begin{eqnarray*}
{\mathcal{O}}^{0,2}_2&=&
\left({\widetilde{R}}^{\wedge^{0,\bullet},B}_{x_0}-{R}^{\wedge^{0,\bullet},B}_{x_0}\right)({\mathcal{R}},e_i)\nabla_{0,e_i}
+\frac{1}{4}(\widetilde{c(dT_{as})}-c(dT_{as}))\\
&&+2(\sum_{l\leq q,m\geq q+1}+\sum_{l\geq q+1,m\leq q})(R^E+\frac{1}{2}R^{\rm det})(w_l,\overline{w}_m)\overline{w}^m\wedge i_{\overline{w}_l}\\
  &&-\frac{\pi}{2}\sqrt{-1}(\sum_{l\leq q,m\geq q+1}+\sum_{l\geq q+1,m\leq q})
 \left<(\nabla^B\nabla^B{\bf J})_{({{\mathcal{R}}},{{\mathcal{R}}})}w_l,\overline{w}_m\right>\overline{w}^m
 \wedge i_{\overline{w}_l}.~~~~~~(4.40)
\end{eqnarray*}
Then ${\mathcal{O}}^{0,1}_2$ is the preserving double degree part of ${\mathcal{O}}^{0}_2$ and
${\mathcal{O}}^{0,2}_2$ is the changing double degree part of ${\mathcal{O}}^{0}_2$.\\
\indent By (4.12), (4.39) and (4.16), we have

\begin{eqnarray*}
&&I_{2j}({\mathcal{L}}_0^{-1}{\mathcal{O}}^{+2}_{2})^{{j-\frac{q}{2}}-k}
({\mathcal{L}}_0^{-1}{\mathcal{O}}^{0,1}_{2})
({\mathcal{L}}_0^{-1}{\mathcal{O}}^{+2}_{2})^{k}P^N
\\
&=&\frac{1}{(4\pi)^k2^kk!}I_{2j}({\mathcal{L}}_0^{-1}{\mathcal{O}}^{+2}_{2})^{{j-\frac{q}{2}}-k}
({\mathcal{L}}_0^{-1}{\mathcal{O}}^{0}_{2})
({\mathcal{R}}^\bot)^k_{x_0}{\mathcal{P}}
\\
&=&\frac{1}{(4\pi)^k2^kk!}I_{2j}({\mathcal{L}}_0^{-1}{\mathcal{O}}^{+2}_{2})^{{j-\frac{q}{2}}-k}
{\mathcal{L}}_0^{-1}\left\{R^E({\mathcal{R}},\partial_{\overline{z}_s})b_s+\widetilde{R}^{\wedge^{0,\bullet},B}_{x_0}
({\mathcal{R}},\partial_{\overline{z}_s})b_s\right.\\
&&-(R^E+\frac{1}{2}R^{\rm det})(w_s,\overline{w}_s)
 +\frac{r^X}{4}-\frac{1}{4}\widetilde{c(dT_{as})}
 -\frac{1}{8}|T_{as}|^2\\
&& +2(\sum_{l,m\leq q}+\sum_{l,m\geq q+1})
(R^E+\frac{1}{2}R^{\rm det})(w_l,\overline{w}_m)\overline{w}^m\wedge i_{\overline{w}_l}\\
&&-\frac{\pi}{2}\sqrt{-1}\left[-\frac{1}{2}\left<(\nabla^B\nabla^B{\bf J})_{({{\mathcal{R}}},{{\mathcal{R}}})}w_l,\overline{w}_l\right>\right.\\
 &&\left.\left.+
 (\sum_{l,m\leq q}+\sum_{l,m\geq q+1})\left<(\nabla^B\nabla^B{\bf J})_{({{\mathcal{R}}},{{\mathcal{R}}})}w_l,\overline{w}_m\right>\overline{w}^m
 \wedge i_{\overline{w}_l}\right]\right\}
 (x_0)({\mathcal{R}}^\bot)^k_{x_0}{{P}^N}\\
&&+\frac{1}{(4\pi)^k2^kk!}I_{2j}({\mathcal{L}}_0^{-1}{\mathcal{O}}^{+2}_{2})^{{j-\frac{q}{2}}-k}
{\mathcal{L}}_0^{-1}({\mathcal{R}}^\bot)^k_{x_0}[{\mathcal{O}}'_2+R^E_{x_0}({\mathcal{R}},e_s)\nabla_{0,e_s}]{{P}^N}
,~~~~(4.41)
\end{eqnarray*}
By (2.20) and (2.21) and (4.17), we have
$$({\mathcal{L}}_0^{-1}{\mathcal{O}}^{+2}_{2})^{{j-\frac{q}{2}}-k}
{\mathcal{L}}_0^{-1}(R^E+\widetilde{R}^{\wedge^{0,\bullet},B})({\mathcal{R}},\partial_{\overline{z}_s})b_s({\mathcal{R}}^\bot)^k_{x_0}{\mathcal{P}}
=({\mathcal{L}}_0^{-1}{\mathcal{O}}^{+2}_{2})^{{j-\frac{q}{2}}-k}
{\mathcal{L}}_0^{-1}$$
$$\cdot\left[(R^E+\widetilde{R}^{\wedge^{0,\bullet},B})(\partial_{z_l},\partial_{\overline{z}_s})({\mathcal{R}}^\bot)^k_{x_0}b_sz_l{\mathcal{P}}
+2(R^E+\widetilde{R}^{\wedge^{0,\bullet},B})(\partial_{z_s},\partial_{\overline{z}_s})({\mathcal{R}}^\bot)^k_{x_0}{\mathcal{P}}\right.$$
$$\left.+\frac{1}{2\pi}(R^E+\widetilde{R}^{\wedge^{0,\bullet},B})(\partial_{\overline{z}_l},\partial_{\overline{z}_s})({\mathcal{R}}^\bot)^k_{x_0}
b_sb_l{\mathcal{P}}
+(R^E+\widetilde{R}^{\wedge^{0,\bullet},B})(\partial_{\overline{z}_l},\partial_{\overline{z}_s})
({\mathcal{R}}^\bot)^k_{x_0}b_s\overline{z}'_l{\mathcal{P}}\right].\eqno(4.42)$$
By (4.17), direct computations show that

$$\left\{({\mathcal{L}}_0^{-1}{\mathcal{O}}^{+2}_{2})^{{j-\frac{q}{2}}-k}
{\mathcal{L}}_0^{-1}
(R^E+\widetilde{R}^{\wedge^{0,\bullet},B})(\partial_{z_l},\partial_{\overline{z}_s})({\mathcal{R}}^\bot)^k_{x_0}b_sz_l{\mathcal{P}}\right\}(0,Z)$$
$$=\frac{-2}{(4\pi)^{{j-\frac{q}{2}}-k+1}}B^{k,{j-\frac{q}{2}}}_1({\mathcal{R}}^\bot)^{{j-\frac{q}{2}}-k}_{x_0}(R^E+\widetilde{R}^{\wedge^{0,\bullet},B})(\partial_{{z}_s},\partial_{\overline{z}_s})
({\mathcal{R}}^\bot)^k_{x_0}{\mathcal{P}}(0,Z);\eqno(4.43)$$

$$2\left\{({\mathcal{L}}_0^{-1}{\mathcal{O}}^{+2}_{2})^{{j-\frac{q}{2}}-k}
{\mathcal{L}}_0^{-1}(R^E+\widetilde{R}^{\wedge^{0,\bullet},B})(\partial_{z_s},\partial_{\overline{z}_s})
({\mathcal{R}}^\bot)^k_{x_0}{\mathcal{P}}\right\}(0,Z)$$
$$=\frac{2}{(4\pi)^{{j-\frac{q}{2}}-k+1}}B^{k,{j-\frac{q}{2}}}_0({\mathcal{R}}^\bot)^{{j-\frac{q}{2}}-k}_{x_0}(R^E+\widetilde{R}^{\wedge^{0,\bullet},B})(\partial_{{z}_s},\partial_{\overline{z}_s})
({\mathcal{R}}^\bot)^k_{x_0}{\mathcal{P}}(0,Z);\eqno(4.44)$$

$$\frac{1}{2\pi}\left\{({\mathcal{L}}_0^{-1}{\mathcal{O}}^{+2}_{2})^{{j-\frac{q}{2}}-k}
{\mathcal{L}}_0^{-1}
(R^E+\widetilde{R}^{\wedge^{0,\bullet},B})(\partial_{\overline{z}_l},\partial_{\overline{z}_s})({\mathcal{R}}^\bot)^k_{x_0}b_sb_l{\mathcal{P}}\right\}(0,Z)$$
$$=\frac{2\pi}{(4\pi)^{{j-\frac{q}{2}}-k+1}}B^{k,{j-\frac{q}{2}}}_2({\mathcal{R}}^\bot)^{{j-\frac{q}{2}}-k}_{x_0}(R^E+\widetilde{R}^{\wedge^{0,\bullet},B})(\partial_{\overline{z}_l},\partial_{\overline{z}_s})
({\mathcal{R}}^\bot)^k_{x_0}\overline{z}_s\overline{z}_l({\mathcal{P}}(0,Z));\eqno(4.45)$$

$$\left\{({\mathcal{L}}_0^{-1}{\mathcal{O}}^{+2}_{2})^{{j-\frac{q}{2}}-k}
{\mathcal{L}}_0^{-1}(R^E+\widetilde{R}^{\wedge^{0,\bullet},B})(\partial_{\overline{z}_l},\partial_{\overline{z}_s})({\mathcal{R}}^\bot)^k_{x_0}b_s\overline{z}'_l{\mathcal{P}}
\right\}(0,Z)$$
$$=\frac{-2\pi}{(4\pi)^{{j-\frac{q}{2}}-k+1}}B^{k,{j-\frac{q}{2}}}_1({\mathcal{R}}^\bot)^{{j-\frac{q}{2}}-k}_{x_0}
(R^E+\widetilde{R}^{\wedge^{0,\bullet},B})(\partial_{\overline{z}_l},\partial_{\overline{z}_s})
({\mathcal{R}}^\bot)^k_{x_0}\overline{z}_s\overline{z}_l({\mathcal{P}}(0,Z)).\eqno(4.46)$$
Similarly
$$I_{2j}({\mathcal{L}}_0^{-1}{\mathcal{O}}^{+2}_{2})^{{j-\frac{q}{2}}-k}
{\mathcal{L}}_0^{-1}\left[
-(R^E+\frac{1}{2}R^{\rm det})(w_s,\overline{w}_s)
 +\frac{r^X}{4}-\frac{1}{4}\widetilde{c(dT_{as})}\right.$$
 $$\left. -\frac{1}{8}|T_{as}|^2
+2(\sum_{l,m\leq q}+\sum_{l,m\geq q+1})(R^E+\frac{1}{2}R^{\rm det})(w_l,\overline{w}_m)\overline{w}^m\wedge i_{\overline{w}_l}
 \right](x_0)({\mathcal{R}}^\bot)^k_{x_0}{\mathcal{P}}$$
 $$=\frac{1}{(4\pi)^{{j-\frac{q}{2}}-k+1}}B^{k,{j-\frac{q}{2}}}_0({\mathcal{R}}^\bot)^{{j-\frac{q}{2}}-k}_{x_0}
\left[
-(R^E+\frac{1}{2}R^{\rm det})(w_s,\overline{w}_s)
 +\frac{r^X}{4}-\frac{1}{4}\widetilde{c(dT_{as})}\right.$$
 $$\left. -\frac{1}{8}|T_{as}|^2
+2(\sum_{l,m\leq q}+\sum_{l,m\geq q+1})(R^E+\frac{1}{2}R^{\rm det})(w_l,\overline{w}_m)\overline{w}^m\wedge i_{\overline{w}_l}
 \right](x_0)({\mathcal{R}}^\bot)^k_{x_0}{\mathcal{P}}.\eqno(4.47)$$
Recall [MM08, (2.7)], by (1.33) and $R^L$ is a $(1,1)$-form, we have
$$b_i\sum_{|\alpha|=2}(\partial_\alpha R^L)_{x_0}({\mathcal{R}},\partial_{\overline{z}_i})
\frac{Z^\alpha}{\alpha!}=b_iz_l\widehat{\Psi}({\mathcal{R}},{\mathcal{R}},\partial_{{z}_l},\partial_{\overline{z}_i}).\eqno(4.48)$$
By (2.20) and (4.24), we have
$$(b_lz_\alpha z_\beta z_\gamma{\mathcal{P}})(0,Z)=0; ~~(b_lb_\gamma z_\alpha z_\beta{\mathcal{P}})(0,Z)=4(\delta_{l\alpha}\delta_{\gamma\beta}+\delta_{l\beta}\delta{\gamma\alpha}){\mathcal{P}}(0,Z);\eqno(4.49)$$
$$(b_lz_\alpha z_\beta \overline{z}_\gamma{\mathcal{P}})(0,Z)=\frac{1}{2\pi}[(b_lb_\gamma z_\alpha z_\beta+2\delta_{\gamma\beta}b_jz_\alpha
+2\delta_{\gamma\alpha}b_lz_\beta){\mathcal{P}}](0,Z);\eqno(4.50)$$
$$\int_{{\mathbb{R}}^{2n}}(b_lz_\alpha \overline{z}_\beta \overline{z}_\gamma{\mathcal{P}})(0,Z){\mathcal{P}}(Z,0)d\nu_{TX}(Z)=0;\eqno(4.51)$$
$$\int_{{\mathbb{R}}^{2n}}(b_l\overline{z}_\alpha \overline{z}_\beta \overline{z}_\gamma{\mathcal{P}})(0,Z){\mathcal{P}}(Z,0)d\nu_{TX}(Z)=0.\eqno(4.52)$$
By (4.48)-(4.52), then
$$\left(b_i\sum_{|\alpha|=2}(\partial_\alpha R^L)_{x_0}({\mathcal{R}},\partial_{\overline{z}_i})
\frac{Z^\alpha}{\alpha!}{\mathcal{P}}\right)(0,Z)$$
$$
=[\widehat{\Psi}(\partial_{z_\beta},\partial_{\overline{z}_\gamma},\partial_{{z}_\alpha},\partial_{\overline{z}_i})+
\widehat{\Psi}(\partial_{\overline{z}_\gamma},\partial_{z_\beta},\partial_{{z}_\alpha},\partial_{\overline{z}_i})](b_i
z_\alpha z_\beta \overline{z}_\gamma{\mathcal{P}})(0,Z).\eqno(4.53)$$
By (4.49),(4.50) and (4.53), then
$$\left\{({\mathcal{L}}_0^{-1}{\mathcal{O}}^{+2}_{2})^{{j-\frac{q}{2}}-k}
{\mathcal{L}}_0^{-1}({\mathcal{R}}^\bot)^k_{x_0}\left[\frac{1}{2}b_i\sum_{|\alpha|=2}(\partial_\alpha R^L)_{x_0}({\mathcal{R}},\partial_{\overline{z}_i})
\frac{Z^\alpha}{\alpha!}\right]{\mathcal{P}}\right\}(0,Z)$$
$$=\frac{1}{(4\pi)^{{j-\frac{q}{2}}-k+2}}(B^{k,{j-\frac{q}{2}}}_2-B^{k,{j-\frac{q}{2}}}_1)
({\mathcal{R}}^\bot)^{j-\frac{q}{2}}_{x_0}(\widehat{\Psi}_{\beta\overline{\beta}\alpha\overline{\alpha}}
+\widehat{\Psi}_{\overline{\alpha}\beta\alpha\overline{\beta}}).\eqno(4.54)$$
It holds that
$$
b_ib_l\left<R^{TX}_{x_0}({\mathcal{R}},\partial_{\overline{z}_i}){\mathcal{R}},\partial_{\overline{z}_l}\right>=
\frac{1}{4}(R^X_{\alpha\overline{i}\beta\overline{l}}b_ib_lz_\alpha z_\beta+R^X_{\overline{\alpha}\overline{i}\overline{\beta}
\overline{l}}b_ib_l\overline{z}_\alpha \overline{ z}_\beta+
2R^X_{\overline{\alpha}\overline{i}{\beta}
\overline{l}}b_ib_l\overline{z}_\alpha { z}_\beta).\eqno(4.55)$$
By (4.49), then
$$\frac{1}{12}\left[({\mathcal{L}}_0^{-1}{\mathcal{O}}^{+2}_{2})^{{j-\frac{q}{2}}-k}
{\mathcal{L}}_0^{-1}({\mathcal{R}}^\bot)^k_{x_0}R^X_{\alpha\overline{i}\beta\overline{l}}b_ib_lz_\alpha z_\beta{\mathcal{P}}\right](0,Z)$$
$$=\frac{1}{3}\frac{1}{(4\pi)^{{j-\frac{q}{2}}-k+1}}B^{{j-\frac{q}{2}},k}_2({\mathcal{R}}^\bot)^{j-\frac{q}{2}}_{x_0}
(R^X_{\alpha\overline{\alpha}\beta\overline{\beta}}
+R^X_{\alpha\overline{\beta}\beta\overline{\alpha}}){\mathcal{P}}(0,Z).\eqno(4.56)$$
By (4.17), then
$$(b_ib_l\overline{z}_\alpha \overline{ z}_\beta{\mathcal{P}})(0,Z)=[(\frac{b_ib_lb_\alpha b_\beta}{4\pi^2}+\frac{1}{2\pi}b_ib_lb_\beta \overline{z}'_\alpha
+\frac{1}{2\pi}b_ib_lb_\alpha \overline{z}'_\beta+b_ib_l \overline{z}'_\alpha\overline{z}'_\beta){\mathcal{P}}](0,Z).\eqno(4.57)$$
By (4.18) and (4.47), we have
$$\int_{{\mathbb{R}}^{2n}}(b_ib_l\overline{z}_\alpha \overline{ z}_\beta{\mathcal{P}})(0,Z){\mathcal{P}}(Z,0)d\nu_{TX}(Z)=0;\eqno(4.58)$$
$$\int_{{\mathbb{R}}^{2n}}\frac{1}{12}\left[({\mathcal{L}}_0^{-1}{\mathcal{O}}^{+2}_{2})^{{j-\frac{q}{2}}-k}
{\mathcal{L}}_0^{-1}({\mathcal{R}}^\bot)^k_{x_0}
R^X_{\overline{\alpha}\overline{i}\overline{\beta}
\overline{l}}b_ib_l\overline{z}_\alpha \overline{z}_\beta{\mathcal{P}}\right](0,Z){\mathcal{P}}(Z,0)d\nu_{TX}(Z)=0;\eqno(4.59)$$
Similar to (4.59), by (2.20) and (4.24), we have
$$(b_ib_l  \overline{z}_\alpha {z}_\beta{\mathcal{P}})(0,Z)=\left[\frac{b_ib_lb_\alpha z_\beta}{2\pi}+\frac{\delta_{\alpha\beta}b_ib_l}{\pi}
+b_ib_lz_\beta\overline{z}'_\alpha
{\mathcal{P}}\right](0,Z),\eqno(4.60)$$
$$\int_{{\mathbb{R}}^{2n}}\frac{1}{6}\left[({\mathcal{L}}_0^{-1}{\mathcal{O}}^{+2}_{2})^{{j-\frac{q}{2}}-k}
{\mathcal{L}}_0^{-1}({\mathcal{R}}^\bot)^k_{x_0}
R^X_{\overline{\alpha}\overline{i}{\beta}
\overline{l}}b_ib_l\overline{z}_\alpha {z}_\beta{\mathcal{P}}\right](0,Z){\mathcal{P}}(Z,0)d\nu_{TX}(Z)=0.\eqno(4.61)$$
By (4.55),(4.56),(4.59) and (4.61), then
$$\frac{1}{3}\left[({\mathcal{L}}_0^{-1}{\mathcal{O}}^{+2}_{2})^{{j-\frac{q}{2}}-k}
{\mathcal{L}}_0^{-1}({\mathcal{R}}^\bot)^k_{x_0}b_ib_l
\left<R^{TX}({\mathcal{R}},\partial_{\overline{z}_i}){\mathcal{R}},\partial_{\overline{z}_l}\right>{\mathcal{P}}\right](0,Z)$$
$$=\frac{1}{3}\frac{1}{(4\pi)^{{j-\frac{q}{2}}-k+1}}B^{{j-\frac{q}{2}},k}_2({\mathcal{R}}^\bot)
^{j-\frac{q}{2}}_{x_0}(R^X_{\alpha\overline{\alpha}\beta\overline{\beta}}
+R^X_{\alpha\overline{\beta}\beta\overline{\alpha}}){\mathcal{P}}(0,Z).\eqno(4.62)$$
We see that
$$\frac{4}{3}\left[({\mathcal{L}}_0^{-1}{\mathcal{O}}^{+2}_{2})^{{j-\frac{q}{2}}-k}
{\mathcal{L}}_0^{-1}({\mathcal{R}}^\bot)^k_{x_0}b_r\left[\left<R^{TX}_{x_0}(\partial_{{z_i}},
\partial_{\overline{z}_i}){\mathcal{R}},\partial_{\overline{z}_r}\right>
-\left<R^{TX}_{x_0}({\mathcal{R}},\partial_{z_i})\partial_{\overline{z}_i},\partial_{\overline{z}_r}\right>\right]
{\mathcal{P}}\right](0,Z)$$
$$=\frac{1}{3}\left\{({\mathcal{L}}_0^{-1}{\mathcal{O}}^{+2}_{2})^{{j-\frac{q}{2}}-k}
{\mathcal{L}}_0^{-1}({\mathcal{R}}^\bot)^k_{x_0}\left[(R^X_{i\overline{i}l\overline{r}}-R^X_{li\overline{i}\overline{r}})b_rz_l+
(R^X_{i\overline{i}\overline{l}\overline{r}}-R^X_{\overline{l}i\overline{i}\overline{r}})b_r\overline{z}_l\right]
{\mathcal{P}}\right\}(0,Z),\eqno(4.63)$$
$$\frac{1}{3}\left\{({\mathcal{L}}_0^{-1}{\mathcal{O}}^{+2}_{2})^{{j-\frac{q}{2}}-k}
{\mathcal{L}}_0^{-1}({\mathcal{R}}^\bot)^k_{x_0}\left[(R^X_{i\overline{i}l\overline{r}}-R^X_{li\overline{i}\overline{r}})b_rz_l
\right]
{\mathcal{P}}\right\}(0,Z)$$
$$=\frac{-2}{3\times (4\pi)^{{j-\frac{q}{2}}-k+1}}B^{k,{j-\frac{q}{2}}}_1({\mathcal{R}}^\bot)^{j-\frac{q}{2}}_{x_0}
(R^X_{i\overline{i}l\overline{l}}-R^X_{li\overline{i}\overline{l}}){\mathcal{P}}(0,Z),\eqno(4.64)$$

$$\int_{{\mathbb{R}}^{2n}}\frac{1}{3}\left\{({\mathcal{L}}_0^{-1}{\mathcal{O}}^{+2}_{2})^{{j-\frac{q}{2}}-k}
{\mathcal{L}}_0^{-1}({\mathcal{R}}^\bot)^k_{x_0}
(R^X_{i\overline{i}\overline{l}\overline{r}}-R^X_{\overline{l}i\overline{i}\overline{r}})b_r\overline{z}_l
{\mathcal{P}}\right\}(0,Z){\mathcal{P}}(Z,0)d\nu_{TX}(Z)=0.\eqno(4.65)$$
By (4.63)-(4.65), then
$$\frac{4}{3}\left[({\mathcal{L}}_0^{-1}{\mathcal{O}}^{+2}_{2})^{{j-\frac{q}{2}}-k}
{\mathcal{L}}_0^{-1}({\mathcal{R}}^\bot)^k_{x_0}b_r\left[\left<R^{TX}_{x_0}(\partial_{{z_i}},
\partial_{\overline{z}_i}){\mathcal{R}},\partial_{\overline{z}_r}\right>
-\left<R^{TX}_{x_0}({\mathcal{R}},\partial_{z_i})\partial_{\overline{z}_i},\partial_{\overline{z}_r}\right>\right]
{\mathcal{P}}\right](0,Z)$$
$$=\frac{-2}{3\times (4\pi)^{{j-\frac{q}{2}}-k+1}}B^{k,{j-\frac{q}{2}}}_1({\mathcal{R}}^\bot)^{j-\frac{q}{2}}_{x_0}
(R^X_{i\overline{i}l\overline{l}}-R^X_{li\overline{i}\overline{l}}){\mathcal{P}}(0,Z),\eqno(4.66)$$
Similar to (4.27) and (4.65), we have
$$-\frac{1}{3}\left[({\mathcal{L}}_0^{-1}{\mathcal{O}}^{+2}_{2})^{{j-\frac{q}{2}}-k}
{\mathcal{L}}_0^{-1}({\mathcal{R}}^\bot)^k_{x_0}
\left[{\mathcal{L}}_0,\left<R^{TX}_{x_0}({\mathcal{R}},\partial_{z_j}){\mathcal{R}},\partial_{\overline{z}_j}\right>\right]
{\mathcal{P}}\right](0,Z)$$
$$=-\frac{1}{6}\left[({\mathcal{L}}_0^{-1}{\mathcal{O}}^{+2}_{2})^{{j-\frac{q}{2}}-k}
{\mathcal{L}}_0^{-1}({\mathcal{R}}^\bot)^k_{x_0}(R^X_{\overline{r}
l\alpha\overline{l}}+R^X_{{\alpha}l\overline{r}\overline{l}})b_rz_\alpha{\mathcal{P}}\right](0,Z)$$
$$-\frac{1}{6}\left[({\mathcal{L}}_0^{-1}{\mathcal{O}}^{+2}_{2})^{{j-\frac{q}{2}}-k}
{\mathcal{L}}_0^{-1}({\mathcal{R}}^\bot)^k_{x_0}(R^X_{\overline{r}
l\overline{\alpha}\overline{l}}+R^X_{\overline{\alpha}l\overline{r}\overline{l}})b_r\overline{z}_\alpha{\mathcal{P}}\right](0,Z)
$$
$$=\frac{1}{3}\frac{1}{(4\pi)^{{j-\frac{q}{2}}-k+1}}B^{k,{j-\frac{q}{2}}}_1({\mathcal{R}}^\bot)^{j-\frac{q}{2}}_{x_0}
(R^X_{\overline{r}
lr\overline{l}}+R^X_{rl\overline{r}\overline{l}}){\mathcal{P}}(0,Z).\eqno(4.67)$$
By (1.34) and
$$z_l\overline{z}_\alpha{\mathcal{P}}=(\frac{b_\alpha z_l}{2\pi}+\frac{\delta_{\alpha l}}{\pi}
+z_l\overline{z}'_{\alpha}){\mathcal{P}},\eqno(4.68)$$
then
$$\left[({\mathcal{L}}_0^{-1}{\mathcal{O}}^{+2}_{2})^{{j-\frac{q}{2}}-k}
{\mathcal{L}}_0^{-1}({\mathcal{R}}^\bot)^k_{x_0}\Phi({\mathcal{R}},{\mathcal{R}}){\mathcal{P}}\right](0,Z)~~~~~~~~~~~~~~~~~~~~~~~~~~~$$
$$=\left\{({\mathcal{L}}_0^{-1}{\mathcal{O}}^{+2}_{2})^{{j-\frac{q}{2}}-k}
{\mathcal{L}}_0^{-1}({\mathcal{R}}^\bot)^k_{x_0}\left([\Phi(\partial_{z_l},\partial_{\overline{z}_\alpha})
+\Phi(\partial_{\overline{z}_\alpha},\partial_{z_l})]z_l\overline{z}_\alpha+
\Phi(\partial_{\overline{z}_l},\partial_{\overline{z}_\alpha})\overline{z}_l\overline{z}_\alpha\right)
{\mathcal{P}}\right\}(0,Z)$$
$$=\frac{1}{2\pi}\frac{1}{(4\pi)^{{j-\frac{q}{2}}-k+1}}(B^{k,{j-\frac{q}{2}}}_0-B^{k,{j-\frac{q}{2}}}_1)({\mathcal{R}}^\bot)^{j-\frac{q}{2}}_{x_0}[
\Phi(w_l,\overline{w}_l)+\Phi(\overline{w}_l,w_l)]{\mathcal{P}}(0,Z).\eqno(4.69)$$
We know that $$(z_lz_mz_\beta z_\alpha{\mathcal{P}})(0,Z)=0,~~(z_lz_mz_\beta \overline{z}_\alpha{\mathcal{P}})(0,Z)=
\frac{1}{2\pi}(b_\alpha z_lz_mz_\beta{\mathcal{P}})(0,Z)=0;\eqno(4.70)$$
$$(z_lz_m\overline{z}_\beta \overline{z}_\alpha{\mathcal{P}})(0,Z)=
\frac{1}{4\pi^2}[(b_\beta b_\alpha z_lz_m+2\delta_{l\alpha}b_\beta z_m+2\delta_{m\alpha}b_\beta z_l$$
$$
+2\delta_{\beta l}b_\alpha z_m+2\delta_{m\beta }b_\alpha z_l
+4\delta_{m\beta}\delta_{l\alpha}+4\delta_{\beta l}\delta_{m\alpha}){\mathcal{P}}](0,Z);\eqno(4.71)$$

$$\int_{{\mathbb{R}}^{2n}}(z_l\overline{z}_m\overline{z}_\beta \overline{z}_\alpha{\mathcal{P}})(0,Z){\mathcal{P}}(Z,0)d\nu_{TX}(Z)=0;\eqno(4.72)$$

$$\int_{{\mathbb{R}}^{2n}}(\overline{z}_l\overline{z}_m\overline{z}_\beta \overline{z}_\alpha{\mathcal{P}})(0,Z){\mathcal{P}}(Z,0)d\nu_{TX}(Z)=0.\eqno(4.73)$$
Let the coefficient of $z_lz_m\overline{z}_\beta\overline{z}_\alpha$ for $l\leq m$ and $\beta\leq\alpha$ in the expansion of
$\Psi({\mathcal{R}},{\mathcal{R}},{\mathcal{R}},{\mathcal{R}})$ be $A_{lm\overline{\beta}\overline{\alpha}}$, then by (4.70)-(4.73), we have
$$\left[({\mathcal{L}}_0^{-1}{\mathcal{O}}^{+2}_{2})^{{j-\frac{q}{2}}-k}
{\mathcal{L}}_0^{-1}({\mathcal{R}}^\bot)^k_{x_0}\Psi({\mathcal{R}},{\mathcal{R}},{\mathcal{R}},{\mathcal{R}})
{\mathcal{P}}\right](0,Z)$$
$$=\left[({\mathcal{L}}_0^{-1}{\mathcal{O}}^{+2}_{2})^{{j-\frac{q}{2}}-k}
{\mathcal{L}}_0^{-1}({\mathcal{R}}^\bot)^k_{x_0}\sum_{l\leq m}\sum_{\beta\leq \alpha}A_{lm\overline{\beta}\overline{\alpha}}
z_lz_m\overline{z}_\beta\overline{z}_\alpha{\mathcal{P}}\right](0,Z)$$
$$=\frac{1}{\pi^2}\frac{1}{(4\pi)^{{j-\frac{q}{2}}-k+1}}(B^{k,{j-\frac{q}{2}}}_0+B^{k,{j-\frac{q}{2}}}_2-2B^{k,{j-\frac{q}{2}}}_1)
({\mathcal{R}}^\bot)^{j-\frac{q}{2}}_{x_0}
(2A_{ll\overline{l}\overline{l}}+\sum_{l<m}A_{lm\overline{l}\overline{m}})$$
$$=\frac{1}{4\pi^2}\frac{1}{(4\pi)^{{j-\frac{q}{2}}-k+1}}(B^{k,{j-\frac{q}{2}}}_0+B^{k,{j-\frac{q}{2}}}_2-2B^{k,{j-\frac{q}{2}}}_1)
({\mathcal{R}}^\bot)^{j-\frac{q}{2}}_{x_0}~~~~~~~~~~~~~~~~~~~~~~~~~~$$
$$~~~~~~~~~~\times(2\Psi_{sl\overline{s}\overline{l}}+2\Psi_{sl\overline{l}\overline{s}}+2\Psi_{s\overline{l}\overline{s}l}
+2\Psi_{\overline{l}ls\overline{s}}+\Psi_{l\overline{l}s\overline{s}}+\Psi_{\overline{l}l\overline{s}s}
+\Psi_{s\overline{l}l\overline{s}}+\Psi_{\overline{l}s\overline{s}l}).\eqno(4.74)$$
By (4.24) and direct computations, we have
\begin{eqnarray*}
&&-\frac{1}{(4\pi)^k2^kk!}\frac{\pi\sqrt{-1}}{2}I_{2j}({\mathcal{L}}_0^{-1}{\mathcal{O}}^{+2}_{2})^{{j-\frac{q}{2}}-k}
{\mathcal{L}}_0^{-1}\left[-\frac{1}{2}\left<(\nabla^B\nabla^B{\bf J})_{({{\mathcal{R}}},{{\mathcal{R}}})}w_l,\overline{w}_l\right>
\right.\\
&&\left.+
 (\sum_{l,m\leq q}+\sum_{l,m\geq q+1})\left<(\nabla^B\nabla^B{\bf J})_{({{\mathcal{R}}},{{\mathcal{R}}})}w_l,\overline{w}_m\right>\overline{w}^m
 \wedge i_{\overline{w}_l}\right]
 (x_0)({\mathcal{R}}^\bot)^k_{x_0}{{P}^N}(0,Z)\\
 &=&\frac{\sqrt{-1}}{2(4\pi)^{j-\frac{q}{2}+1}}B_0^{1,k}(B_1^{k,j-\frac{q}{2}}-B_0^{k,j-\frac{q}{2}})
({\mathcal{R}}^\bot)^{j-\frac{q}{2}-k}_{x_0}\\
&&
\left[-\frac{1}{2}\left<(\nabla^B\nabla^B{\bf J})_{(\partial_{z_\alpha},\partial_{\overline{z}_\alpha})}w_l,\overline{w}_l\right>
-\frac{1}{2}\left<(\nabla^B\nabla^B{\bf J})_{(\partial_{\overline{z}_\alpha},\partial_{z_\alpha})}w_l,\overline{w}_l\right>\right.\\
&&
+ (\sum_{l,m\leq q}+\sum_{l,m\geq q+1})\left<(\nabla^B\nabla^B{\bf J})_{(\partial_{z_\alpha},\partial_{\overline{z}_\alpha})}w_l,\overline{w}_m\right>\overline{w}^m
 \wedge i_{\overline{w}_l}\\
 &&\left.
 +(\sum_{l,m\leq q}+\sum_{l,m\geq q+1})\left<(\nabla^B\nabla^B{\bf J})_{(\partial_{\overline{z}_\alpha},\partial_{z_\alpha})}w_l,\overline{w}_m\right>\overline{w}^m
 \wedge i_{\overline{w}_l}
 \right]({\mathcal{R}}^\bot)^k_{x_0}{{P}^N}(0,Z).~(4.75)
\end{eqnarray*}
By (4.41)-(4.47),(4.54),(4.56),(4.62),(4.66),(4.67),(4.69),(4.74) and (4.75), we get

\begin{eqnarray*}
&&I_{2j}({\mathcal{L}}_0^{-1}{\mathcal{O}}^{+2}_{2})^{{j-\frac{q}{2}}-k}
({\mathcal{L}}_0^{-1}{\mathcal{O}}^{0,1}_{2})
({\mathcal{L}}_0^{-1}{\mathcal{O}}^{+2}_{2})^{k}P^N
({\mathcal{O}}^{-2}_{2}{\mathcal{L}}_0^{-1})^{{j-\frac{q}{2}}}I_{2j}~~~~~~~~~~~~~~~~~~(4.76)\\
&=&\frac{1}{(4\pi)^{2j-q+1}2^{{j-\frac{q}{2}}+k}({j-\frac{q}{2}})!k!}I_{2j}({\mathcal{R}}^\bot)(x)^{{j-\frac{q}{2}}-k}\\
&&
\cdot\left\{-B^{k,{j-\frac{q}{2}}}_1R^E(w_s,\overline{w_s})+(B^{k,{j-\frac{q}{2}}}_0-B^{k,{j-\frac{q}{2}}}_1)
\widetilde{R}^{\wedge^{0,\bullet},B}(w_s,\overline{w_s})\right.\\
&&+B^{k,{j-\frac{q}{2}}}_0\left[2(\sum_{l,m\leq q}+\sum_{l,m\geq q+1})(R^E+\frac{1}{2}R^{\rm det})(w_s,\overline{w_m})
\overline{w^m}\wedge i_{\overline{w_s}}\right.\\
&&\left.-\frac{1}{4}\widetilde{c(dT_{as})}
+\frac{r^X}{4}-\frac{1}{8}|T_{as}|^2-\frac{1}{2}R^{{\rm det}}(w_s,\overline{w_s})\right]
+\frac{\sqrt{-1}}{2}(B_1^{k,j-\frac{q}{2}}-B_0^{k,j-\frac{q}{2}})\\
&&
\cdot\left[-\frac{1}{2}\left<(\nabla^B\nabla^B{\bf J})_{(\partial_{z_\alpha},\partial_{\overline{z}_\alpha})}w_l,\overline{w}_l\right>
-\frac{1}{2}\left<(\nabla^B\nabla^B{\bf J})_{(\partial_{\overline{z}_\alpha},\partial_{z_\alpha})}w_l,\overline{w}_l\right>\right.\\
&&
+ (\sum_{l,m\leq q}+\sum_{l,m\geq q+1})\left<(\nabla^B\nabla^B{\bf J})_{(\partial_{z_\alpha},\partial_{\overline{z}_\alpha})}w_l,\overline{w}_m\right>\overline{w}^m
 \wedge i_{\overline{w}_l}\\
  &&\left.\left.
 +(\sum_{l,m\leq q}+\sum_{l,m\geq q+1})\left<(\nabla^B\nabla^B{\bf J})_{(\partial_{\overline{z}_\alpha},\partial_{z_\alpha})}w_l,\overline{w}_m\right>\overline{w}^m
 \wedge i_{\overline{w}_l}
 \right]
\right\}
\\
&&\cdot({\mathcal{R}}^\bot)(x)^kI_{{\rm det}(\overline{W}^*)\otimes E}({\mathcal{R}}^{\bot,*})(x)^{j-\frac{q}{2}}I_{2j}\\,
&&+\frac{1}{(4\pi)^{2j-q+1}2^{j-\frac{q}{2}+k}(j-\frac{q}{2})!k!}I_{2j}\left\{(B^{k,{j-\frac{q}{2}}}_0
-\frac{1}{3}B^{k,{j-\frac{q}{2}}}_1)R^{TX}_{\alpha\beta\overline{\alpha}\overline{\beta}}\right.\\
&&
+(\frac{1}{3}B^{k,{j-\frac{q}{2}}}_2-\frac{2}{3}B^{k,{j-\frac{q}{2}}}_1)R^{TX}_{\alpha\overline{\alpha}\beta\overline{\beta}}
+(-\frac{1}{3}B^{k,{j-\frac{q}{2}}}_2+\frac{1}{3}B^{k,{j-\frac{q}{2}}}_1)R^{TX}_{\overline{\alpha}\beta\alpha\overline{\beta}}\\
&&
+\frac{1}{2\pi}(B^{k,{j-\frac{q}{2}}}_0-B^{k,{j-\frac{q}{2}}}_1)[\Phi(w_l,\overline{w_l})+\Phi(\overline{w_l},w_l)]\\
&&
+\frac{1}{4\pi^2}(B^{k,{j-\frac{q}{2}}}_0+B^{k,{j-\frac{q}{2}}}_2-2B^{k,{j-\frac{q}{2}}}_1)\\.
&&\cdot
(2\Psi_{sl\overline{s}\overline{l}}+2\Psi_{sl\overline{l}\overline{s}}+2\Psi_{s\overline{l}\overline{s}l}
+2\Psi_{\overline{l}ls\overline{s}}+\Psi_{l\overline{l}s\overline{s}}+\Psi_{\overline{l}l\overline{s}s}
+\Psi_{s\overline{l}l\overline{s}}+\Psi_{\overline{l}s\overline{s}l})\\
&&\left.
+(B^{k,{j-\frac{q}{2}}}_2-B^{k,{j-\frac{q}{2}}}_1)(\widehat{\Psi}_{\beta\overline{\beta}\alpha\overline{\alpha}}
+\widehat{\Psi}_{\overline{\alpha}\beta\alpha\overline{\beta}})\right\}
({\mathcal{R}}^\bot)(x)^{j-\frac{q}{2}}I_{{\rm det}(\overline{W}^*)\otimes E}({\mathcal{R}}^{\bot,*})(x)^{j-\frac{q}{2}}I_{2j}.
\end{eqnarray*}

In the following, we assume that $n=4$, $q=2$, $j=2$. Then by ${\mathcal{O}}^{0,2}_{2}$ changing the double degree and preserving the total degree, so
\begin{eqnarray*}
&&\sum_{k=0}^{j-\frac{q}{2}}I_{2j}({\mathcal{L}}_0^{-1}{\mathcal{O}}^{+2}_{2})^{{j-\frac{q}{2}}-k}
({\mathcal{L}}_0^{-1}{\mathcal{O}}^{0,2}_{2})
({\mathcal{L}}_0^{-1}{\mathcal{O}}^{+2}_{2})^{k}P^N
({\mathcal{O}}^{-2}_{2}{\mathcal{L}}_0^{-1})^{{j-\frac{q}{2}}}I_{2j}\\
&=&I_{4}({\mathcal{L}}_0^{-1}{\mathcal{O}}^{+2}_{2})
({\mathcal{L}}_0^{-1}{\mathcal{O}}^{0,2}_{2})
P^N
({\mathcal{O}}^{-2}_{2}{\mathcal{L}}_0^{-1})I_{4}
+I_{4}
({\mathcal{L}}_0^{-1}{\mathcal{O}}^{0,2}_{2})
({\mathcal{L}}_0^{-1}{\mathcal{O}}^{+2}_{2})P^N
({\mathcal{O}}^{-2}_{2}{\mathcal{L}}_0^{-1})I_{4}\\
&=&I_{4}({\mathcal{L}}_0^{-1}{\mathcal{O}}^{+2}_{2})
({\mathcal{L}}_0^{-1}{\mathcal{O}}^{0,2}_{2})
P^N
({\mathcal{O}}^{-2}_{2}{\mathcal{L}}_0^{-1})I_{4}.~~~~~~~~~~~~~~~~~~~~~~~~~~~~~~~~~~~~~~~~~~~~~~(4.77)
\end{eqnarray*}
By (4.40), we have
\begin{eqnarray*}
{\mathcal{O}}^{0,2}_2P^N&=&
\left\{-\sum_{l\leq q,m\geq q+1}\left<R^B
({\mathcal{R}},e_i)w_l,\overline{w}_m\right>{\overline{w}}^m\wedge i_{\overline{w}_l}\nabla_{0,e_i}\right.\\
&&
+2\sum_{l\leq q,m\geq q+1}(R^E+\frac{1}{2}R^{\rm det})(w_l,\overline{w}_m)\overline{w}^m\wedge i_{\overline{w}_l}\\
&&-\frac{\pi}{2}\sqrt{-1}\sum_{l\leq q,m\geq q+1}
 \left<(\nabla^B\nabla^B{\bf J})_{({{\mathcal{R}}},{{\mathcal{R}}})}w_l,\overline{w}_m\right>\overline{w}^m
 \wedge i_{\overline{w}_l}\\
&& +\frac{1}{2}\sum_{l\leq q,m\geq q+1}dT_{as}(w_l,\overline{w}_m,w_k,\overline{w}_k)\overline{w}^m
 \wedge i_{\overline{w}_l}\\
 &&\left.
 -\frac{1}{4}
\sum_{1\leq i,j'\leq q}~\sum_{k~{\rm or}~l\geq q+1}dT_{as}(w_i,w_{j'},\overline{w}_k,\overline{w}_l)\overline{w}^k\wedge\overline{w}^l\wedge
i_{\overline{w}_i}i_{\overline{w}_{j'}}\right\}P^N.~~(4.78)
\end{eqnarray*}
By (2.21), we have
$$I_{4}({\mathcal{L}}_0^{-1}{\mathcal{O}}^{+2}_{2})
{\mathcal{L}}_0^{-1}\left[2\sum_{l\leq q,m\geq q+1}(R^E+\frac{1}{2}R^{\rm det})(w_l,\overline{w}_m)\overline{w}^m\wedge i_{\overline{w}_l}\right.$$
$$\left. +\frac{1}{2}\sum_{l\leq q,m\geq q+1}dT_{as}(w_l,\overline{w}_m,w_k,\overline{w}_k)\overline{w}^m
 \wedge i_{\overline{w}_l}\right]P^N$$
 $$
 =\frac{1}{64\pi^2}{\mathcal{R}}^0\left[2\sum_{l\leq q,m\geq q+1}(R^E+\frac{1}{2}R^{\rm det})(w_l,\overline{w}_m)\overline{w}^m\wedge i_{\overline{w}_l}
 \right.$$
 $$\left.
 +\frac{1}{2}\sum_{l\leq q,m\geq q+1}dT_{as}(w_l,\overline{w}_m,w_k,\overline{w}_k)\overline{w}^m
 \wedge i_{\overline{w}_l}\right]P^N.\eqno(4.79)$$
\noindent Similar to (4.75), it holds that
$$I_{4}({\mathcal{L}}_0^{-1}{\mathcal{O}}^{+2}_{2})
{\mathcal{L}}_0^{-1}\left[-\frac{\pi}{2}\sqrt{-1}\sum_{l\leq q,m\geq q+1}
 \left<(\nabla^B\nabla^B{\bf J})_{({{\mathcal{R}}},{{\mathcal{R}}})}w_l,\overline{w}_m\right>\overline{w}^m
 \wedge i_{\overline{w}_l}\right]P^N(0,Z)$$
 $$=\frac{-5\sqrt{-1}}{32\times 36\pi^2}{\mathcal{R}}^0\left[\sum_{l\leq q,m\geq q+1}
 \left<(\nabla^B\nabla^B{\bf J})_{(\partial_{z_\alpha},\partial_{{\overline{z}_\alpha}})}w_l,\overline{w}_m\right>\overline{w}^m
 \wedge i_{\overline{w}_l}\right.$$
 $$\left.
 + \left<(\nabla^B\nabla^B{\bf J})_{(\partial_{{\overline{z}_\alpha}},\partial_{z_\alpha})}w_l,\overline{w}_m\right>\overline{w}^m
 \wedge i_{\overline{w}_l}\right]P^N(0,Z).\eqno(4.80)$$
 By (4.16) and (4.17), we get
 $$-I_{4}({\mathcal{L}}_0^{-1}{\mathcal{O}}^{+2}_{2})
{\mathcal{L}}_0^{-1}\left[\sum_{l\leq q,m\geq q+1}\left<R^B
({\mathcal{R}},e_i)w_l,\overline{w}_m\right>{\overline{w}}^m\wedge i_{\overline{w}_l}\nabla_{0,e_i}\right]P^N(0,Z)$$
$$=\frac{5}{288\pi^2}{\mathcal{R}}^0\left[\sum_{l\leq q,m\geq q+1}\left<R^B(\partial_{z_\alpha},\partial_{{\overline{z}_\alpha}})w_l,\overline{w}_m\right>{\overline{w}}^m\wedge i_{\overline{w}_l}\right]P^N(0,Z).\eqno(4.81)$$
 Similarly
 $$-\frac{1}{4}I_{4}({\mathcal{L}}_0^{-1}{\mathcal{O}}^{+2}_{2})
{\mathcal{L}}_0^{-1}\sum_{1\leq i,j'\leq q}~\sum_{k~{\rm or}~l\geq q+1}dT_{as}(w_i,w_{j'},\overline{w}_k,\overline{w}_l)\overline{w}^k\wedge\overline{w}^l\wedge
i_{\overline{w}_i}i_{\overline{w}_{j'}}P^N(0,Z)$$
$$=-\frac{1}{512\pi^2}{\mathcal{R}}^\top\sum_{1\leq i,j'\leq q}\sum_{k,l\geq q+1}
 dT_{as}(w_i,w_{j'},\overline{w}_k,\overline{w}_l)\overline{w}^k\wedge\overline{w}^l\wedge
i_{\overline{w}_i}i_{\overline{w}_{j'}}P^N(0,Z)$$
$$-\frac{1}{256\pi^2}{\mathcal{R}}^0\sum_{1\leq i,j'\leq q}(\sum_{k\leq q,l\geq q+1}+\sum_{l\leq q,k\geq q+1})
 dT_{as}(w_i,w_{j'},\overline{w}_k,\overline{w}_l)\overline{w}^k\wedge\overline{w}^l\wedge
i_{\overline{w}_i}i_{\overline{w}_{j'}}P^N(0,Z).\eqno(4.82)$$
 When $n=4$, $q=2$, $j=2$, we have
 $${\rm III}_a(0,0)=\frac{1}{4\pi}[(4.79)+(4.80)+(4.81)+(4.82)]I_{{\rm det}(\overline{W}^*)\otimes E}
( {\mathcal{R}}^\bot)^*I_4+\sum_{k=0}^1(4.76).\eqno(4.83)$$

\noindent{\bf 4.3 The computations of the term V}\\

\indent In this section, we will compute the term V. By the discussions after (4.2),
for $l=2$, then $4j\leq 2q+2k-4\leq 4j$, so $k=2j+2-q$. There are
$2j-q$ ${\mathcal{O}}_{r_i}$ equal to ${\mathcal{O}}_{2}$ and $2$ equal to ${\mathcal{O}}_{1}$;
By (3.16) and (3.17), when $l=2$ and $l_1=0$, then $i_0=j-\frac{q}{2}+3$, when $l=2$ and $l_1=1$, then $i_0=j-\frac{q}{2}+2$,
when $l=2$ and $l_1=2$, then $i_0=j-\frac{q}{2}+1$, so we only have one $\eta_{i_0}=N$.
These three cases correspond to the terms ${\rm V}_a$, ${\rm V}_b$, ${\rm V}^*_a$ and
$$ {\rm V}_a=\sum_{0\leq m_1+m_2\leq {j-\frac{q}{2}}}I_{2j}({\mathcal{L}}^{-1}_0{\mathcal{O}}^{+2}_2)^{m_1}({\mathcal{L}}^{-1}_0P^{N^\bot}{\mathcal{O}}_1)({\mathcal{L}}^{-1}_0{\mathcal{O}}^{+2}_2)^{m_2}$$
$$\times({\mathcal{L}}^{-1}_0P^{N^\bot}{\mathcal{O}}_1)
({\mathcal{L}}^{-1}_0{\mathcal{O}}^{+2}_2)^{{j-\frac{q}{2}}-m_1-m_2}P^N({\mathcal{O}}^{-2}_2{\mathcal{L}}^{-1}_0)^{j-\frac{q}{2}}I_{2j},
\eqno(4.84)$$
$${\rm V}_b=
\sum_{0\leq m_1,m_2\leq {j-\frac{q}{2}}}I_{2j}({\mathcal{L}}^{-1}_0{\mathcal{O}}^{+2}_2)^{m_1}
({\mathcal{L}}^{-1}_0P^{N^\bot}{\mathcal{O}}_1)({\mathcal{L}}^{-1}_0{\mathcal{O}}^{+2}_2)^{{j-\frac{q}{2}}-m_1}
P^N$$
$$\times({\mathcal{O}}^{-2}_2{\mathcal{L}}^{-1}_0)^{m_2}
({\mathcal{O}}_1{\mathcal{L}}^{-1}_0P^{N^\bot})({\mathcal{O}}^{-2}_2{\mathcal{L}}^{-1}_0)^{{j-\frac{q}{2}}-m_2}I_{2j},\eqno(4.85)$$
$${\rm V}={\rm V}_a+{\rm V}_b+{\rm V}^*_a.\eqno(4.86)$$
By (2.12) and (2.3) and Proposition 2.1 in [LuW1], we have
$${\mathcal{O}}_1={\mathcal{O}}'_1+{\mathcal{O}}''_1,\eqno(4.87)$$
where $${\mathcal{O}}'_1=-\frac{2}{3}(\partial_sR^L)_{x_0}({\mathcal{R}},e_i)Z_s\nabla_{0,e_i}-\frac{1}{3}(\partial_sR^L)_{x_0}({\mathcal{R}},e_s)
;\eqno(4.88)$$
$${\mathcal{O}}''_1=
-4\pi\sqrt{-1}(\sum_{l\leq q}\sum_{q+1\leq m}+\sum_{m\leq q}\sum_{q+1\leq l})\left<(\nabla^B_{{\mathcal{R}}}{\bf J})w_l,\overline{w}_m\right>\overline{w}^m\wedge i_{\overline{w}_l}.
\eqno(4.89)$$
Firstly, we compute the term
$$ {\rm V}^0_a:=\sum_{0\leq m_1+m_2\leq {j-\frac{q}{2}}}I_{2j}({\mathcal{L}}^{-1}_0{\mathcal{O}}^{+2}_2)^{m_1}({\mathcal{L}}^{-1}_0
P^{N^\bot}{\mathcal{O}}'_1)({\mathcal{L}}^{-1}_0{\mathcal{O}}^{+2}_2)^{m_2}$$
$$\times({\mathcal{L}}^{-1}_0P^{N^\bot}{\mathcal{O}}'_1)
({\mathcal{L}}^{-1}_0{\mathcal{O}}^{+2}_2)^{{j-\frac{q}{2}}-m_1-m_2}P^N({\mathcal{O}}^{-2}_2{\mathcal{L}}^{-1}_0)^{j-\frac{q}{2}}I_{2j},
\eqno(4.90)$$
By
$$({\mathcal{L}}^{-1}_0{\mathcal{O}}^{+2}_2)^{{j-\frac{q}{2}}-m_1-m_2}P^N=\frac{1}{(4\pi)^{{j-\frac{q}{2}}-m_1-m_2}}B^{1,{j-\frac{q}{2}}-m_1-m_2}_0
( {\mathcal{R}}^\bot)^{{j-\frac{q}{2}}-m_1-m_2}{{P}^N},\eqno(4.91)$$
$$\frac{2}{3}(\partial_lR^L)_{x_0}(\partial_{z_\alpha},\partial_{\overline{z}_i}
)Z_lz_\alpha b_i
=\frac{2}{3}(\partial_{z_\beta }R^L)_{x_0}(\partial_{z_\alpha},\partial_{\overline{z}_i})z_\alpha z_\beta b_i
+\frac{2}{3}(\partial_{\overline{z}_\beta }R^L)_{x_0}(\partial_{z_\alpha},\partial_{\overline{z}_i})z_\alpha \overline{z}_\beta b_i,\eqno(4.92)$$

and (2.12), (4.16), (4.17) and $R^L$ being a $(1,1)$ form, then
$${\mathcal{O}}'_1({\mathcal{L}}^{-1}_0{\mathcal{O}}^{+2}_2)^{{j-\frac{q}{2}}-m_1-m_2}P^N=\frac{1}{(4\pi)^{{j-\frac{q}{2}}-m_1-m_2}}
B^{1,{j-\frac{q}{2}}-m_1-m_2}_0
{\mathcal{R}}^{{j-\frac{q}{2}}-m_1-m_2}$$
$$\times\left[\frac{2}{3}(\partial_{z_{j'} }R^L)_{x_0}(\partial_{z_\alpha},\partial_{\overline{z}_l})z_\alpha z_{j'} b_l
+\frac{2}{3}(\partial_{\overline{z}_r }R^L)_{x_0}(\partial_{z_\alpha},\partial_{\overline{z}_l})z_\alpha \overline{z}_rb_l \right.$$
$$\left.-\frac{1}{3}(\partial_{j'}R^L)_{x_0}(\partial_{z_\alpha},e_{j'})z_\alpha-
\frac{1}{3}(\partial_{j'}R^L)_{x_0}(\partial_{\overline{z}_\alpha},e_{j'})(\frac{b_\alpha}{2\pi}+\overline{z}'_\alpha)\right]
P^N.\eqno(4.93)$$
By
$$z_\alpha z_jb_l=b_lz_\alpha z_j+2\delta_{l\alpha}z_j+2\delta_{jl}z_\alpha,\eqno(4.94)$$
and
$$z_\alpha\overline{z}_rb_l{\mathcal{P}}=(\frac{b_lb_rz_\alpha}{2\pi}+\frac{\delta_{\alpha r}b_l}{\pi}+\frac{\delta_{\alpha l}b_r}{\pi}+
b_lz_\alpha\overline{z}'_r+2\delta_{\alpha l}\overline{z}'_r){\mathcal{P}},\eqno(4.95)$$
then
\begin{eqnarray*}
&&({\mathcal{L}}^{-1}_0{\mathcal{O}}^{+2}_2)^{m_2}
({\mathcal{L}}^{-1}_0{\mathcal{O}}'_1)
({\mathcal{L}}^{-1}_0{\mathcal{O}}^{+2}_2)^{{j-\frac{q}{2}}-m_1-m_2}P^N\\
&=&\frac{1}{(4\pi)^{{j-\frac{q}{2}}-m_1+1}}B^{1,{j-\frac{q}{2}}-m_1-m_2}_0B^{{j-\frac{q}{2}}-m_1-m_2,{j-\frac{q}{2}}-m_1}_0
( {\mathcal{R}}^\bot)^{{j-\frac{q}{2}}-m_1}\\
&&
\cdot\left[\frac{4}{3}(\partial_{\overline{z}_r }R^L)_{x_0}(\partial_{z_l},\partial_{\overline{z}_l})\overline{z}'_r
+\frac{4}{3}(\partial_{{z}_{j'} }R^L)_{x_0}(\partial_{z_\alpha},\partial_{\overline{z}_\alpha})z_j\right.\\
&&
+\frac{4}{3}(\partial_{{z}_{j'} }R^L)_{x_0}(\partial_{z_\alpha},\partial_{\overline{z}_j})z_\alpha-
\frac{1}{3}(\partial_{j'}R^L)_{x_0}(\partial_{z_\alpha},e_{j'})z_\alpha\\
&&\left.-
\frac{1}{3}(\partial_{j'}R^L)_{x_0}(\partial_{\overline{z}_\alpha},e_{j'})\overline{z}'_\alpha \right]{{P}^N}\\
&&+\frac{1}{(4\pi)^{{j-\frac{q}{2}}-m_1+1}}B^{1,{j-\frac{q}{2}}-m_1-m_2}_0B^{{j-\frac{q}{2}}-m_1-m_2,{j-\frac{q}{2}}-m_1}_1
( {\mathcal{R}}^\bot)^{{j-\frac{q}{2}}-m_1}\\
&&\cdot
\left[\frac{2}{3\pi}(\partial_{\overline{z}_r }R^L)_{x_0}(\partial_{z_r},\partial_{\overline{z}_l})b_l+\frac{2}{3\pi}(\partial_{\overline{z}_r }R^L)_{x_0}(\partial_{z_\alpha},\partial_{\overline{z}_\alpha})b_r\right.\\
&&
+\frac{2}{3}(\partial_{\overline{z}_r }R^L)_{x_0}(\partial_{z_\alpha},\partial_{\overline{z}_l})b_lz_\alpha \overline{z}'_\gamma
+
\frac{2}{3}(\partial_{{z}_{j'} }R^L)_{x_0}(\partial_{z_\alpha},\partial_{\overline{z}_l})b_lz_\alpha z_{j'}\\
&&\left.-
\frac{1}{3}(\partial_{j'}R^L)_{x_0}(\partial_{\overline{z}_\alpha},e_{j'})\frac{b_\alpha}{2\pi}\right]{{P}^N}\\
&&+\frac{1}{(4\pi)^{{j-\frac{q}{2}}-m_1+1}}B^{1,{j-\frac{q}{2}}-m_1-m_2}_0B^{{j-\frac{q}{2}}-m_1-m_2,j-m_1}_2
( {\mathcal{R}}^\bot)^{{j-\frac{q}{2}}-m_1}\\
&&\times\frac{1}{3\pi}(\partial_{\overline{z}_r }R^L)_{x_0}(\partial_{z_\alpha},\partial_{\overline{z}_l})b_lb_rz_\alpha
{{P}^N}.
~~~~~~~~~~~~~~~~~~~~~~~~~~~~~~~~~~~~~~~~~~~~~~~~~~~(4.96)
\end{eqnarray*}
By (4.16) and (4.17) and $R^L$ being a $(1,1)$-form, we get
$${\mathcal{O}}'_1=A_1+A_2+A_3+A_4+A_5,\eqno(4.97)$$
where
$$A_1=-\frac{2}{3}(\partial_{{z}_\beta }R^L)_{x_0}(\partial_{\overline{z}_\alpha},\partial_{{z}_i})\overline{z}_\alpha z_\beta b^+_i-
\frac{2}{3}(\partial_{\overline{z}_r }R^L)_{x_0}(\partial_{\overline{z}_\alpha},\partial_{{z}_i})
\overline{z}_\alpha \overline{z}_r b^+_i,\eqno(4.98)$$
$$A_2=[\frac{4}{3}(\partial_{{z}_{j'} }R^L)_{x_0}(\partial_{z_\alpha},\partial_{\overline{z}_\alpha})+
\frac{4}{3}(\partial_{{z}_{\alpha} }R^L)_{x_0}(\partial_{z_{j'}},\partial_{\overline{z}_\alpha})-
\frac{1}{3}(\partial_{\alpha}R^L)_{x_0}(\partial_{z_{j'}},e_{\alpha})]z_{j'}$$
$$+
[\frac{4}{3}(\partial_{\overline{z}_r }R^L)_{x_0}(\partial_{z_l},\partial_{\overline{z}_l})-
\frac{1}{3}(\partial_{j'}R^L)_{x_0}(\partial_{\overline{z}_r},e_{j'})]\overline{z}'_r,\eqno(4.99)$$
$$A_3=[\frac{2}{3\pi}(\partial_{\overline{z}_r }R^L)_{x_0}(\partial_{z_r},\partial_{\overline{z}_\alpha})
+\frac{2}{3\pi}(\partial_{\overline{z}_\alpha }R^L)_{x_0}(\partial_{z_r},\partial_{\overline{z}_r})-
\frac{1}{6\pi}(\partial_{j'}R^L)_{x_0}(\partial_{\overline{z}_\alpha},e_{j'})]b_\alpha$$
$$
+\frac{2}{3}(\partial_{\overline{z}_r }R^L)_{x_0}(\partial_{z_\alpha},\partial_{\overline{z}_l})b_lz_\alpha \overline{z}'_\gamma
+
\frac{2}{3}(\partial_{{z}_{j'} }R^L)_{x_0}(\partial_{z_\alpha},\partial_{\overline{z}_l})b_lz_\alpha z_{j'},\eqno(4.100)$$
$$A_4=\frac{1}{3\pi}(\partial_{\overline{z}_r }R^L)_{x_0}(\partial_{z_\alpha},\partial_{\overline{z}_l})b_lb_rz_\alpha. \eqno(4.101)$$
$$A_5=
\frac{2}{3}(\partial_{\overline{z}_r }R^L)_{x_0}(\partial_{z_\alpha},\partial_{\overline{z}_l})(b_lz_\alpha+2\delta_{\alpha l})
(\overline{z}_r-\frac{b_r}{2\pi}-\overline{z}'_r)$$
$$
-\frac{1}{3}(\partial_{j'}R^L)_{x_0}(\partial_{\overline{z}_r},e_{j'})(\overline{z}_r-\frac{b_r}{2\pi}-\overline{z}'_r).\eqno(4.102)$$
By
$$\int_{{\mathbb{R}}^{2n}}(b_l\overline{z}'_r{\mathcal{P}})(0,Z){\mathcal{P}}(Z,0)d\nu_{TX}(Z)=0;\eqno(4.103)$$
$$\int_{{\mathbb{R}}^{2n}}(b_lz_\alpha \overline{z}'_r\overline{z}'_s{\mathcal{P}})(0,Z){\mathcal{P}}(Z,0)d\nu_{TX}(Z)=0;\eqno(4.104)$$
$$\int_{{\mathbb{R}}^{2n}}(b_lz_\alpha\overline{z}'_rz_{j'}{\mathcal{P}})(Z,0){\mathcal{P}}(Z,0)d\nu_{TX}(Z)=0;\eqno(4.105)$$
$$(b_lz_{j'}{\mathcal{P}})(0,Z)=-2\delta_{j'l}{\mathcal{P}}(0,Z),\eqno(4.106)$$
$$(b_lz_\alpha z_{j'}z_s{\mathcal{P}})(0,Z)=0~~
(b_lz_\alpha z_{j'}\overline{z}'_s{\mathcal{P}})(0,Z)=0,\eqno(4.107)$$
then
$$A_1A_2{\mathcal{P}}=0;\eqno(4.108)$$
$$\int_{{\mathbb{R}}^{2n}}(A^2_2{\mathcal{P}})(0,Z){\mathcal{P}}(Z,0)d\nu_{TX}(Z)=0;\eqno(4.109)$$
\begin{eqnarray*}
&&\int_{{\mathbb{R}}^{2n}}\frac{1}{(4\pi)^{{j-\frac{q}{2}}-m_1+1}}B^{1,{j-\frac{q}{2}}-m_1-m_2}_0
B^{{j-\frac{q}{2}}-m_1-m_2,j-m_1}_0({\mathcal{L}}^{-1}_0{\mathcal{O}}^{+2}_2)^{m_1}\\
&&\cdot{\mathcal{L}}^{-1}_0
( {\mathcal{R}}^\bot)^{{j-\frac{q}{2}}-m_1}
(A_3A_2{\mathcal{P}})(0,Z){\mathcal{P}}(Z,0)d\nu_{TX}(Z)\\
&=&\frac{-2}{(4\pi)^{{j-\frac{q}{2}}+2}}B^{1,{j-\frac{q}{2}}-m_1-m_2}_0B^{{j-\frac{q}{2}}-m_1-m_2,j-m_1}_0
B^{{j-\frac{q}{2}}-m_1,{j-\frac{q}{2}}}_1( {\mathcal{R}}^\bot)^{{j-\frac{q}{2}}}\\
&&\times
[\frac{2}{3\pi}(\partial_{\overline{z}_r }R^L)_{x_0}(\partial_{z_r},\partial_{\overline{z}_\alpha})
+\frac{2}{3\pi}(\partial_{\overline{z}_\alpha }R^L)_{x_0}(\partial_{z_r},\partial_{\overline{z}_r})
-
\frac{1}{6\pi}(\partial_{j'}R^L)_{x_0}(\partial_{\overline{z}_\alpha},e_{j'})]\\
&&
\times
[\frac{4}{3}(\partial_{{z}_{\alpha} }R^L)_{x_0}(\partial_{z_{\alpha'}},\partial_{\overline{z}_{\alpha'}})+
\frac{4}{3}(\partial_{{z}_{\alpha'} }R^L)_{x_0}(\partial_{z_{\alpha}},\partial_{\overline{z}_{\alpha'}})
-
\frac{1}{3}(\partial_{\alpha'}R^L)_{x_0}(\partial_{z_{\alpha}},e_{\alpha'})];~(4.110)
\end{eqnarray*}

\begin{eqnarray*}
&&\int_{{\mathbb{R}}^{2n}}\frac{1}{(4\pi)^{{j-\frac{q}{2}}-m_1+1}}B^{1,{j-\frac{q}{2}}-m_1-m_2}_0
B^{{j-\frac{q}{2}}-m_1-m_2,{j-\frac{q}{2}}-m_1}_0({\mathcal{L}}^{-1}_0{\mathcal{O}}^{+2}_2)^{m_1}\\
&&\cdot{\mathcal{L}}^{-1}_0
( {\mathcal{R}}^\bot)^{{j-\frac{q}{2}}-m_1}
(A_4A_2{\mathcal{P}})(0,Z){\mathcal{P}}(Z,0)d\nu_{TX}(Z)\\
&=&\frac{1}{(4\pi)^{{j-\frac{q}{2}}+2}}B^{1,{j-\frac{q}{2}}-m_1-m_2}_0B^{{j-\frac{q}{2}}-m_1-m_2,{j-\frac{q}{2}}-m_1}_0
B^{{j-\frac{q}{2}}-m_1,{j-\frac{q}{2}}}_2( {\mathcal{R}}^\bot)^{{j-\frac{q}{2}}}\\
&&\cdot\left\{
\frac{4}{3\pi}(\partial_{\overline{z}_r }R^L)_{x_0}(\partial_{z_l},\partial_{\overline{z}_l})\left[\frac{4}{3}(\partial_{{z}_{r} }R^L)_{x_0}(\partial_{z_\alpha},\partial_{\overline{z}_\alpha})\right.\right.\\
&&\left.+
\frac{4}{3}(\partial_{{z}_{\alpha} }R^L)_{x_0}(\partial_{z_{r}},\partial_{\overline{z}_\alpha})-
\frac{1}{3}(\partial_{\alpha}R^L)_{x_0}(\partial_{z_{r}},e_{\alpha})\right]\\
&&+\frac{4}{3\pi}(\partial_{\overline{z}_\alpha }R^L)_{x_0}(\partial_{z_\alpha},\partial_{\overline{z}_l})
\left[\frac{4}{3}(\partial_{{z}_{l} }R^L)_{x_0}(\partial_{z_\alpha},\partial_{\overline{z}_\alpha})\right.\\
&&\left.\left.+
\frac{4}{3}(\partial_{{z}_{\alpha} }R^L)_{x_0}(\partial_{z_{l}},\partial_{\overline{z}_\alpha})-
\frac{1}{3}(\partial_{\alpha}R^L)_{x_0}(\partial_{z_{l}},e_{\alpha})\right]\right\}.
~~~~~~~~~~~~~~~~~~~~~~~~~~~~~~(4.111)
\end{eqnarray*}

By (2.20) and (4.17), we have
$$(\overline{z}_\alpha z_\beta b^+_ib_{\alpha'}{\mathcal{P}})(0,Z)=[(2\delta_{i\alpha'}b_\alpha z_\beta+4\delta_{i\alpha'}\delta_{\alpha\beta})
{\mathcal{P}}](0,Z),\eqno(4.112)$$
$$(\overline{z}_\alpha z_\beta b^+_ib_lz_{\alpha'}\overline{z}'_r{\mathcal{P}})(0,Z)=0,\eqno(4.113)$$
$$(\overline{z}_\alpha z_\beta b^+_ib_lz_{\alpha'}z_{j'}{\mathcal{P}})(0,Z)=0,\eqno(4.114)$$
$$\int_{{\mathbb{R}}^{2n}}(\overline{z}_\alpha \overline{z}_\beta b^+_ib_{\alpha'}{\mathcal{P}})(0,Z)
{\mathcal{P}}(Z,0)d\nu_{TX}(Z)=0,\eqno(4.115)$$
$$\int_{{\mathbb{R}}^{2n}}(\overline{z}_\alpha \overline{z}_\beta b^+_ib_lz_{\alpha'}\overline{z}'_r{\mathcal{P}})(0,Z)
{\mathcal{P}}(Z,0)d\nu_{TX}(Z)=0,\eqno(4.116)$$
$$(\overline{z}_\alpha \overline{z}_\beta b^+_ib_lz_{\alpha'}z_{j'}{\mathcal{P}})(0,Z)
=\frac{\delta_{li}}{\pi^2}(4\delta_{j'\alpha}\delta_{\alpha'\beta}+2\delta_{j'\alpha}b_\beta z_{\alpha'}+b_\alpha b_\beta z_{\alpha'}z_{j'}$$
$$+
2\delta_{\alpha'\beta}b_\alpha z_{j'}+2\delta_{\alpha\alpha'}b_\beta z_{j'}+2\delta_{j'\beta}b_\alpha z_{\alpha'}+4\delta_{j'\beta}\delta_{\alpha\alpha'}).\eqno(4.117)$$
By (4.98), (4.100) and (4.112)-(4.117), we have
\begin{eqnarray*}
&&\int_{{\mathbb{R}}^{2n}}\frac{1}{(4\pi)^{{j-\frac{q}{2}}-m_1+1}}B^{1,{j-\frac{q}{2}}-m_1-m_2}_0
B^{{j-\frac{q}{2}}-m_1-m_2,j-m_1}_1({\mathcal{L}}^{-1}_0{\mathcal{O}}^{+2}_2)^{m_1}\\
&&\cdot{\mathcal{L}}^{-1}_0
( {\mathcal{R}}^\bot)^{{j-\frac{q}{2}}-m_1}
(A_1A_3{\mathcal{P}})(0,Z){\mathcal{P}}(Z,0)d\nu_{TX}(Z)\\
&=&-\frac{8}{3}\frac{1}{(4\pi)^{{j-\frac{q}{2}}+2}}B^{1,{j-\frac{q}{2}}-m_1-m_2}_0B^{{j-\frac{q}{2}}-m_1-m_2,{j-\frac{q}{2}}-m_1}_1\\
&&
\cdot(B^{{j-\frac{q}{2}}-m_1,{j-\frac{q}{2}}}_0-B^{{j-\frac{q}{2}}-m_1,{j-\frac{q}{2}}}_1)( {\mathcal{R}}^\bot)^{{j-\frac{q}{2}}}\\
&&\times
(\partial_{{z}_\alpha }R^L)_{x_0}(\partial_{\overline{z}_\alpha},\partial_{{z}_{\alpha'}})
\left[\frac{2}{3\pi}(\partial_{\overline{z}_r }R^L)_{x_0}(\partial_{z_r},\partial_{\overline{z}_{\alpha'}})\right.\\
&&\left.
+\frac{2}{3\pi}(\partial_{\overline{z}_{\alpha'} }R^L)_{x_0}(\partial_{z_r},\partial_{\overline{z}_r})-
\frac{1}{6\pi}(\partial_{j'}R^L)_{x_0}(\partial_{\overline{z}_{\alpha'}},e_{j'})\right]\\
&&-\frac{16}{9\pi^2}\frac{1}{(4\pi)^{{j-\frac{q}{2}}+2}}B^{1,{j-\frac{q}{2}}-m_1-m_2}_0B^{{j-\frac{q}{2}}-m_1-m_2,{j-\frac{q}{2}}-m_1}_1\\
&&\cdot
(B^{{j-\frac{q}{2}}-m_1,{j-\frac{q}{2}}}_2-B^{j-m_1,{j-\frac{q}{2}}}_1+B^{{j-\frac{q}{2}}-m_1,{j-\frac{q}{2}}}_0)
( {\mathcal{R}}^\bot)^{{j-\frac{q}{2}}}\\
&&\times[(\partial_{\overline{z}_\beta }R^L)_{x_0}(\partial_{\overline{z}_\alpha},\partial_{{z}_{l}})
(\partial_{{z}_\beta }R^L)_{x_0}(\partial_{{z}_\alpha},\partial_{\overline{z}_{l}})\\
&&+
(\partial_{\overline{z}_\beta }R^L)_{x_0}(\partial_{\overline{z}_\alpha},\partial_{{z}_{l}})
(\partial_{{z}_\alpha }R^L)_{x_0}(\partial_{{z}_\beta},\partial_{\overline{z}_{l}})].~~~~~~~~~~~~~~~~~~~~~~~~~~~~~~~~~~~~~~~~~(4.118)
\end{eqnarray*}
In the following, we write $T=0$ if
$$\int_{{\mathbb{R}}^{2n}}T(0,Z){\mathcal{P}}(Z,0)d\nu_{TX}(Z)=0.\eqno(4.119)$$
We know that
$$(\overline{z}'_rb_{\alpha'}{\mathcal{P}})(0,Z)=(\overline{z}'_rb_lz_{\alpha'}\overline{z}'_{r'}{\mathcal{P}})(0,Z)=
(\overline{z}'_rb_lz_{\alpha}z_{j'}{\mathcal{P}})(0,Z)$$
$$=(z_{j'}b_lz_{\alpha'}\overline{z}'_{r'}{\mathcal{P}})(0,Z)=
(z_{j'}b_lz_{\alpha}z_{j''}{\mathcal{P}})(0,Z)=0,\eqno(4.120)$$
$$(z_{j'}b_{\alpha'}{\mathcal{P}})(0,Z)=[(b_{\alpha'}z_{j'}+2\delta_{j'\alpha'}){\mathcal{P}}](0,Z),\eqno(4.121)$$
so
\begin{eqnarray*}
&&\int_{{\mathbb{R}}^{2n}}\frac{1}{(4\pi)^{{j-\frac{q}{2}}-m_1+1}}B^{1,{j-\frac{q}{2}}-m_1-m_2}_0
B^{{j-\frac{q}{2}}-m_1-m_2,{j-\frac{q}{2}}-m_1}_1({\mathcal{L}}^{-1}_0{\mathcal{O}}^{+2}_2)^{m_1}~~~~~~
~~(4.122)\\
&&\cdot{\mathcal{L}}^{-1}_0
( {\mathcal{R}}^\bot)^{{j-\frac{q}{2}}-m_1}
(A_2A_3{\mathcal{P}})(0,Z){\mathcal{P}}(Z,0)d\nu_{TX}(Z)\\
&&=\frac{2}{(4\pi)^{{j-\frac{q}{2}}+2}}B^{1,{j-\frac{q}{2}}-m_1-m_2}_0B^{{j-\frac{q}{2}}-m_1-m_2,{j-\frac{q}{2}}-m_1}_1\\
&&\cdot(B^{{j-\frac{q}{2}}-m_1,{j-\frac{q}{2}}}_0-B^{{j-\frac{q}{2}}-m_1,{j-\frac{q}{2}}}_1)( {\mathcal{R}}^\bot)^{{j-\frac{q}{2}}}\\
&&\times[\frac{4}{3}(\partial_{{z}_{j'} }R^L)_{x_0}(\partial_{z_\alpha},\partial_{\overline{z}_{\alpha}})+
\frac{4}{3}(\partial_{{z}_{\alpha} }R^L)_{x_0}(\partial_{z_{j'}},\partial_{\overline{z}_\alpha})-
\frac{1}{3}(\partial_{\alpha}R^L)_{x_0}(\partial_{z_{j'}},e_{\alpha})]\\
&&\times [\frac{2}{3\pi}(\partial_{\overline{z}_r }R^L)_{x_0}(\partial_{z_r},\partial_{\overline{z}_{j'}})
+\frac{2}{3\pi}(\partial_{\overline{z}_{j'} }R^L)_{x_0}(\partial_{z_r},\partial_{\overline{z}_r})-
\frac{1}{6\pi}(\partial_{j''}R^L)_{x_0}(\partial_{\overline{z}_{j'}},e_{j''})].
\end{eqnarray*}
It holds that
$$b_lb_{l'}{\mathcal{P}}=b_lb_{l'}z_\alpha \overline{z}'_r{\mathcal{P}}=b_lz_\alpha\overline{z}'_rb_{l'}{\mathcal{P}}=
b_lz_\alpha\overline{z}'_rb_{l'}z_{\alpha'}\overline{z}'_{r'}{\mathcal{P}}$$
$$=b_lz_\alpha\overline{z}'_rb_{l'}z_{\alpha'}z_{j'}{\mathcal{P}}=
b_lz_\alpha{z}_jb_{l'}z_{\alpha'}\overline{z}'_{r'}{\mathcal{P}}=
b_lz_\alpha{z}_jb_{l'}z_{\alpha'}z_{j'}{\mathcal{P}}=0,\eqno(4.123)$$
$$(b_lz_\alpha z_{j'}b_{l'}{\mathcal{P}})(0,Z)=[(2\delta_{j'l'}b_lz_\alpha+2\delta_{\alpha l'}b_lz_{j'}+b_lb_{l'}z_\alpha z_{j'}){\mathcal{P}}](0,Z).\eqno(4.124)$$
By (4.49), (4.123) and (4.124), we have
\begin{eqnarray*}
&&\int_{{\mathbb{R}}^{2n}}\frac{1}{(4\pi)^{{j-\frac{q}{2}}-m_1+1}}B^{1,{j-\frac{q}{2}}-m_1-m_2}_0
B^{{j-\frac{q}{2}}-m_1-m_2,{j-\frac{q}{2}}-m_1}_1({\mathcal{L}}^{-1}_0{\mathcal{O}}^{+2}_2)^{m_1}~~~~~~
\\
&&\cdot{\mathcal{L}}^{-1}_0
( {\mathcal{R}}^\bot)^{{j-\frac{q}{2}}-m_1}
(A^2_3{\mathcal{P}})(0,Z){\mathcal{P}}(Z,0)d\nu_{TX}(Z)\\
&&=\frac{8}{3}\frac{1}{(4\pi)^{{j-\frac{q}{2}}+2}}B^{1,{j-\frac{q}{2}}-m_1-m_2}_0B^{{j-\frac{q}{2}}-m_1-m_2,j-m_1}_1\\
&&(2B^{{j-\frac{q}{2}}-m_1,{j-\frac{q}{2}}}_2-B^{{j-\frac{q}{2}}-m_1,{j-\frac{q}{2}}}_1)( {\mathcal{R}}^\bot)^{{j-\frac{q}{2}}}\\
&&\times\left\{(\partial_{{z}_{j'} }R^L)_{x_0}(\partial_{z_\alpha},\partial_{\overline{z}_{j'}})
[\frac{2}{3\pi}(\partial_{\overline{z}_r }R^L)_{x_0}(\partial_{z_r},\partial_{\overline{z}_{\alpha}})\right.\\
&&
+\frac{2}{3\pi}(\partial_{\overline{z}_{\alpha} }R^L)_{x_0}(\partial_{z_r},\partial_{\overline{z}_r})-
\frac{1}{6\pi}(\partial_{j'}R^L)_{x_0}(\partial_{\overline{z}_{\alpha}},e_{j'})]\\
&&+(\partial_{{z}_{j'} }R^L)_{x_0}(\partial_{z_\alpha},\partial_{\overline{z}_{\alpha}})
[\frac{2}{3\pi}(\partial_{\overline{z}_r }R^L)_{x_0}(\partial_{z_r},\partial_{\overline{z}_{j'}})\\
&&\left.
+\frac{2}{3\pi}(\partial_{\overline{z}_{j'} }R^L)_{x_0}(\partial_{z_r},\partial_{\overline{z}_r})-
\frac{1}{6\pi}(\partial_{j''}R^L)_{x_0}(\partial_{\overline{z}_{j'}},e_{j''})]\right\}~~~~~~~~~~~~~~~~~~~~~~~~~~(4.125)\\
\end{eqnarray*}
We know that
$$b_lb_rz_\alpha b_{l'}{\mathcal{P}}=b_lb_rz_\alpha b_{l'}z_{\alpha'}\overline{z}'_{r'}{\mathcal{P}}=0,\eqno(4.126)$$
$$b_lb_rz_\alpha b_{l'}z_{\alpha'}{z}_{j'}{\mathcal{P}}=(2\delta_{\alpha l'}b_lb_rz_{\alpha'}z_{j'}+b_lb_rb_{l'}z_\alpha z_{\alpha'} z_{j'}){\mathcal{P}}.
\eqno(4.127)$$
$$(b_lb_rb_{l'}z_\alpha z_{\alpha'} z_{j'}{\mathcal{P}})(0,Z)=-8
[\delta_{\alpha l'}(\delta_{l\alpha'}\delta_{j'r}+\delta_{r\alpha'}\delta_{lj'})$$
$$+
\delta_{\alpha r}(\delta_{lj'}\delta_{l'\alpha'}+\delta_{l\alpha'}\delta_{j'l'})+
\delta_{\alpha l}(\delta_{r\alpha'}\delta_{l'j'}+\delta_{rj'}\delta_{l'\alpha'})].\eqno(4.128)$$
By (4.126)-(4.128), we get
\begin{eqnarray*}
&&\int_{{\mathbb{R}}^{2n}}\frac{1}{(4\pi)^{{j-\frac{q}{2}}-m_1+1}}B^{1,{j-\frac{q}{2}}-m_1-m_2}_0
B^{{j-\frac{q}{2}}-m_1-m_2,{j-\frac{q}{2}}-m_1}_1({\mathcal{L}}^{-1}_0{\mathcal{O}}^{+2}_2)^{m_1}~~~~~~
\\
&&\cdot{\mathcal{L}}^{-1}_0
( {\mathcal{R}}^\bot)^{{j-\frac{q}{2}}-m_1}
(A_4A_3{\mathcal{P}})(0,Z){\mathcal{P}}(Z,0)d\nu_{TX}(Z)\\
&&=\frac{16}{9\pi}\frac{1}{(4\pi)^{{j-\frac{q}{2}}+2}}B^{1,{j-\frac{q}{2}}-m_1-m_2}_0B^{{j-\frac{q}{2}}-m_1-m_2,{j-\frac{q}{2}}-m_1}_1
( {\mathcal{R}}^\bot)^{{j-\frac{q}{2}}}\\
&&\times\left\{(B^{{j-\frac{q}{2}}-m_1,{j-\frac{q}{2}}}_2-B^{{j-\frac{q}{2}}-m_1,{j-\frac{q}{2}}}_3)\left[
(\partial_{\overline{z}_r }R^L)_{x_0}(\partial_{z_\alpha},\partial_{\overline{z}_{l}})
(\partial_{{z}_r }R^L)_{x_0}(\partial_{z_l},\partial_{\overline{z}_{\alpha}})\right.\right.\\
&&\left.
+(\partial_{\overline{z}_r }R^L)_{x_0}(\partial_{z_\alpha},\partial_{\overline{z}_{l}})
(\partial_{{z}_l }R^L)_{x_0}(\partial_{z_r},\partial_{\overline{z}_{\alpha}})\right]\\
&&-B^{j-m_1,j}_3\left[
(\partial_{\overline{z}_\alpha }R^L)_{x_0}(\partial_{z_\alpha},\partial_{\overline{z}_{j'}})
(\partial_{{z}_{j'} }R^L)_{x_0}(\partial_{z_{l'}},\partial_{\overline{z}_{l'}})\right.\\
&&
+(\partial_{\overline{z}_\alpha }R^L)_{x_0}(\partial_{z_\alpha},\partial_{\overline{z}_{l}})
(\partial_{{z}_{l'} }R^L)_{x_0}(\partial_{z_l},\partial_{\overline{z}_{l'}})\\
&&
+(\partial_{\overline{z}_r }R^L)_{x_0}(\partial_{z_l},\partial_{\overline{z}_{l}})
(\partial_{{z}_{l'} }R^L)_{x_0}(\partial_{z_r},\partial_{\overline{z}_{l'}})\\
&&\left.\left.+(\partial_{\overline{z}_r }R^L)_{x_0}(\partial_{z_l},\partial_{\overline{z}_{l}})
(\partial_{{z}_r }R^L)_{x_0}(\partial_{z_{\alpha'}},\partial_{\overline{z}_{\alpha'}})\right]\right\}
~~~~~~~~~~~~~~~~~~~~~~~~~~~~~~~~~~~~~~~(4.129)
\end{eqnarray*}
It holds that
$$(\overline{z}_\alpha z_\beta b^+_ib_lb_rz_{\alpha'}{\mathcal{P}})(0,Z)=
\frac{2}{\pi}[\delta_{ri}(b_lb_\alpha z_\beta  z_{\alpha'}+2\delta_{\alpha\beta}b_lz_{\alpha'}+2\delta_{\alpha\alpha'}b_lz_\beta+
2\delta_{l\beta}b_\alpha z_{\alpha'}+4\delta_{l\beta}\delta_{\alpha\alpha'})$$
$$+\delta_{li}(b_rb_\alpha z_\beta  z_{\alpha'}+2\delta_{\alpha\beta}b_rz_{\alpha'}+2\delta_{\alpha\alpha'}b_rz_\beta+
2\delta_{r\beta}b_\alpha z_{\alpha'}+4\delta_{r\beta}\delta_{\alpha\alpha'})],\eqno(4.130)$$
$$\overline{z}_\alpha \overline{z}_{r}b^+_ib_lb_{r'}z_{\alpha'}{\mathcal{P}}=0.\eqno(4.131)$$
By (4.130) and (4.131), we have
\begin{eqnarray*}
&&\int_{{\mathbb{R}}^{2n}}\frac{1}{(4\pi)^{{j-\frac{q}{2}}-m_1+1}}B^{1,{j-\frac{q}{2}}-m_1-m_2}_0
B^{{j-\frac{q}{2}}-m_1-m_2,{j-\frac{q}{2}}-m_1}_2({\mathcal{L}}^{-1}_0{\mathcal{O}}^{+2}_2)^{m_1}~~~~~~
\\
&&\cdot{\mathcal{L}}^{-1}_0
( {\mathcal{R}}^\bot)^{{j-\frac{q}{2}}-m_1}
(A_1A_4{\mathcal{P}})(0,Z){\mathcal{P}}(Z,0)d\nu_{TX}(Z)\\
&&=-\frac{16}{9\pi^2}\frac{1}{(4\pi)^{{j-\frac{q}{2}}+2}}B^{1,{j-\frac{q}{2}}-m_1-m_2}_0B^{{j-\frac{q}{2}}-m_1-m_2,j-m_1}_2
( {\mathcal{R}}^\bot)^{{j-\frac{q}{2}}}\\
&&\left[(B^{{j-\frac{q}{2}}-m_1,{j-\frac{q}{2}}}_2-B^{{j-\frac{q}{2}}-m_1,{j-\frac{q}{2}}}_1)
(\partial_{{z}_\alpha }R^L)_{x_0}(\partial_{\overline{z}_\alpha},\partial_{{z}_{r}})
(\partial_{\overline{z}_{r} }R^L)_{x_0}(\partial_{z_{l}},\partial_{\overline{z}_{l}})\right.\\
&&+
(B^{j-m_1,j}_2+B^{j-m_1,j}_0-B^{j-m_1,j}_1)
(\partial_{{z}_l }R^L)_{x_0}(\partial_{\overline{z}_\alpha},\partial_{{z}_{r}})
(\partial_{\overline{z}_{r} }R^L)_{x_0}(\partial_{z_{\alpha}},\partial_{\overline{z}_{l}})\\
&&+
(B^{j-m_1,j}_2+B^{j-m_1,j}_0-B^{j-m_1,j}_1)
(\partial_{{z}_r }R^L)_{x_0}(\partial_{\overline{z}_\alpha},\partial_{{z}_{l}})
(\partial_{\overline{z}_{r} }R^L)_{x_0}(\partial_{z_{\alpha}},\partial_{\overline{z}_{l}})\\
&&+
(B^{j-m_1,j}_2-B^{j-m_1,j}_1)
(\partial_{{z}_\alpha }R^L)_{x_0}(\partial_{\overline{z}_\alpha},\partial_{{z}_{l}})
(\partial_{\overline{z}_{r} }R^L)_{x_0}(\partial_{z_{r}},\partial_{\overline{z}_{l}})\\
&&-
B^{j-m_1,j}_1
(\partial_{{z}_l }R^L)_{x_0}(\partial_{\overline{z}_\alpha},\partial_{{z}_{r}})
(\partial_{\overline{z}_{r} }R^L)_{x_0}(\partial_{z_{\alpha}},\partial_{\overline{z}_{l}})\\
&&\left.-
B^{j-m_1,j}_1
(\partial_{{z}_r }R^L)_{x_0}(\partial_{\overline{z}_\alpha},\partial_{{z}_{l}})
(\partial_{\overline{z}_{r} }R^L)_{x_0}(\partial_{z_{\alpha}},\partial_{\overline{z}_{l}})\right].
~~~~~~~~~~~~~~~~~~~~~~~~~~~~~~~(4.132)
\end{eqnarray*}
It holds that
$$\overline{z}'_{r'}b_lb_rz_\alpha {\mathcal{P}}=0,\eqno(4.133)$$
$$z_{j'}b_lb_rz_\alpha {\mathcal{P}}=(b_lb_rz_{j'}z_\alpha+2\delta_{jr'}b_lz_\alpha+2\delta_{j'l}b_rz_\alpha){\mathcal{P}}.\eqno(4.134)$$
By (4.133) and (4.134), we get
\begin{eqnarray*}
&&\int_{{\mathbb{R}}^{2n}}\frac{1}{(4\pi)^{{j-\frac{q}{2}}-m_1+1}}B^{1,{j-\frac{q}{2}}-m_1-m_2}_0B^{{j-\frac{q}{2}}-m_1-m_2,
{j-\frac{q}{2}}-m_1}_2({\mathcal{L}}^{-1}_0{\mathcal{O}}^{+2}_2)^{m_1}
~~~~~~~~(4.135)
\\
&&\cdot{\mathcal{L}}^{-1}_0
( {\mathcal{R}}^\bot)^{{j-\frac{q}{2}}-m_1}
(A_2A_4{\mathcal{P}})(0,Z){\mathcal{P}}(Z,0)d\nu_{TX}(Z)\\
&&=\frac{4}{3\pi}\frac{1}{(4\pi)^{{j-\frac{q}{2}}+2}}B^{1,{j-\frac{q}{2}}-m_1-m_2}_0B^{{j-\frac{q}{2}}-m_1-m_2,{j-\frac{q}{2}}-m_1}_2\\
&&\cdot(B^{{j-\frac{q}{2}}-m_1,{j-\frac{q}{2}}}_2-B^{{j-\frac{q}{2}}-m_1,{j-\frac{q}{2}}}_1)( {\mathcal{R}}^\bot)^{{j-\frac{q}{2}}}\\
&&\left\{
(\partial_{\overline{z}_\alpha }R^L)_{x_0}(\partial_{{z}_\alpha},\partial_{\overline{z}_{l}})
[\frac{4}{3}(\partial_{{z}_{l} }R^L)_{x_0}(\partial_{z_\alpha},\partial_{\overline{z}_\alpha})+
\frac{4}{3}(\partial_{{z}_{\alpha} }R^L)_{x_0}(\partial_{z_l},\partial_{\overline{z}_\alpha})\right.\\
&&-
\frac{1}{3}(\partial_{\alpha}R^L)_{x_0}(\partial_{z_{l}},e_{\alpha})]+
(\partial_{\overline{z}_r}R^L)_{x_0}(\partial_{{z}_l},\partial_{\overline{z}_{l}})\\
&&\left.\times
[\frac{4}{3}(\partial_{{z}_r }R^L)_{x_0}(\partial_{z_\alpha},\partial_{\overline{z}_\alpha})+
\frac{4}{3}(\partial_{{z}_{\alpha} }R^L)_{x_0}(\partial_{z_{r}},\partial_{\overline{z}_\alpha})-
\frac{1}{3}(\partial_{\alpha}R^L)_{x_0}(\partial_{z_{r}},e_{\alpha})]\right\}
\end{eqnarray*}
We know that
$$b_{l'}b_lb_rz_\alpha{\mathcal{P}}=b_{l'}z_\alpha\overline{z}'_{r'} b_lb_rz_\alpha{\mathcal{P}}=b_{l'}b_{r'}z_{\alpha'}b_lb_rz_\alpha{\mathcal{P}}=0,\eqno(4.136)$$
$$b_{l'}z_{\alpha'}z_{j'}b_lb_rz_\alpha{\mathcal{P}}=(b_{l'}b_lb_rz_{\alpha'}z_{j'}z_\alpha+2\delta_{\alpha'r}b_{l'}b_lz_{j'}z_\alpha+
2\delta_{j'r}b_{l'}b_lz_{\alpha'}z_\alpha$$
$$+2\delta_{\alpha'l}b_{l'}b_rz_{j'}z_\alpha
+2\delta_{j'l}b_{l'}b_rz_{\alpha'}z_\alpha
+4\delta_{\alpha'l}\delta_{j'r}b_{l'}z_\alpha+4\delta_{j'l}\delta_{r\alpha'}b_{l'}z_\alpha){\mathcal{P}}.\eqno(4.137)$$

By (4.136) and (4.137), we get
\begin{eqnarray*}
&&\int_{{\mathbb{R}}^{2n}}\frac{1}{(4\pi)^{{j-\frac{q}{2}}-m_1+1}}B^{1,{j-\frac{q}{2}}-m_1-m_2}_0
B^{{j-\frac{q}{2}}-m_1-m_2,{j-\frac{q}{2}}-m_1}_2({\mathcal{L}}^{-1}_0{\mathcal{O}}^{+2}_2)^{m_1}
\\
&&\cdot{\mathcal{L}}^{-1}_0
( {\mathcal{R}}^\bot)^{{j-\frac{q}{2}}-m_1}
(A_3A_4{\mathcal{P}})(0,Z){\mathcal{P}}(Z,0)d\nu_{TX}(Z)\\
&=&\frac{16}{9\pi}\frac{1}{(4\pi)^{{j-\frac{q}{2}}+2}}B^{1,{j-\frac{q}{2}}-m_1-m_2}_0B^{{j-\frac{q}{2}}-m_1-m_2,{j-\frac{q}{2}}-m_1}_2
( {\mathcal{R}}^\bot)^{{j-\frac{q}{2}}}\\
&&
\cdot\left\{(-B^{{j-\frac{q}{2}}-m_1,{j-\frac{q}{2}}}_3+2B^{{j-\frac{q}{2}}-m_1,{j-\frac{q}{2}}}_2-B^{{j-\frac{q}{2}}-m_1,{j-\frac{q}{2}}}_1)\right.\\
&&
\times[(\partial_{{z}_r }R^L)_{x_0}(\partial_{{z}_l},\partial_{\overline{z}_{\alpha}})
(\partial_{\overline{z}_r }R^L)_{x_0}(\partial_{{z}_\alpha},\partial_{\overline{z}_{l}})+
(\partial_{{z}_l }R^L)_{x_0}(\partial_{{z}_r},\partial_{\overline{z}_{\alpha}})
(\partial_{\overline{z}_r }R^L)_{x_0}(\partial_{{z}_\alpha},\partial_{\overline{z}_{l}})]\\
&&+
(-B^{{j-\frac{q}{2}}-m_1,{j-\frac{q}{2}}}_3+
B^{{j-\frac{q}{2}}-m_1,{j-\frac{q}{2}}}_2)
[(\partial_{{z}_l }R^L)_{x_0}(\partial_{{z}_{l'}},\partial_{\overline{z}_{l'}})
(\partial_{\overline{z}_\alpha}R^L)_{x_0}(\partial_{{z}_\alpha},\partial_{\overline{z}_{l}})\\
&&+
(\partial_{{z}_{l'} }R^L)_{x_0}(\partial_{{z}_l},\partial_{\overline{z}_{l'}})
(\partial_{\overline{z}_\alpha }R^L)_{x_0}(\partial_{{z}_\alpha},\partial_{\overline{z}_{l}})
+
(\partial_{{z}_{l'} }R^L)_{x_0}(\partial_{z_r},\partial_{\overline{z}_{l'}})
(\partial_{\overline{z}_r }R^L)_{x_0}(\partial_{z_l},\partial_{\overline{z}_{l}})\\
&&\left.+
(\partial_{{z}_{r} }R^L)_{x_0}(\partial_{{z}_l},\partial_{\overline{z}_{l'}})
(\partial_{\overline{z}_r}R^L)_{x_0}(\partial_{{z}_{l}},\partial_{\overline{z}_{l}})]\right\}.~~~~~~~~~~~~~~~~~~~~~~~~~~~~~~~~~~~~~~~(4.138)
\end{eqnarray*}

$$\int_{{\mathbb{R}}^{2n}}\frac{1}{(4\pi)^{{j-\frac{q}{2}}-m_1+1}}B^{1,{j-\frac{q}{2}}-m_1-m_2}_0
B^{{j-\frac{q}{2}}-m_1-m_2,{j-\frac{q}{2}}-m_1}_2({\mathcal{L}}^{-1}_0{\mathcal{O}}^{+2}_2)^{m_1}
$$
$$\times{\mathcal{L}}^{-1}_0
( {\mathcal{R}}^\bot)^{{j-\frac{q}{2}}-m_1}
(A^2_4{\mathcal{P}})(0,Z){\mathcal{P}}(Z,0)d\nu_{TX}(Z)=0.\eqno(4,139)$$
By the basic assumption, then
 $$(\nabla^X_UR^L)(V,W)=\left<(\nabla^X_U\widetilde{{ J}})V,W\right>.\eqno(4.140)$$
 Direct computations show that
 $$(b_lz_\alpha+2\delta_{\alpha l})
(\overline{z}_r-\frac{b_r}{2\pi}-\overline{z}'_r)z_{j'}{\mathcal{P}}=\frac{\delta_{rj'}}{\pi}(b_lz_\alpha
+2\delta_{\alpha l}){\mathcal{P}};\eqno(4.141)$$
$$(\overline{z}_r-\frac{b_r}{2\pi}-\overline{z}'_r)z_{j'}{\mathcal{P}}=\frac{\delta_{rj'}}{\pi}{\mathcal{P}};\eqno(4.142)$$
$$(b_lz_\alpha+2\delta_{\alpha l})
(\overline{z}_r-\frac{b_r}{2\pi}-\overline{z}'_r)\overline{z}_{j'}{\mathcal{P}}=
 (\overline{z}_r-\frac{b_r}{2\pi}-\overline{z}'_r)\overline{z}_{j'}{\mathcal{P}}=0;\eqno(4.143)$$
 $$(b_lz_\alpha+2\delta_{\alpha l})
(\overline{z}_r-\frac{b_r}{2\pi}-\overline{z}'_r)b_{j'}{\mathcal{P}}=
 (\overline{z}_r-\frac{b_r}{2\pi}-\overline{z}'_r)b_{j'}{\mathcal{P}}=0;\eqno(4.144)$$
 $$(b_lz_\alpha+2\delta_{\alpha l})
(\overline{z}_r-\frac{b_r}{2\pi}-\overline{z}'_r)b_{l'}z_{\alpha'}\overline{z}'_{r'}{\mathcal{P}}=
 [\frac{b_lb_{l'}z_\alpha\overline{z}'_{r'}}{\pi}\delta_{rr'}+\frac{2\delta_{\alpha l'}\delta_{rr'}}{\pi}b_l\overline{z}'_{r'}
 +\frac{2\delta_{\alpha l}\delta_{rr'}}{\pi}b_{l'}\overline{z}'_{r'}]{\mathcal{P}};\eqno(4.145)$$
 $$(\overline{z}_r-\frac{b_r}{2\pi}-\overline{z}'_r)b_{l'}z_{\alpha'}\overline{z}'_{r'}{\mathcal{P}}=\frac{\delta_{rr'}}{\pi}
 b_{l'}\overline{z}'_{r'}{\mathcal{P}};\eqno(4.146)$$
 $$(b_lz_\alpha+2\delta_{\alpha l})
(\overline{z}_r-\frac{b_r}{2\pi}-\overline{z}'_r)b_{l'}z_{\alpha'}{z}_{j'}{\mathcal{P}}=$$$$
[ \frac{b_lb_{l'}z_\alpha z_{\alpha'}}{\pi}\delta_{rj'}
 +\frac{2\delta_{\alpha l'}\delta_{rj'}}{\pi}b_l{z}_{\alpha'}
 + \frac{b_lb_{l'}z_\alpha z_{j'}}{\pi}\delta_{r\alpha'}$$
  $$
 +\frac{2\delta_{\alpha l'}\delta_{r\alpha'}}{\pi}b_l{z}_{j'}
  +\frac{2\delta_{\alpha l}\delta_{rj'}}{\pi}b_{l'}{z}_{\alpha'}
 +\frac{2\delta_{\alpha l}\delta_{r\alpha'}}{\pi}b_{l'}{z}_{j'}]{\mathcal{P}};\eqno(4.147)$$
$$ (\overline{z}_r-\frac{b_r}{2\pi}-\overline{z}'_r)b_{l'}z_{\alpha'}{z}_{j'}{\mathcal{P}}=
 [\frac{\delta_{r j'}}{\pi}b_{l'}{z}_{\alpha'}
 +\frac{\delta_{r\alpha'}}{\pi}b_{l'}{z}_{j'}]{\mathcal{P}};\eqno(4.148)$$
$$(b_lz_\alpha+2\delta_{\alpha l})(\overline{z}_r-\frac{b_r}{2\pi}-\overline{z}'_r)b_{l'}b_{r'}z_{\alpha'}{\mathcal{P}}=
 [\frac{\delta_{r\alpha'}}{\pi}b_lb_{l'}b_{r'}z_\alpha$$
 $$+\frac{2\delta_{\alpha r'}\delta_{r\alpha'}}{\pi}b_lb_{l'}
 +\frac{2\delta_{\alpha l'}\delta_{r\alpha'}}{\pi}b_lb_{r'}
 +\frac{2\delta_{\alpha l}\delta_{r\alpha'}}{\pi}b_{l'}b_{r'}]{\mathcal{P}};\eqno(4.149)$$
 $$(\overline{z}_r-\frac{b_r}{2\pi}-\overline{z}'_r)b_{l'}b_{r'}z_{\alpha'}{\mathcal{P}}=\frac{\delta_{r\alpha'}}{\pi}b_{l'}b_{r'}
 {\mathcal{P}}.\eqno(4.150)$$

By (4.96),(4.102) and (4.141)-(4.150), we get

\begin{eqnarray*}
 &&\sum_{0\leq m_1+m_2\leq {j-\frac{q}{2}}}I_{2j}({\mathcal{L}}^{-1}_0{\mathcal{O}}^{+2}_2)^{m_1}({\mathcal{L}}^{-1}_0
P^{N^\bot}A_5)({\mathcal{L}}^{-1}_0{\mathcal{O}}^{+2}_2)^{m_2}\\
&&\times({\mathcal{L}}^{-1}_0P^{N^\bot}{\mathcal{O}}'_1)
({\mathcal{L}}^{-1}_0{\mathcal{O}}^{+2}_2)^{{j-\frac{q}{2}}-m_1-m_2}P^N({\mathcal{O}}^{-2}_2{\mathcal{L}}^{-1}_0)^{j-\frac{q}{2}}I_{2j}
\\
&&=\sum_{0\leq m_1+m_2\leq {j-\frac{q}{2}}}\frac{1}{(4\pi)^{2j-q+2}}B_0^{1,j-\frac{q}{2}}B_0^{1,j-\frac{q}{2}-m_1-m_2}B_0^{j-\frac{q}{2}-m_1-m_2,j-\frac{q}{2}-m_1}
({{\mathcal{R}}^\bot})^{j-\frac{q}{2}}\\
&&\left\{(B_0^{j-\frac{q}{2}-m_1,j-\frac{q}{2}}-B_1^{j-\frac{q}{2}-m_1,j-\frac{q}{2}})
\frac{4}{3\pi}(\partial_{\overline{z}_r }R^L)_{x_0}(\partial_{z_\alpha},\partial_{\overline{z}_\alpha})\right.\\
&&\cdot[\frac{4}{3}(\partial_{{z}_{r} }R^L)_{x_0}(\partial_{z_{\alpha'}},\partial_{\overline{z}_{\alpha'}})+
\frac{4}{3}(\partial_{{z}_{\alpha'} }R^L)_{x_0}(\partial_{z_{r}},\partial_{\overline{z}_{\alpha'}})-
\frac{1}{3}(\partial_{\alpha'}R^L)_{x_0}(\partial_{z_{r}},e_{\alpha'})]\\
&&-\frac{1}{3\pi}B_0^{j-\frac{q}{2}-m_1,j-\frac{q}{2}}
(\partial_{j'}R^L)_{x_0}(\partial_{\overline{z}_r},e_{j'})[\frac{4}{3}(\partial_{{z}_{r} }R^L)_{x_0}(\partial_{z_{\alpha'}},\partial_{\overline{z}_{\alpha'}})\\
&&\left.+
\frac{4}{3}(\partial_{{z}_{\alpha'} }R^L)_{x_0}(\partial_{z_{r}},\partial_{\overline{z}_{\alpha'}})-
\frac{1}{3}(\partial_{\alpha'}R^L)_{x_0}(\partial_{z_{r}},e_{\alpha'})]\right\}I_{{\rm det}(\overline{W}^*)\otimes E}
({{\mathcal{R}}^{\bot,*}})^{j-\frac{q}{2}}\\
&&+\sum_{0\leq m_1+m_2\leq {j-\frac{q}{2}}}\frac{1}{(4\pi)^{2j-q+2}}B_0^{1,j-\frac{q}{2}}B_0^{1,j-\frac{q}{2}-m_1-m_2}B_1^{j-\frac{q}{2}-m_1-m_2,j-\frac{q}{2}-m_1}
({{\mathcal{R}}^\bot})^{j-\frac{q}{2}}\\
&&\left\{\frac{16}{9\pi}(B_2^{j-\frac{q}{2}-m_1,j-\frac{q}{2}}-B_1^{j-\frac{q}{2}-m_1,j-\frac{q}{2}})
[(\partial_{\overline{z}_r }R^L)_{x_0}(\partial_{z_\alpha},\partial_{\overline{z}_\alpha})
(\partial_{{z}_r }R^L)_{x_0}(\partial_{z_{\alpha'}},\partial_{\overline{z}_{\alpha'}})\right.\\
&&
+(\partial_{\overline{z}_r }R^L)_{x_0}(\partial_{z_\alpha},\partial_{\overline{z}_l})
(\partial_{{z}_r }R^L)_{x_0}(\partial_{z_{l}},\partial_{\overline{z}_{\alpha}})
+(\partial_{\overline{z}_r }R^L)_{x_0}(\partial_{z_\alpha},\partial_{\overline{z}_\alpha})
(\partial_{{z}_{l'} }R^L)_{x_0}(\partial_{z_{r}},\partial_{\overline{z}_{l'}})\\
&&
+(\partial_{\overline{z}_r }R^L)_{x_0}(\partial_{z_\alpha},\partial_{\overline{z}_l})
(\partial_{{z}_l }R^L)_{x_0}(\partial_{z_{r}},\partial_{\overline{z}_{\alpha}})]\\
&&+\frac{4}{9\pi}B_1^{j-\frac{q}{2}-m_1,j-\frac{q}{2}}[(\partial_{j'}R^L)_{x_0}(\partial_{\overline{z}_r},e_{j'})
(\partial_{{z}_r }R^L)_{x_0}(\partial_{z_{\alpha'}},\partial_{\overline{z}_{\alpha'}})\\
&&\left.
+(\partial_{j'}R^L)_{x_0}(\partial_{\overline{z}_r},e_{j'})
(\partial_{{z}_{l'} }R^L)_{x_0}(\partial_{z_{r}},\partial_{\overline{z}_{l'}})]\right\}I_{{\rm det}(\overline{W}^*)\otimes E}
({{\mathcal{R}}^{\bot,*}})^{j-\frac{q}{2}}.~~~~~~~~~~~~~~~~(4.151)
\end{eqnarray*}

Write $B_{l_1}B_{l_2}B_{l_3}$ for $B^{1,{j-\frac{q}{2}}-m_1-m_2}_{l_1}B^{{j-\frac{q}{2}}-m_1-m_2,{j-\frac{q}{2}}-m_1}_{l_2}B^{j-\frac{q}{2}-m_1,{j-\frac{q}{2}}}_{l_3}$.
By (4.96)-(4.101),(4.108)-(4.111),(4.118),(4.112),(4.125),(4.129),(4.132),(4.135),(4.138) and (4.139), we get

\begin{eqnarray*}
&&\sum_{0\leq m_1+m_2\leq {j-\frac{q}{2}}}I_{2j}({\mathcal{L}}^{-1}_0{\mathcal{O}}^{+2}_2)^{m_1}({\mathcal{L}}^{-1}_0
P^{N^\bot}(A_1+A_2+A_3+A_4))\\
&&({\mathcal{L}}^{-1}_0{\mathcal{O}}^{+2}_2)^{m_2}({\mathcal{L}}^{-1}_0P^{N^\bot}{\mathcal{O}}'_1)
({\mathcal{L}}^{-1}_0{\mathcal{O}}^{+2}_2)^{{j-\frac{q}{2}}-m_1-m_2}P^N({\mathcal{O}}^{-2}_2{\mathcal{L}}^{-1}_0)^{j-\frac{q}{2}}I_{2j}
\\
&=&\sum_{0\leq m_1+m_2\leq {j-\frac{q}{2}}}\frac{1}{2^{j-\frac{q}{2}}({j-\frac{q}{2}})!(4\pi)^{2j-q+2}}I_{2j}\left\{
\frac{1}{9\pi}(B_0B_0B_2-B_0B_0B_1)\right.\\
&&\times(2A_{\overline{l}ls}A_{sr\overline{r}}-2A_{\overline{l}ls}A_{r\overline{r}s}
+2A_{\overline{s}l\overline{l}}A_{sr\overline{r}}\\
&&-2A_{\overline{s}l\overline{l}}A_{r\overline{r}s}
-A_{\overline{l}l\overline{s}}A_{e_rse_r}-A_{\overline{s}l\overline{l}}A_{e_rse_r})\\
&&
+[-\frac{4}{9\pi}B_0B_1B_0-\frac{2}{3\pi}B_0B_1B_1-\frac{4}{9\pi}B_0B_1B_2\\
&&
+\frac{2}{9\pi}B_0B_1B_3
-(\frac{2}{9\pi^2}+\frac{4}{9\pi})B_0B_2B_2-(\frac{2}{9\pi^2}
+\frac{2}{9\pi})B_0B_2B_1\\
&&
+\frac{2}{9\pi}B_0B_2B_3](A_{l\overline{l}s}A_{\overline{r}r\overline{s}}
+A_{l\overline{l}s}A_{\overline{s}r\overline{r}})\\
&&
+[(-\frac{2}{9\pi^2}-\frac{2}{9\pi})B_0B_1B_2+\frac{2}{9\pi^2}B_0B_1B_1-\frac{2}{9\pi^2}B_0B_1B_0\\
&&+\frac{2}{9\pi}B_0B_1B_3
+(-\frac{2}{9\pi^2}-\frac{4}{9\pi})B_0B_2B_2-\frac{2}{9\pi^2}B_0B_2B_0\\
&&+(\frac{4}{9\pi^2}+\frac{2}{9\pi})B_0B_2B_1+\frac{2}{9\pi}B_0B_2B_3]
(A_{\overline{l}\overline{s}r}A_{ls\overline{r}}
+A_{\overline{l}\overline{s}r}A_{sl\overline{r}})\\
&&+\frac{1}{9\pi}(2B_0B_1B_0-4B_0B_1B_1+4B_0B_1B_2-2B_0B_1B_3\\
&&
-2B_0B_2B_1-2B_0B_2B_3+4B_0B_2B_2)(
A_{sl\overline{l}}A_{\overline{r}r\overline{{s}}}
+A_{sl\overline{l}}A_{\overline{s}r\overline{r}})\\
&&
+\frac{1}{9\pi}(-B_0B_1B_0+B_0B_1B_1-B_0B_2B_0\\
&&+B_0B_2B_1)(A_{e_rse_r}A_{\overline{r}r\overline{s}}+A_{e_rse_r}A_{\overline{s}r\overline{r}})\\
&&+\frac{1}{9\pi}(2B_0B_1B_0-3B_0B_1B_1+2B_0B_1B_2-B_0B_0B_1)A_{e_r\overline{s}e_r}A_{l\overline{l}{s}}\\
&&+\frac{1}{9\pi}(-B_0B_1B_0+2B_0B_1B_1+2B_0B_1B_2+B_0B_0B_1)A_{e_r\overline{s}e_r}A_{sl\overline{l}}\\
&&\left.+\frac{1}{18\pi}(B_0B_1B_0-B_0B_1B_1-B_0B_0B_1)A_{e_rse_r}A_{e_l\overline{s}e_l}\right\}\\
&&
\cdot({\mathcal{R}}^\bot)(x)^{{j-\frac{q}{2}}}I_{{\rm det}(\overline{W}^*)\otimes E}({\mathcal{R}}^{\bot,*})(x)^{{j-\frac{q}{2}}}I_{2j}
.~~~~~~~~~~~~~~~~~~~~~~~~~~~~~~~~~~~~~~~~(4.152)
\end{eqnarray*}
By (4.97), then
$$V^0_a=(4.151)+(4.152).\eqno(4.153)$$
Set
$$V^1_a=\sum_{0\leq m_1+m_2\leq {j-\frac{q}{2}}}I_{2j}({\mathcal{L}}^{-1}_0{\mathcal{O}}^{+2}_2)^{m_1}({\mathcal{L}}^{-1}_0
P^{N^\bot}{\mathcal{O}}''_1)
({\mathcal{L}}^{-1}_0{\mathcal{O}}^{+2}_2)^{m_2}$$
$$\cdot({\mathcal{L}}^{-1}_0P^{N^\bot}{\mathcal{O}}''_1)({\mathcal{L}}^{-1}_0{\mathcal{O}}^{+2}_2)^{{j-\frac{q}{2}}-m_1-m_2}
P^N({\mathcal{O}}^{-2}_2{\mathcal{L}}^{-1}_0)^{j-\frac{q}{2}}I_{2j}.
\eqno(4.154)$$
Set $n=4$, $j=2$, $q=2$. By (4.89) and $({\mathcal{L}}^{-1}_0{\mathcal{O}}^{+2}_2)({\mathcal{L}}^{-1}_0P^{N^\bot}{\mathcal{O}}''_1)P^N$
and $({\mathcal{L}}^{-1}_0{\mathcal{O}}^{+2}_2)
P^N$ being $(0, (2,2))$ forms, then
\begin{eqnarray*}
V^1_a&=&I_{4}[({\mathcal{L}}^{-1}_0{\mathcal{O}}^{+2}_2)({\mathcal{L}}^{-1}_0
P^{N^\bot}{\mathcal{O}}''_1)
({\mathcal{L}}^{-1}_0P^{N^\bot}{\mathcal{O}}''_1)\\
&&
+({\mathcal{L}}^{-1}_0
P^{N^\bot}{\mathcal{O}}''_1)
({\mathcal{L}}^{-1}_0{\mathcal{O}}^{+2}_2)({\mathcal{L}}^{-1}_0P^{N^\bot}{\mathcal{O}}''_1)\\
&&+({\mathcal{L}}^{-1}_0
P^{N^\bot}{\mathcal{O}}''_1)
({\mathcal{L}}^{-1}_0P^{N^\bot}{\mathcal{O}}''_1)({\mathcal{L}}^{-1}_0{\mathcal{O}}^{+2}_2)]
P^N({\mathcal{O}}^{-2}_2{\mathcal{L}}^{-1}_0)I_{4}\\
&
=&I_{4}({\mathcal{L}}^{-1}_0{\mathcal{O}}^{+2}_2)({\mathcal{L}}^{-1}_0
P^{N^\bot}{\mathcal{O}}''_1)
({\mathcal{L}}^{-1}_0P^{N^\bot}{\mathcal{O}}''_1)P^N({\mathcal{O}}^{-2}_2{\mathcal{L}}^{-1}_0)I_{4}.
~~~~~~~~~~~~~(4.155)
\end{eqnarray*}

Let $$\widehat{\Phi}(U)=\sum_{l\leq q}\sum_{q+1\leq m}\left<(\nabla^B_U{\bf J})_{x_0}w_l,\overline{w}_m\right>\overline{w}^m
 \wedge i_{\overline{w}_l},~~{\rm for}~~U\in TX.\eqno(4.156)$$
 $$\widetilde{\Phi}(U)=\sum_{m\leq q}\sum_{q+1\leq l}\left<(\nabla^B_U{\bf J})_{x_0}w_l,\overline{w}_m\right>\overline{w}^m
 \wedge i_{\overline{w}_l},~~{\rm for}~~U\in TX.\eqno(4.157)$$
By (4.89) and (4.156), then
$$({\mathcal{L}}^{-1}_0P^{N^\bot}{\mathcal{O}}''_1)P^N
=-\frac{\sqrt{-1}}{2}\widehat{\Phi}(\partial_{z_s})z_s-\frac{\sqrt{-1}}{2}\widehat{\Phi}(\partial_{\overline{z}_s})\overline{z}'_s
-\frac{\sqrt{-1}}{6\pi}\widehat{\Phi}(\partial_{\overline{z}_s})b_s.\eqno(4.158)$$
By (4.89) and (4.156)-(4.158), similar to the computations of $V^0_a$, we get
$$V^1_a=\left[-\frac{7}{96\times 120\pi^3}{\mathcal{R}}^\top\widehat{\Phi}(\partial_{z_s})\widehat{\Phi}(\partial_{\overline{z}_s})
+\frac{1}{288\pi^3}{\mathcal{R}}^\bot\widetilde{\Phi}(\partial_{z_s})\widehat{\Phi}(\partial_{\overline{z}_s})\right.$$
$$\left.
-\frac{7}{32\times 120\pi^3}{\mathcal{R}}^\top\widehat{\Phi}(\partial_{\overline{z}_s}){\mathcal{R}}^\top\widehat{\Phi}(\partial_{z_s})
+\frac{1}{96\pi^3}{\mathcal{R}}^\bot\widetilde{\Phi}(\partial_{\overline{z}_s})\widehat{\Phi}(\partial_{z_s})\right]
I_{{\rm det}(\overline{W}^*)\otimes E}({\mathcal{R}}^{\bot,*})(x)I_{4}.\eqno(4.159)
$$
Set
$$V^2_a=\sum_{0\leq m_1+m_2\leq {j-\frac{q}{2}}}I_{2j}({\mathcal{L}}^{-1}_0{\mathcal{O}}^{+2}_2)^{m_1}({\mathcal{L}}^{-1}_0
P^{N^\bot}{\mathcal{O}}'_1)
({\mathcal{L}}^{-1}_0{\mathcal{O}}^{+2}_2)^{m_2}$$
$$\cdot({\mathcal{L}}^{-1}_0P^{N^\bot}{\mathcal{O}}''_1)({\mathcal{L}}^{-1}_0{\mathcal{O}}^{+2}_2)^{{j-\frac{q}{2}}-m_1-m_2}
P^N({\mathcal{O}}^{-2}_2{\mathcal{L}}^{-1}_0)^{j-\frac{q}{2}}I_{2j}.
\eqno(4.160)$$
By $n=4$, $j=2$, $q=2$, then
\begin{eqnarray*}
V^2_a&=&I_{4}[({\mathcal{L}}^{-1}_0{\mathcal{O}}^{+2}_2)({\mathcal{L}}^{-1}_0
P^{N^\bot}{\mathcal{O}}'_1)
({\mathcal{L}}^{-1}_0P^{N^\bot}{\mathcal{O}}''_1)\\
&&
+({\mathcal{L}}^{-1}_0
P^{N^\bot}{\mathcal{O}}'_1)
({\mathcal{L}}^{-1}_0{\mathcal{O}}^{+2}_2)({\mathcal{L}}^{-1}_0P^{N^\bot}{\mathcal{O}}''_1)]
P^N({\mathcal{O}}^{-2}_2{\mathcal{L}}^{-1}_0)I_{4}.~~~~~~~~~~~~~~~~~(4.161)
\end{eqnarray*}
By (4.158), similar to  the computations of $V^0_a$, then
$$({\mathcal{L}}^{-1}_0
P^{N^\bot}A_1)
({\mathcal{L}}^{-1}_0{\mathcal{O}}^{+2}_2)({\mathcal{L}}^{-1}_0P^{N^\bot}{\mathcal{O}}''_1)P^N(0,Z)$$
$$
=\frac{\sqrt{-1}}{648\pi^3}
(\partial_{{z}_\alpha }R^L)_{x_0}(\partial_{\overline{z}_{\alpha}},\partial_{{z}_{s}}){\mathcal{R}}^{0}\widehat{\Phi}(\partial_{\overline{z}_s})
I_{{\rm det}(\overline{W}^*)\otimes E}P^N(0,Z);\eqno(4.162)$$

$$({\mathcal{L}}^{-1}_0
P^{N^\bot}A_2)
({\mathcal{L}}^{-1}_0{\mathcal{O}}^{+2}_2)({\mathcal{L}}^{-1}_0P^{N^\bot}{\mathcal{O}}''_1)P^N(0,Z)$$
$$=\frac{-\sqrt{-1}}{864\pi^3}[\frac{4}{3}(\partial_{{z}_{s} }R^L)_{x_0}(\partial_{z_\alpha},\partial_{\overline{z}_\alpha})+
\frac{4}{3}(\partial_{{z}_{\alpha} }R^L)_{x_0}(\partial_{z_{s}},\partial_{\overline{z}_\alpha})$$
$$-
\frac{1}{3}(\partial_{\alpha}R^L)_{x_0}(\partial_{z_{s}},e_{\alpha})]
{\mathcal{R}}^{0}\widehat{\Phi}(\partial_{\overline{z}_s})
I_{{\rm det}(\overline{W}^*)\otimes E}P^N(0,Z);\eqno(4.163)$$
\begin{eqnarray*}
&&({\mathcal{L}}^{-1}_0
P^{N^\bot}A_3)
({\mathcal{L}}^{-1}_0{\mathcal{O}}^{+2}_2)({\mathcal{L}}^{-1}_0P^{N^\bot}{\mathcal{O}}''_1)P^N(0,Z)\\
&&=\frac{\sqrt{-1}}{96\pi^2}[\frac{2}{3\pi}(\partial_{\overline{z}_r }R^L)_{x_0}(\partial_{z_r},\partial_{\overline{z}_\alpha})
+\frac{2}{3\pi}(\partial_{\overline{z}_\alpha }R^L)_{x_0}(\partial_{z_r},\partial_{\overline{z}_r})\\
&&-
\frac{1}{6\pi}(\partial_{j'}R^L)_{x_0}(\partial_{\overline{z}_\alpha},e_{j'})]
{\mathcal{R}}^{0}\widehat{\Phi}(\partial_{z_\alpha})I_{{\rm det}(\overline{W}^*)\otimes E}\\
&&+\frac{\sqrt{-1}}{1296\pi^3}[(\partial_{{z}_{s} }R^L)_{x_0}(\partial_{z_\alpha},\partial_{\overline{z}_\alpha})
{\mathcal{R}}^{0}\widehat{\Phi}(\partial_{\overline{z}_s})\\
&&+
(\partial_{{z}_{l'} }R^L)_{x_0}(\partial_{z_s},\partial_{\overline{z}_{l'}})
{\mathcal{R}}^{0}\widehat{\Phi}(\partial_{\overline{z}_s})]I_{{\rm det}(\overline{W}^*)\otimes E}P^N(0,Z);~~~~~~~~~~~~~~~~~~~~~~~~~~~(4.164)
\end{eqnarray*}

$$({\mathcal{L}}^{-1}_0
P^{N^\bot}A_4)
({\mathcal{L}}^{-1}_0{\mathcal{O}}^{+2}_2)({\mathcal{L}}^{-1}_0P^{N^\bot}{\mathcal{O}}''_1)P^N(0,Z)$$
$$=\frac{-\sqrt{-1}}{192\pi^3}[(\partial_{\overline{z}_s }R^L)_{x_0}(\partial_{z_\alpha},\partial_{\overline{z}_\alpha})
{\mathcal{R}}^{0}\widehat{\Phi}(\partial_{{z}_s})$$
$$+
(\partial_{\overline{z}_\alpha }R^L)_{x_0}(\partial_{z_\alpha},\partial_{\overline{z}_s})
{\mathcal{R}}^{0}\widehat{\Phi}(\partial_{{z}_s})]I_{{\rm det}(\overline{W}^*)\otimes E}P^N(0,Z);\eqno(4.165)$$

$$({\mathcal{L}}^{-1}_0
P^{N^\bot}A_5)
({\mathcal{L}}^{-1}_0{\mathcal{O}}^{+2}_2)({\mathcal{L}}^{-1}_0P^{N^\bot}{\mathcal{O}}''_1)P^N(0,Z)$$
$$=\frac{-\sqrt{-1}}{96\pi^3}[\frac{1}{3}
(\partial_{\overline{z}_r }R^L)_{x_0}(\partial_{z_\alpha},\partial_{\overline{z}_\alpha})
{\mathcal{R}}^{0}\widehat{\Phi}(\partial_{{z}_r})$$
$$-
\frac{1}{4}
(\partial_{j'}R^L)_{x_0}(\partial_{\overline{z}_r},e_{j'})
{\mathcal{R}}^{0}\widehat{\Phi}(\partial_{{z}_r})]P^N(0,Z).\eqno(4.166)$$

$$({\mathcal{L}}^{-1}_0{\mathcal{O}}^{+2}_2)({\mathcal{L}}^{-1}_0
P^{N^\bot}A_1)
({\mathcal{L}}^{-1}_0P^{N^\bot}{\mathcal{O}}''_1)P^N(0,Z)$$
$$
=\frac{5\sqrt{-1}}{1296\pi^3}(\partial_{{z}_\alpha }R^L)_{x_0}(\partial_{\overline{z}_\alpha},\partial_{{z}_s})
{\mathcal{R}}^{0}\widehat{\Phi}(\partial_{\overline{z}_s})I_{{\rm det}(\overline{W}^*)\otimes E}P^N(0,Z);\eqno(4.167)$$

$$({\mathcal{L}}^{-1}_0{\mathcal{O}}^{+2}_2)({\mathcal{L}}^{-1}_0
P^{N^\bot}A_2)
({\mathcal{L}}^{-1}_0P^{N^\bot}{\mathcal{O}}''_1)P^N(0,Z)$$
$$
=\frac{-\sqrt{-1}}{1728\pi^3}
[\frac{4}{3}(\partial_{{z}_{s} }R^L)_{x_0}(\partial_{z_\alpha},\partial_{\overline{z}_\alpha})+
\frac{4}{3}(\partial_{{z}_{\alpha} }R^L)_{x_0}(\partial_{z_{s}},\partial_{\overline{z}_\alpha})$$
$$-
\frac{1}{3}(\partial_{\alpha}R^L)_{x_0}(\partial_{z_{s}},e_{\alpha})]
{\mathcal{R}}^{0}\widehat{\Phi}(\partial_{\overline{z}_s})
I_{{\rm det}(\overline{W}^*)\otimes E}P^N(0,Z);\eqno(4.168)$$

\begin{eqnarray*}
&&({\mathcal{L}}^{-1}_0{\mathcal{O}}^{+2}_2)({\mathcal{L}}^{-1}_0
P^{N^\bot}A_3)
({\mathcal{L}}^{-1}_0P^{N^\bot}{\mathcal{O}}''_1)P^N(0,Z)\\
&&=\frac{\sqrt{-1}}{144\pi^2}[\frac{2}{3\pi}(\partial_{\overline{z}_r }R^L)_{x_0}(\partial_{z_r},\partial_{\overline{z}_\alpha})
+\frac{2}{3\pi}(\partial_{\overline{z}_\alpha }R^L)_{x_0}(\partial_{z_r},\partial_{\overline{z}_r})\\
&&-
\frac{1}{6\pi}(\partial_{j'}R^L)_{x_0}(\partial_{\overline{z}_\alpha},e_{j'})]
{\mathcal{R}}^{0}\widehat{\Phi}(\partial_{z_\alpha})I_{{\rm det}(\overline{W}^*)\otimes E}\\
&&+\frac{7\sqrt{-1}}{36\times 144\pi^3}[(\partial_{{z}_{s} }R^L)_{x_0}(\partial_{z_\alpha},\partial_{\overline{z}_\alpha})
{\mathcal{R}}^{0}\widehat{\Phi}(\partial_{\overline{z}_s})\\
&&+
(\partial_{{z}_{l'} }R^L)_{x_0}(\partial_{z_s},\partial_{\overline{z}_{l'}})
{\mathcal{R}}^{0}\widehat{\Phi}(\partial_{\overline{z}_s})]I_{{\rm det}(\overline{W}^*)\otimes E}P^N(0,Z);~~~~~~~~~~~~~~~~~~~~~~~~~~~(4.169)
\end{eqnarray*}

$$({\mathcal{L}}^{-1}_0{\mathcal{O}}^{+2}_2)({\mathcal{L}}^{-1}_0
P^{N^\bot}A_4)
({\mathcal{L}}^{-1}_0P^{N^\bot}{\mathcal{O}}''_1)P^N(0,Z)$$
$$=\frac{-\sqrt{-1}}{384\pi^3}[(\partial_{\overline{z}_s }R^L)_{x_0}(\partial_{z_\alpha},\partial_{\overline{z}_\alpha})
{\mathcal{R}}^{0}\widehat{\Phi}(\partial_{{z}_s})$$
$$+
(\partial_{\overline{z}_\alpha }R^L)_{x_0}(\partial_{z_\alpha},\partial_{\overline{z}_s})
{\mathcal{R}}^{0}\widehat{\Phi}(\partial_{{z}_s})]I_{{\rm det}(\overline{W}^*)\otimes E}P^N(0,Z);\eqno(4.170)$$

$$({\mathcal{L}}^{-1}_0{\mathcal{O}}^{+2}_2)({\mathcal{L}}^{-1}_0
P^{N^\bot}A_5)
({\mathcal{L}}^{-1}_0P^{N^\bot}{\mathcal{O}}''_1)P^N(0,Z)$$
$$=\frac{-\sqrt{-1}}{96\pi^3}[\frac{5}{9}
(\partial_{\overline{z}_r }R^L)_{x_0}(\partial_{z_\alpha},\partial_{\overline{z}_\alpha})
{\mathcal{R}}^{0}\widehat{\Phi}(\partial_{{z}_r})$$
$$-
\frac{1}{4}
(\partial_{j'}R^L)_{x_0}(\partial_{\overline{z}_r},e_{j'})
{\mathcal{R}}^{0}\widehat{\Phi}(\partial_{{z}_r})]P^N(0,Z).\eqno(4.171)$$
By (4.97) and (4.161), then
$$V^2_a=\frac{1}{4\pi}[(4.162)+\cdots+(4.171)]\times {\mathcal{R}}^{\bot,*}I_4.\eqno(4.172)$$

Set
$$V^3_a=\sum_{0\leq m_1+m_2\leq {j-\frac{q}{2}}}I_{2j}({\mathcal{L}}^{-1}_0{\mathcal{O}}^{+2}_2)^{m_1}({\mathcal{L}}^{-1}_0
P^{N^\bot}{\mathcal{O}}''_1)
({\mathcal{L}}^{-1}_0{\mathcal{O}}^{+2}_2)^{m_2}$$
$$\cdot({\mathcal{L}}^{-1}_0P^{N^\bot}{\mathcal{O}}'_1)({\mathcal{L}}^{-1}_0{\mathcal{O}}^{+2}_2)^{{j-\frac{q}{2}}-m_1-m_2}
P^N({\mathcal{O}}^{-2}_2{\mathcal{L}}^{-1}_0)^{j-\frac{q}{2}}I_{2j}.
\eqno(4.173)$$
By $n=4$, $j=2$, $q=2$, then
$$
V^3_a=I_{4}({\mathcal{L}}^{-1}_0{\mathcal{O}}^{+2}_2)({\mathcal{L}}^{-1}_0
P^{N^\bot}{\mathcal{O}}''_1)
({\mathcal{L}}^{-1}_0P^{N^\bot}{\mathcal{O}}'_1)
P^N({\mathcal{O}}^{-2}_2{\mathcal{L}}^{-1}_0)I_{4}.\eqno(4.174)
$$
By (4.89) and
$({\mathcal{L}}^{-1}_0P^{N^\bot}{\mathcal{O}}'_1)
P^N$ being a $(0,(2,0))$-form and ${\mathcal{L}}_0A_2P^N=0$, then
$${\mathcal{O}}''_1
({\mathcal{L}}^{-1}_0P^{N^\bot}{\mathcal{O}}'_1)
P^N
=-4\pi\sqrt{-1}[\widehat{\Phi}(\partial_{{z}_s})z_s+\widehat{\Phi}(\partial_{\overline{z}_s})\overline{z}_s]
(\frac{1}{4\pi}A_3+\frac{1}{8\pi}A_4)P^N.\eqno(4.175)$$
Direct computations show that
$$z_sb_\alpha=b_\alpha z_s+2\delta_{\alpha s},~~z_sb_lz_\alpha\overline{z}'_r=2\delta_{sl}z_\alpha\overline{z}'_r
+b_lz_sz_\alpha\overline{z}'_r;\eqno(4.176)$$
$$z_sb_lz_\alpha z_{j'}=2\delta_{sl}z_\alpha z_{j'}+b_lz_sz_\alpha z_{j'}, ~~z_sb_lb_rz_\alpha=2\delta_{sl}b_rz_\alpha
+2\delta_{sr}b_lz_\alpha+b_lb_rz_sz_\alpha;\eqno(4.177)$$
$$\overline{z}_sb_\alpha{\mathcal{P}}=(\frac{b_\alpha b_s}{2\pi}+b_\alpha \overline{z}'_s){\mathcal{P}};\eqno(4.178)$$
$$\overline{z}_sb_lz_\alpha\overline{z}'_r{\mathcal{P}}=[\frac{1}{2\pi}b_lb_sz_\alpha\overline{z}'_r+
\frac{\delta_{\alpha s}}{\pi}b_l\overline{z}'_r+b_lz_\alpha\overline{z}'_r\overline{z}'_s]{\mathcal{P}};\eqno(4.179)$$
$$\overline{z}_sb_lz_\alpha z_{j'}{\mathcal{P}}=[\frac{b_lb_sz_\alpha z_{j'}}{2\pi}+\frac{\delta_{\alpha s}}{\pi}b_lz_{j'}
+\frac{\delta_{j' s}}{\pi}b_lz_{\alpha}+
b_lz_\alpha z_{j'}\overline{z}'_s]{\mathcal{P}};\eqno(4.180)$$
$$\overline{z}_sb_lb_rz_\alpha{\mathcal{P}}=[\frac{b_lb_rb_sz_\alpha}{2\pi}+\frac{\delta_{\alpha s}}{\pi}b_lz_{j'}+
\frac{\delta_{j' s}}{\pi}b_lz_{\alpha}+b_lz_\alpha z_{j'}\overline{z}'_s]{\mathcal{P}};\eqno(4.181)$$
By (4.174)-(4.181), we can get
$$
V^3_a=\frac{-\sqrt{-1}}{4\pi}{\mathcal{R}}^0\left\{\frac{11}{8\times 144\pi^3}\widehat{\Phi}(\partial_{z_s})
[(\partial_{\overline{z}_\alpha }R^L)_{x_0}(\partial_{{z}_\alpha},\partial_{\overline{z}_s})+
(\partial_{\overline{z}_s }R^L)_{x_0}(\partial_{{z}_\alpha},\partial_{\overline{z}_\alpha})]\right.$$
$$
-\frac{5}{6\times 288\pi^3}\widehat{\Phi}(\partial_{z_s})
(\partial_{j'}R^L)_{x_0}(\partial_{\overline{z}_s},e_{j'})-
\frac{7}{12\times 144\pi^3}\widehat{\Phi}(\partial_{\overline{z}_s})[
(\partial_{{z}_s }R^L)_{x_0}(\partial_{z_\alpha},\partial_{\overline{z}_\alpha})$$
$$\left.+(\partial_{{z}_l }R^L)_{x_0}(\partial_{z_s},\partial_{\overline{z}_l})]\right\}
{\mathcal{R}}^{\bot,*}I_4.\eqno(4.182)$$
By (4.84),(4.90),(4.154),(4.160) and (4.173), we get
$$V_a=V^0_a+V^1_a+V^2_a+V^3_a.\eqno(4.183)$$
\indent Nextly, we compute $V_b$. Set\\
$${\rm V}^0_b=
\sum_{0\leq m_1,m_2\leq {j-\frac{q}{2}}}I_{2j}({\mathcal{L}}^{-1}_0{\mathcal{O}}^{+2}_2)^{m_1}
({\mathcal{L}}^{-1}_0P^{N^\bot}{\mathcal{O}}'_1)({\mathcal{L}}^{-1}_0{\mathcal{O}}^{+2}_2)^{{j-\frac{q}{2}}-m_1}
P^N$$
$$\times({\mathcal{O}}^{-2}_2{\mathcal{L}}^{-1}_0)^{m_2}
({\mathcal{O}}'_1{\mathcal{L}}^{-1}_0P^{N^\bot})({\mathcal{O}}^{-2}_2{\mathcal{L}}^{-1}_0)^{{j-\frac{q}{2}}-m_2}I_{2j}.\eqno(4.184)$$
\indent By (4.184), then
$${\rm V}^0_b=I_{2j}CC^*I_{2j},~~C=\sum_{0\leq m_1\leq {j-\frac{q}{2}}}C_{m_1},\eqno(4.185)$$
where $$C_{m_1}=({\mathcal{L}}^{-1}_0{\mathcal{O}}^{+2}_2)^{m_1}
({\mathcal{L}}^{-1}_0P^{N^\bot}{\mathcal{O}}'_1)({\mathcal{L}}^{-1}_0{\mathcal{O}}^{+2}_2)^{{j-\frac{q}{2}}-m_1}
P^N.\eqno(4.186)$$
\noindent Similar to (4.96), then
\begin{eqnarray*}
&&C_{m_1}(0,Z)\\
&=&\frac{1}{(4\pi)^{{j-\frac{q}{2}}+1}}
( {\mathcal{R}}^\bot)^{j-\frac{q}{2}}\left\{(\frac{4}{3}B^{1,{j-\frac{q}{2}}-m_1}_0B^{{j-\frac{q}{2}}-m_1,{j-\frac{q}{2}}}_0\right.~~~~
~~~~~~~~~~~~~~~~~~~~~~~~~~(4.187)\\
&&
-\frac{8}{3}B^{1,{j-\frac{q}{2}}-m_1}_0B^{{j-\frac{q}{2}}-m_1,{j-\frac{q}{2}}}_1
+
\frac{4}{3}B^{1,{j-\frac{q}{2}}-m_1}_0B^{{j-\frac{q}{2}}-m_1,{j-\frac{q}{2}}}_2)(\partial_{\overline{z}_{r} }R^L)_{x_0}(\partial_{{z}_l},\partial_{\overline{z}_{l}})\\
&&+(
\frac{4}{3}B^{1,{j-\frac{q}{2}}-m_1}_0B^{{j-\frac{q}{2}}-m_1,{j-\frac{q}{2}}}_2-
\frac{4}{3}B^{1,{j-\frac{q}{2}}-m_1}_0B^{{j-\frac{q}{2}}-m_1,{j-\frac{q}{2}}}_1)(\partial_{\overline{z}_{\alpha} }R^L)_{x_0}(\partial_{{z}_\alpha},\partial_{\overline{z}_{r}})\\
&&+(\frac{1}{3}B^{1,{j-\frac{q}{2}}-m_1}_0B^{{j-\frac{q}{2}}-m_1,{j-\frac{q}{2}}}_1-\frac{1}{3}B^{1,{j-\frac{q}{2}}-m_1}_0
B^{{j-\frac{q}{2}}-m_1,{j-\frac{q}{2}}}_0)\\
&&\left.\cdot(\partial_{l }R^L)_{x_0}(\partial_{\overline{z}_r},e_l)\right\}\overline{z}_r[{\mathcal{P}}(0,Z)].
\end{eqnarray*}
So
\begin{eqnarray*}
&&(C^*_{m_1})(Z,0)\\
&=&-\frac{1}{(4\pi)^{{j-\frac{q}{2}}+1}}
\left\{(\frac{4}{3}B^{1,{j-\frac{q}{2}}-m_1}_0B^{{j-\frac{q}{2}}-m_1,{j-\frac{q}{2}}}_0-\frac{8}{3}B^{1,{j-\frac{q}{2}}-m_1}_0
B^{{j-\frac{q}{2}}-m_1,{j-\frac{q}{2}}}_1
\right.\\
&&+
\frac{4}{3}B^{1,{j-\frac{q}{2}}-m_1}_0B^{{j-\frac{q}{2}}-m_1,{j-\frac{q}{2}}}_2)(\partial_{{z}_{r} }R^L)_{x_0}(\partial_{\overline{z}_l},\partial_{{z}_{l}})\\
&&+(
\frac{4}{3}B^{1,{j-\frac{q}{2}}-m_1}_0B^{{j-\frac{q}{2}}-m_1,{j-\frac{q}{2}}}_2-
\frac{4}{3}B^{1,{j-\frac{q}{2}}-m_1}_0B^{{j-\frac{q}{2}}-m_1,{j-\frac{q}{2}}}_1)(\partial_{{z}_{\alpha} }R^L)_{x_0}(\partial_{\overline{z}_\alpha},\partial_{z_r})\\
&&+(\frac{1}{3}B^{1,{j-\frac{q}{2}}-m_1}_0B^{{j-\frac{q}{2}}-m_1,{j-\frac{q}{2}}}_1-\frac{1}{3}B^{1,{j-\frac{q}{2}}-m_1}_0
B^{{j-\frac{q}{2}}-m_1,{j-\frac{q}{2}}}_0)\\
&&\left.\cdot
(\partial_{l }R^L)_{x_0}(\partial_{z_r},e_l)\right\}( {\mathcal{R}}^{\bot,*})^{j-\frac{q}{2}}z_r[{\mathcal{P}}(Z,0)].~~~~~~~~~~~~~~~~~~~~~~~~~~~~~~~~~~~~~~~~~(4.188)
\end{eqnarray*}
By (4.24), (4.185)-(4.188), we get
\begin{eqnarray*}
{\rm V}^0_b
&=&\frac{-2}{9\times(4\pi)^{2j-q+3}}I_{2j}\sum_{0\leq m_1,m_2\leq {j-\frac{q}{2}}}(B^{1,{j-\frac{q}{2}}-m_1}_0)^2
\left[(B^{{j-\frac{q}{2}}-m_1,{j-\frac{q}{2}}}_0\right.~~~~~~~~~~~(4.189)\\
&&-2B^{{j-\frac{q}{2}}-m_1,{j-\frac{q}{2}}}_1+B^{{j-\frac{q}{2}}-m_1,{j-\frac{q}{2}}}_2)
A_{\overline{r}s\overline{s}}
\\
&&\left.
+(B^{{{j-\frac{q}{2}}}
-m_1,{{j-\frac{q}{2}}}}_2-B^{{{j-\frac{q}{2}}}-m_1,{{j-\frac{q}{2}}}}_1)
A_{\overline{s}s\overline{r}}\right.\\
&&\left.+\frac{1}{2}(B^{{{j-\frac{q}{2}}}-m_1,{{j-\frac{q}{2}}}}_1-B^{{{j-\frac{q}{2}}}-m_1,
{{j-\frac{q}{2}}}}_0)
A_{e_l\overline{r}e_l}\right]\\
&&\times \left[(B^{{{j-\frac{q}{2}}}-m_2,{{j-\frac{q}{2}}}}_0-2B^{{{j-\frac{q}{2}}}-m_2,{j-\frac{q}{2}}}_1
+B^{{j-\frac{q}{2}}-m_2,{j-\frac{q}{2}}}_2)A_{{r}\overline{s'}s'}\right.\\
&&
+(B^{{j-\frac{q}{2}}-m_2,{j-\frac{q}{2}}}_2-B^{{j-\frac{q}{2}}-m_2,{j-\frac{q}{2}}}_1)A_{s'\overline{s'}{r}}\\
&&\left.
+\frac{1}{2}(B^{{j-\frac{q}{2}}-m_2,{j-\frac{q}{2}}}_1-B^{{j-\frac{q}{2}}-m_2,{j-\frac{q}{2}}}_0)A_{e_{l'}{r}e_{l'}}\right]
({\mathcal{R}}^\bot)(x)^{{j-\frac{q}{2}}}( {\mathcal{R}}^{\bot,*})(x)^{{j-\frac{q}{2}}}I_{2j}.
\end{eqnarray*}
When $n=4$, $j=2$ and $q=2$, by (4.158), then
$$
\sum_{0\leq m_1\leq {j-\frac{q}{2}}}I_{2j}({\mathcal{L}}^{-1}_0{\mathcal{O}}^{+2}_2)^{m_1}
({\mathcal{L}}^{-1}_0P^{N^\bot}{\mathcal{O}}''_1)({\mathcal{L}}^{-1}_0{\mathcal{O}}^{+2}_2)^{{j-\frac{q}{2}}-m_1}
P^N(0,Z)$$
$$=I_{4}({\mathcal{L}}^{-1}_0{\mathcal{O}}^{+2}_2)
({\mathcal{L}}^{-1}_0P^{N^\bot}{\mathcal{O}}''_1)
P^N(0,Z)=\frac{-5\sqrt{-1}}{144\pi}{\mathcal{R}}^0\widehat{\Phi}(\partial_{\overline{z}_s})\overline{z}_sP^N(0,Z).\eqno(4.190)$$
Set
$${\rm V}^1_b=
\sum_{0\leq m_1,m_2\leq {j-\frac{q}{2}}}I_{2j}({\mathcal{L}}^{-1}_0{\mathcal{O}}^{+2}_2)^{m_1}
({\mathcal{L}}^{-1}_0P^{N^\bot}{\mathcal{O}}''_1)({\mathcal{L}}^{-1}_0{\mathcal{O}}^{+2}_2)^{{j-\frac{q}{2}}-m_1}
P^N$$
$$\times({\mathcal{O}}^{-2}_2{\mathcal{L}}^{-1}_0)^{m_2}
({\mathcal{O}}'_1{\mathcal{L}}^{-1}_0P^{N^\bot})({\mathcal{O}}^{-2}_2{\mathcal{L}}^{-1}_0)^{{j-\frac{q}{2}}-m_2}I_{2j}.\eqno(4.191)$$
$${\rm V}^2_b=
\sum_{0\leq m_1,m_2\leq {j-\frac{q}{2}}}I_{2j}({\mathcal{L}}^{-1}_0{\mathcal{O}}^{+2}_2)^{m_1}
({\mathcal{L}}^{-1}_0P^{N^\bot}{\mathcal{O}}'_1)({\mathcal{L}}^{-1}_0{\mathcal{O}}^{+2}_2)^{{j-\frac{q}{2}}-m_1}
P^N$$
$$\times({\mathcal{O}}^{-2}_2{\mathcal{L}}^{-1}_0)^{m_2}
({\mathcal{O}}''_1{\mathcal{L}}^{-1}_0P^{N^\bot})({\mathcal{O}}^{-2}_2{\mathcal{L}}^{-1}_0)^{{j-\frac{q}{2}}-m_2}I_{2j}.\eqno(4.192)$$
$${\rm V}^3_b=
\sum_{0\leq m_1,m_2\leq {j-\frac{q}{2}}}I_{2j}({\mathcal{L}}^{-1}_0{\mathcal{O}}^{+2}_2)^{m_1}
({\mathcal{L}}^{-1}_0P^{N^\bot}{\mathcal{O}}''_1)({\mathcal{L}}^{-1}_0{\mathcal{O}}^{+2}_2)^{{j-\frac{q}{2}}-m_1}
P^N$$
$$\times({\mathcal{O}}^{-2}_2{\mathcal{L}}^{-1}_0)^{m_2}
({\mathcal{O}}''_1{\mathcal{L}}^{-1}_0P^{N^\bot})({\mathcal{O}}^{-2}_2{\mathcal{L}}^{-1}_0)^{{j-\frac{q}{2}}-m_2}I_{2j}.\eqno(4.193)$$
By (4.187),(4.188) and (4.190), then
\begin{eqnarray*}
{\rm V}^1_b&=&\frac{5\sqrt{-1}}{144\times 16\pi^3}{\mathcal{R}}^0\widehat{\Phi}(\partial_{\overline{z}_r})
\left\{(\frac{4}{3}B^{1,1}_0B^{1,1}_0-\frac{8}{3}B^{1,1}_0
B^{1,1}_1
\right.\\
&&+
\frac{4}{3}B^{1,1}_0B^{1,1}_2)(\partial_{{z}_{r} }R^L)_{x_0}(\partial_{\overline{z}_l},\partial_{{z}_{l}})\\
&&+(
\frac{4}{3}B^{1,1}_0B^{1,1}_2-
\frac{4}{3}B^{1,1}_0B^{1,1}_1)(\partial_{{z}_{\alpha} }R^L)_{x_0}(\partial_{\overline{z}_\alpha},\partial_{z_r})\\
&&+(\frac{1}{3}B^{1,1}_0B^{1,1}_1-\frac{1}{3}B^{1,1}_0
B^{1,1}_0)\\
&&\left.\cdot
(\partial_{l }R^L)_{x_0}(\partial_{z_r},e_l)\right\}
I_{{\rm det}(\overline{W}^*)\otimes E}{\mathcal{R}}^{\bot,*}I_4~~~~~~~~~~~~~~~~~~~~~~~~~~~~~~~~~~~~~~~~~~~(4.194)
\end{eqnarray*}
\begin{eqnarray*}
{\rm V}^2_b&=&\frac{5\sqrt{-1}}{144\times 16\pi^3}
\left\{(\frac{4}{3}B^{1,1}_0B^{1,1}_0\right.\\
&&
-\frac{8}{3}B^{1,1}_0B^{1,1}_1
+
\frac{4}{3}B^{1,1}_0B^{1,1}_2)(\partial_{\overline{z}_{r} }R^L)_{x_0}(\partial_{{z}_l},\partial_{\overline{z}_{l}})\\
&&+(
\frac{4}{3}B^{1,1}_0B^{1,1}_2-
\frac{4}{3}B^{1,1}_0B^{1,1}_1)(\partial_{\overline{z}_{\alpha} }R^L)_{x_0}(\partial_{{z}_\alpha},\partial_{\overline{z}_{r}})\\
&&+(\frac{1}{3}B^{1,1}_0B^{1,1}_1-\frac{1}{3}B^{1,1}_0
B^{1,1}_0)\\
&&\left.\cdot(\partial_{l }R^L)_{x_0}(\partial_{\overline{z}_r},e_l)\right\}{\mathcal{R}}^\bot
I_{{\rm det}(\overline{W}^*)\otimes E}\widehat{\Phi}(\partial_{{z}_r})^*{\mathcal{R}}^{0.*}I_4
~~~~~~~~~~~~~~~~~~~~~~~~~~~~~~(4.195)
\end{eqnarray*}
$$V^3_b=\frac{25}{(144\pi)^2}{\mathcal{R}}^0\widehat{\Phi}(\partial_{\overline{z}_r})
I_{{\rm det}(\overline{W}^*)\otimes E}\widehat{\Phi}^*(\partial_{{z}_r}){\mathcal{R}}^{0.*}I_4.\eqno(4.196)$$
By (4.85),(4.184),(4.191)-(4.193), we get
$$V_b=V^0_b+V^1_b+V^2_b+V^3_b.\eqno(4.197)$$\\

\noindent{\bf 4.4 The computations of the term IV}\\

\indent In this section, we will compute the term IV. By the discussions after (4.2), for $l=1$, then $4j\leq 2q+2k-2\leq 4j+1$, so $k=2j-q+1$. There are $2j - 1-q$ ${\mathcal{O}}_{r_i}$ equal to ${\mathcal{O}}_{2}$ and $1$ equal to ${\mathcal{O}}_{3}$ and $1$ equal to ${\mathcal{O}}_{1}$. When
$l_1=1$, by (3.16) and (3.17), then $i_0=j+1-\frac{q}{2}$ is unique. When $l_1=0$, then $i_0=j+2-\frac{q}{2}$. So
$${\rm }IV={\rm IV}_a+{\rm IV}_b+{\rm IV}_c+{\rm IV}^*_a+{\rm IV}^*_b+{\rm IV}^*_c,\eqno(4.198)$$
where
$${\rm IV}_a=\sum_{0\leq m_1+m_2\leq {j-\frac{q}{2}}-1}I_{2j}({\mathcal{L}}_0^{-1}{\mathcal{O}}^{+2}_{2})^{{j-\frac{q}{2}}}P^N ({\mathcal{O}}^{-2}_{2}{\mathcal{L}}_0^{-1})^{m_1}({\mathcal{O}}_{1}{\mathcal{L}}_0^{-1})$$
$$\cdot({\mathcal{O}}^{-2}_{2}{\mathcal{L}}_0^{-1})^{m_2}
({\mathcal{O}}^{-2}_{3}{\mathcal{L}}_0^{-1})({\mathcal{O}}^{-2}_{2}{\mathcal{L}}_0^{-1})^{{j-\frac{q}{2}}-1-m_1-m_2}
I_{2j};\eqno(4.199)$$
$${\rm IV}_b=\sum_{0\leq m_1+m_2\leq {j-\frac{q}{2}}-1}I_{2j}({\mathcal{L}}_0^{-1}{\mathcal{O}}^{+2}_{2})^{{j-\frac{q}{2}}}P^N ({\mathcal{O}}^{-2}_{2}{\mathcal{L}}_0^{-1})^{m_1}({\mathcal{O}}^{-2}_{3}{\mathcal{L}}_0^{-1})$$
$$\cdot({\mathcal{O}}^{-2}_{2}{\mathcal{L}}_0^{-1})^{m_2}
({\mathcal{O}}_{1}{\mathcal{L}}_0^{-1})({\mathcal{O}}^{-2}_{2}{\mathcal{L}}_0^{-1})^{{j-\frac{q}{2}}-1-m_1-m_2}
I_{2j};\eqno(4.200)$$
$${\rm IV}_c=\sum_{0\leq m_1\leq {j-\frac{q}{2}}-1}\sum_{0\leq m_2\leq {j-\frac{q}{2}}}I_{2j}
({\mathcal{L}}_0^{-1}{\mathcal{O}}^{+2}_{2})^{m_1}
({\mathcal{L}}_0^{-1}{\mathcal{O}}^{+2}_{3})({\mathcal{L}}_0^{-1}{\mathcal{O}}^{+2}_{2})^{{j-\frac{q}{2}}-1-m_1}P^N$$
$$\cdot({\mathcal{O}}^{-2}_{2}{\mathcal{L}}_0^{-1})^{m_2}
({\mathcal{O}}_{1}{\mathcal{L}}_0^{-1})({\mathcal{O}}^{-2}_{2}{\mathcal{L}}_0^{-1})^{{j-\frac{q}{2}}-m_2}
I_{2j}.\eqno(4.201)$$
By (3.15) in [PZ], then
$${\mathcal{O}}^{+2}_{3}=z_{j'}\frac{\partial{\mathcal{R}}.}{\partial z_{j'}}(0)+
\overline{z}_{j'}\frac{\partial{\mathcal{R}}.}{\partial \overline{z}_{j'}}(0).\eqno(4.202)$$
By (4.199), we have
$${\rm IV}_a=\sum_{0\leq m_1+m_2\leq {j-\frac{q}{2}}-1}I_{2j}({\mathcal{L}}_0^{-1}{\mathcal{O}}^{+2}_{2})^{{j-\frac{q}{2}}}P^N
[({\mathcal{L}}_0^{-1}{\mathcal{O}}^{+2}_{2})^{{j-\frac{q}{2}}-1-m_1-m_2}$$
$$\cdot({\mathcal{L}}_0^{-1}{\mathcal{O}}^{+2}_{3})({\mathcal{L}}_0^{-1}{\mathcal{O}}^{+2}_{2})^{m_2}
({\mathcal{L}}_0^{-1}{\mathcal{O}}_{1})({\mathcal{L}}_0^{-1}{\mathcal{O}}^{+2}_{2})^{m_1}P^N]^*I_{2j}.\eqno(4.203)$$
Set
$${\rm IV}^1_a=\sum_{0\leq m_1+m_2\leq {j-\frac{q}{2}}-1}I_{2j}({\mathcal{L}}_0^{-1}{\mathcal{O}}^{+2}_{2})^{{j-\frac{q}{2}}}P^N
[({\mathcal{L}}_0^{-1}{\mathcal{O}}^{+2}_{2})^{{j-\frac{q}{2}}-1-m_1-m_2}$$
$$\cdot({\mathcal{L}}_0^{-1}{\mathcal{O}}^{+2}_{3})({\mathcal{L}}_0^{-1}{\mathcal{O}}^{+2}_{2})^{m_2}
({\mathcal{L}}_0^{-1}{\mathcal{O}}'_{1})({\mathcal{L}}_0^{-1}{\mathcal{O}}^{+2}_{2})^{m_1}P^N]^*I_{2j}.\eqno(4.204)$$
$${\rm IV}^2_a=\sum_{0\leq m_1+m_2\leq {j-\frac{q}{2}}-1}I_{2j}({\mathcal{L}}_0^{-1}{\mathcal{O}}^{+2}_{2})^{{j-\frac{q}{2}}}P^N
[({\mathcal{L}}_0^{-1}{\mathcal{O}}^{+2}_{2})^{{j-\frac{q}{2}}-1-m_1-m_2}$$
$$\cdot({\mathcal{L}}_0^{-1}{\mathcal{O}}^{+2}_{3})({\mathcal{L}}_0^{-1}{\mathcal{O}}^{+2}_{2})^{m_2}
({\mathcal{L}}_0^{-1}{\mathcal{O}}''_{1})({\mathcal{L}}_0^{-1}{\mathcal{O}}^{+2}_{2})^{m_1}P^N]^*I_{2j}.\eqno(4.205)$$
By (4.99)-(4.101) and direct computations, then
$$\left[({\mathcal{L}}_0^{-1}{\mathcal{O}}^{+2}_{2})^{j-\frac{q}{2}-1-m_1-m_2}{\mathcal{L}}_0^{-1}(z_{j'}\frac{\partial{\mathcal{R}}.}{\partial z_{j'}}(0))
A_2{\mathcal{P}}\right](0,Z)=0,\eqno(4.206)$$

\begin{eqnarray*}
&&\left[({\mathcal{L}}_0^{-1}{\mathcal{O}}^{+2}_{2})^{{j-\frac{q}{2}}-1-m_1-m_2}{\mathcal{L}}_0^{-1}(z_{j'}\frac{\partial{\mathcal{R}}.}{\partial z_{j'}}(0))
{\mathcal{R}}^{m_1+m_2}\right.\\
&&\left.\cdot\frac{1}{(4\pi)^{m_1+m_2+1}}B^{1,m_1}_0B^{m_1,m_1+m_2}_1A_3{\mathcal{P}}\right](0,Z)\\
&&=\frac{2}{(4\pi)^{{j-\frac{q}{2}}+1}}B^{1,m_1}_0B^{m_1,m_1+m_2}_1(B^{m_1+m_2+1,{j-\frac{q}{2}}}_0-B^{m_1+m_2+1,{j-\frac{q}{2}}}_1)\\
&&
\cdot({\mathcal{R}}^\bot)^{{j-\frac{q}{2}}-1-m_1-m_2}\frac{\partial{\mathcal{R}}^\bot.}{\partial z_{j'}}(0)({\mathcal{R}}^\bot)^{m_1+m_2}
[\frac{2}{3\pi}(\partial_{\overline{z}_{r} }R^L)_{x_0}(\partial_{{z}_r},\partial_{\overline{z}_{j'}})\\
&&+\frac{2}{3\pi}(\partial_{\overline{z}_{j'} }R^L)_{x_0}(\partial_{{z}_r},\partial_{\overline{z}_{r}})-
\frac{1}{6\pi}(\partial_{j'' }R^L)_{x_0}(\partial_{\overline{z}_{j'}},e_{j''})]{\mathcal{P}}(0,Z),~~~~~~~~~~~~~~~~~~~(4.207)
\end{eqnarray*}

\begin{eqnarray*}
&&\left[({\mathcal{L}}_0^{-1}{\mathcal{O}}^{+2}_{2})^{{j-\frac{q}{2}}-1-m_1-m_2}{\mathcal{L}}_0^{-1}(z_{j'}\frac{\partial{\mathcal{R}}^\bot.}{\partial z_{j'}}(0))\right.\\
&&\left.\cdot({\mathcal{R}}^\bot)^{m_1+m_2}\cdot\frac{1}{(4\pi)^{m_1+m_2+1}}B^{1,m_1}_0B^{m_1,m_1+m_2}_2A_4{\mathcal{P}}\right](0,Z)\\
&&=\frac{4}{3\pi\times(4\pi)^{{j-\frac{q}{2}}+1}}B^{1,m_1}_0B^{m_1,m_1+m_2}_2(B^{m_1+m_2+1,{j-\frac{q}{2}}}_2-B^{m_1+m_2+1,{j-\frac{q}{2}}}_1)\\
&&
\cdot{\mathcal{R}}^{{j-\frac{q}{2}}-1-m_1-m_2}\frac{\partial{\mathcal{R}}^\bot.}{\partial z_{j'}}(0)({\mathcal{R}}^\bot)^{m_1+m_2}\\
&&
\cdot[(\partial_{\overline{z}_{r} }R^L)_{x_0}(\partial_{{z}_r},\partial_{\overline{z}_{j'}})+
(\partial_{\overline{z}_{j'} }R^L)_{x_0}(\partial_{{z}_r},\partial_{\overline{z}_{r}})]{\mathcal{P}}(0,Z),~~~~~~~~~~~~~~~~~~~~~~~~~~~(4.208)
\end{eqnarray*}

\begin{eqnarray*}
&&\left[({\mathcal{L}}_0^{-1}{\mathcal{O}}^{+2}_{2})^{{j-\frac{q}{2}}-1-m_1-m_2}{\mathcal{L}}_0^{-1}(\overline{z}_{j'}
\frac{\partial{\mathcal{R}}^\bot.}{\partial \overline{z}_{j'}}(0))
({\mathcal{R}}^\bot)^{m_1+m_2}\right.\\
&&\left.\cdot\frac{1}{(4\pi)^{m_1+m_2+1}}B^{1,m_1}_0B^{m_1,m_1+m_2}_0A_2{\mathcal{P}}\right](0,Z)\\
&&=\frac{1}{3\pi\times(4\pi)^{{j-\frac{q}{2}}+1}}B^{1,m_1}_0B^{m_1,m_1+m_2}_0(B^{m_1+m_2+1,{j-\frac{q}{2}}}_0-B^{m_1+m_2+1,{j-\frac{q}{2}}}_1)\\
&&
\cdot({\mathcal{R}}^\bot)^{{j-\frac{q}{2}}-1-m_1-m_2}\frac{\partial{\mathcal{R}}^\bot.}{\partial \overline{z}_{j'}}(0)
({\mathcal{R}}^\bot)^{m_1+m_2}
[4(\partial_{z_{j'}}R^L)_{x_0}(\partial_{z_r},\partial_{\overline{z}_{r}})\\
&&+4(\partial_{z_r}R^L)_{x_0}(\partial_{{z}_{j'}},\partial_{\overline{z}_{r}})-
(\partial_{j'' }R^L)_{x_0}(\partial_{{z}_{j'}},e_{j''})]{\mathcal{P}}(0,Z),~~~~~~~~~~~~~~~~~~~~~~~~~~~(4.209)
\end{eqnarray*}
\begin{eqnarray*}
&&\left[({\mathcal{L}}_0^{-1}{\mathcal{O}}^{+2}_{2})^{{j-\frac{q}{2}}-1-m_1-m_2}{\mathcal{L}}_0^{-1}(\overline{z}_{j'}
\frac{\partial{\mathcal{R}}^\bot.}{\partial \overline{z}_{j'}}(0))
({\mathcal{R}}^\bot)^{m_1+m_2}\right.\\
&&\left.\cdot\frac{1}{(4\pi)^{m_1+m_2+1}}B^{1,m_1}_0B^{m_1,m_1+m_2}_1A_3{\mathcal{P}}\right](0,Z)\\
&&=\frac{4}{3\pi\times(4\pi)^{{j-\frac{q}{2}}+1}}B^{1,m_1}_0B^{m_1,m_1+m_2}_1(B^{m_1+m_2+1,{j-\frac{q}{2}}}_2-B^{m_1+m_2+1,{j-\frac{q}{2}}}_1)\\
&&\cdot({\mathcal{R}}^\bot)^{{j-\frac{q}{2}}-1-m_1-m_2}
\frac{\partial{\mathcal{R}}^\bot.}{\partial \overline{z}_{j'}}(0)({\mathcal{R}}^\bot)^{m_1+m_2}\\
&&\cdot
[(\partial_{z_{j'}}R^L)_{x_0}(\partial_{z_r},\partial_{\overline{z}_{r}})+
(\partial_{z_r}R^L)_{x_0}(\partial_{{z}_{j'}},\partial_{\overline{z}_{r}})],~~~~~~~~~~~~~~~~~~~~~~~~~~~~~~~~~~~~~~~(4.210)
\end{eqnarray*}
\begin{eqnarray*}
&&\left[({\mathcal{L}}_0^{-1}{\mathcal{O}}^{+2}_{2})^{{j-\frac{q}{2}}-1-m_1-m_2}{\mathcal{L}}_0^{-1}(\overline{z}_{j'}
\frac{\partial{\mathcal{R}}^\bot.}{\partial \overline{z}_{j'}}(0))
({\mathcal{R}}^\bot)^{m_1+m_2}\right.\\
&&\left.\cdot\frac{1}{(4\pi)^{m_1+m_2+1}}B^{1,m_1}_0B^{m_1,m_1+m_2}_2A_4{\mathcal{P}}\right](0,Z)=0.~~~~~~~~~~~~~~~~~~~~~~~~~~~~~~~~~~(4.211)
\end{eqnarray*}
By (4.96),(4.202),(4.204) and (4.206)-(4.211), we get
\begin{eqnarray*}
{\rm IV}^1_a&=&\frac{-1}{2^{j-\frac{q}{2}}({j-\frac{q}{2}})!(4\pi)^{2j-q+1}\times 3\pi}\sum_{0\leq m_1+m_2\leq j-\frac{q}{2}-1}B^{1,m_1}_0
I_{2j}({\mathcal{R}}^\bot)(x)^{{j-\frac{q}{2}}}\\
&&\times\left[(B^{m_1,m_1+m_2}_1B^{m_1+m_2+1,{j-\frac{q}{2}}}_0-B^{m_1,{j-\frac{q}{2}}}_1\right.\\
&&
+B^{m_1,{j-\frac{q}{2}}}_2-B^{m_1,m_1+m_2}_2B^{m_1+m_2+1,{j-\frac{q}{2}}}_1)\\
&&\left.\cdot({\mathcal{R}}^{\bot,*}(x))^{m_1+m_2}(\nabla^{\wedge^{0,\bullet}\otimes E}_{w_l}({\mathcal{R}}_.))^{\bot,*}
({\mathcal{R}}^{\bot,*}(x))^{{j-\frac{q}{2}}-1-m_1-m_2}(A_{r\overline{r}l}+A_{l\overline{r}r})\right.\\
&&-\frac{1}{2}B^{m_1,m_1+m_2}_1(B^{m_1+m_2+1,{j-\frac{q}{2}}}_0-B^{m_1+m_2+1,{j-\frac{q}{2}}}_1)\\
&&\cdot({\mathcal{R}}^{\bot,*}(x))^{m_1+m_2}(\nabla^{\wedge^{0,\bullet}\otimes E}_{w_l}({\mathcal{R}}_.))^{\bot,*}
({\mathcal{R}}^{\bot,*}(x))^{{j-\frac{q}{2}}-1-m_1-m_2}A_{e_rle_r}\\
&&+(B^{m_1,m_1+m_2}_1B^{m_1+m_2+1,{j-\frac{q}{2}}}_2-B^{m_1,{j-\frac{q}{2}}}_1\\
&&+B^{m_1,{j-\frac{q}{2}}}_0-B^{m_1,m_1+m_2}_0
B^{m_1+m_2+1,{j-\frac{q}{2}}}_1)\\
&&\cdot({\mathcal{R}}^{\bot,*}(x))^{m_1+m_2}(\nabla^{\wedge^{0,\bullet}\otimes E}_{\overline{w_l}}({\mathcal{R}}_.))^{\bot,*}
({\mathcal{R}}^{\bot,*}(x))^{{j-\frac{q}{2}}-1-m_1-m_2}(A_{\overline{l}\overline{r}r}+A_{\overline{r}\overline{l}r})\\
&&-\frac{1}{2}B^{m_1,m_1+m_2}_0(B^{m_1+m_2+1,{j-\frac{q}{2}}}_0-B^{m_1+m_2+1,{j-\frac{q}{2}}}_1)\\
&&\left.\cdot({\mathcal{R}}^{\bot,*}(x))^{m_1+m_2}(\nabla^{\wedge^{0,\bullet}\otimes E}_{\overline{w_l}}({\mathcal{R}}_.))^{\bot,*}
({\mathcal{R}}^{\bot,*}(x))^{{j-\frac{q}{2}}-1-m_1-m_2}A_{e_r\overline{l}e_r}\right]I_{2j},~~(4.212)
\end{eqnarray*}
When $n=4$, $q=2$ and $j=2$, we have
$$IV^2_a=I_4
({\mathcal{L}}_0^{-1}{\mathcal{O}}^{+2}_{2})P^N
[({\mathcal{L}}_0^{-1}{\mathcal{O}}^{+2}_{3})
({\mathcal{L}}_0^{-1}P^{N,\bot}{\mathcal{O}}''_{1})P^N]^*I_{4}.\eqno(4.213)$$
By (4.202) and (4.158), we know that
$$IV^2_a=
\frac{\sqrt{-1}}{192\pi^3}{\mathcal{R}^\bot}I_{{\rm det}(\overline{W}^*)\otimes E}[\frac{1}{3}\widehat{\Phi}^*(\partial_{{z}_r})
\frac{\partial {\mathcal{R}}^{0,*}}{\partial{\overline{z}_r}}
+\frac{1}{2}\widehat{\Phi}(\partial_{\overline{z}_r})\frac{\partial {\mathcal{R}}^{0,*}}{\partial{{z}_r}}]I_4.\eqno(4.214)$$
By (4.203)-(4.205), we get
$$IV_a=IV^1_a+IV^2_a.\eqno(4.215)$$
\indent Nextly, we compute the term ${\rm IV}_b$.\\
\indent By (4.200), we have
$${\rm IV}_b=\sum_{0\leq m_1+m_2\leq {j-\frac{q}{2}}-1}I_{2j}({\mathcal{L}}_0^{-1}{\mathcal{O}}^{+2}_{2})^{{j-\frac{q}{2}}}P^N
[({\mathcal{L}}_0^{-1}{\mathcal{O}}^{+2}_{2})^{{j-\frac{q}{2}}-1-m_1-m_2}$$
$$\cdot({\mathcal{L}}_0^{-1}{\mathcal{O}}_{1})({\mathcal{L}}_0^{-1}{\mathcal{O}}^{+2}_{2})^{m_2}
({\mathcal{L}}_0^{-1}{\mathcal{O}}^{+2}_{3})({\mathcal{L}}_0^{-1}{\mathcal{O}}^{+2}_{2})^{m_1}P^N]^*I_{2j}.\eqno(4.216)$$
By (4.202), then
\begin{eqnarray*}
&&({\mathcal{L}}_0^{-1}{\mathcal{O}}^{+2}_{2})^{m_2}
({\mathcal{L}}_0^{-1}{\mathcal{O}}^{+2}_{3})({\mathcal{L}}_0^{-1}{\mathcal{O}}^{+2}_{2})^{m_1}{\mathcal{P}}\\
&=&
\frac{1}{(4\pi)^{m_1+m_2+1}}B^{1,m_1+m_2+1}_0
{\mathcal{R}}^{m_2}\frac{\partial{\mathcal{R}}.}{\partial z_{j'}}(0){\mathcal{R}}^{m_1}z_{j'}{\mathcal{P}}\\
&&+\frac{1}{2\pi\times(4\pi)^{m_1+m_2+1}}B^{1,m_1}_0B^{m_1+1,m_1+m_2+1}_1{\mathcal{R}}^{m_2}\frac{\partial{\mathcal{R}}.}{\partial \overline{z}_{j'}}(0){\mathcal{R}}^{m_1}b_{j'}{\mathcal{P}}\\
&&+\frac{1}{(4\pi)^{m_1+m_2+1}}B^{1,m_1+m_2+1}_0{\mathcal{R}}^{m_2}\frac{\partial{\mathcal{R}}.}{\partial \overline{z}_{j'}}(0){\mathcal{R}}^{m_1}\overline{z}'_{j'}{\mathcal{P}}.~~~~~~~~~~~~~~~~~~~~~~~~~~~~~~(4.217)
\end{eqnarray*}
By (4.97) and (4.217) and direct computations, we know that after integration
\begin{eqnarray*}
&&(A_1+A_2+A_3+A_4)({\mathcal{L}}_0^{-1}{\mathcal{O}}^{+2}_{2})^{m_2}
({\mathcal{L}}_0^{-1}{\mathcal{O}}^{+2}_{3})({\mathcal{L}}_0^{-1}{\mathcal{O}}^{+2}_{2})^{m_1}{\mathcal{P}}\\
&=&\frac{1}{2\pi\times(4\pi)^{m_1+m_2+1}}B^{1,m_1}_0B^{m_1+1,m_1+m_2+1}_1({\mathcal{R}}^\bot)^{m_2}\frac{\partial{\mathcal{R}}^\bot.}{\partial \overline{z}_{j''}}(0)({\mathcal{R}}^\bot)^{m_1}(b_{j''}z_{j'}+2\delta_{j'j''}){\mathcal{P}}\\
&&\times[\frac{4}{3}(\partial_{{z}_{j'} }R^L)_{x_0}(\partial_{z_\alpha},\partial_{\overline{z}_\alpha})+
\frac{4}{3}(\partial_{{z}_{\alpha} }R^L)_{x_0}(\partial_{z_{j'}},\partial_{\overline{z}_\alpha})-
\frac{1}{3}(\partial_{\alpha}R^L)_{x_0}(\partial_{z_{j'}},e_{\alpha})]\\
&&-\frac{1}{3\pi\times(4\pi)^{m_1+m_2+1}}B^{1,m_1}_0B^{m_1+1,m_1+m_2+1}_1({\mathcal{R}}^\bot)^{m_2}\frac{\partial{\mathcal{R}}^\bot.}{\partial \overline{z}_{j'}}(0)({\mathcal{R}}^\bot)^{m_1}\\
&&\cdot(\partial_{{z}_{\beta} }R^L)_{x_0}(\partial_{\overline{z}_\alpha},\partial_{z_l})(2\delta_{lj'}b_\alpha z_\beta+4\delta_{lj'}\delta_{\alpha\beta}){\mathcal{P}}\\
&&+\frac{1}{(4\pi)^{m_1+m_2+1}}B^{1,m_1+m_2+1}_0({\mathcal{R}}^\bot)^{m_2}\frac{\partial{\mathcal{R}}^\bot.}{\partial {z}_{j'}}(0)({\mathcal{R}}^\bot)^{m_1}b_{\alpha}z_{j'}{\mathcal{P}}\\
&&\times[\frac{2}{3\pi}(\partial_{\overline{z}_r }R^L)_{x_0}(\partial_{z_r},\partial_{\overline{z}_\alpha})
+\frac{2}{3\pi}(\partial_{\overline{z}_\alpha }R^L)_{x_0}(\partial_{z_r},\partial_{\overline{z}_r})-
\frac{1}{6\pi}(\partial_{j''}R^L)_{x_0}(\partial_{\overline{z}_\alpha},e_{j''})]\\
&&+\frac{1}{3\pi\times(4\pi)^{m_1+m_2+1}}B^{1,m_1}_0B^{m_1+1,m_1+m_2+1}_1({\mathcal{R}}^\bot)^{m_2}\frac{\partial{\mathcal{R}}^\bot.}{\partial \overline{z}_{j'}}(0)({\mathcal{R}}^\bot)^{m_1}\\
&&\cdot(\partial_{{z}_{j''} }R^L)_{x_0}(\partial_{{z}_\alpha},\partial_{\overline{z}_l})(b_lb_{j'}z_\alpha z_{j''}
+2\delta_{\alpha j'}b_lz_{j''}+2\delta_{j'j''}b_lz_\alpha){\mathcal{P}}\\
&&+\frac{1}{3\pi\times(4\pi)^{m_1+m_2+1}}B^{1,m_1+m_2+1}_0({\mathcal{R}}^\bot)^{m_2}\frac{\partial{\mathcal{R}}^\bot.}{\partial {z}_{j'}}(0)({\mathcal{R}}^\bot)^{m_1}\\
&&\cdot(\partial_{\overline{z}_{r} }R^L)_{x_0}(\partial_{{z}_\alpha},\partial_{\overline{z}_l})b_lb_{r}z_\alpha z_{j'}
{\mathcal{P}}~~~~~~~~~~~~~~~~~~~~~~~~~~~~~~~~~~~~~~~~~~~~~~~~~~~~~~~(4.218)
\end{eqnarray*}
Set
$${\rm IV}^1_b=\sum_{0\leq m_1+m_2\leq {j-\frac{q}{2}}-1}I_{2j}({\mathcal{L}}_0^{-1}{\mathcal{O}}^{+2}_{2})^{{j-\frac{q}{2}}}P^N
[({\mathcal{L}}_0^{-1}{\mathcal{O}}^{+2}_{2})^{{j-\frac{q}{2}}-1-m_1-m_2}$$
$$\cdot({\mathcal{L}}_0^{-1}(A_1+A_2+A_3+A_4))({\mathcal{L}}_0^{-1}{\mathcal{O}}^{+2}_{2})^{m_2}
({\mathcal{L}}_0^{-1}{\mathcal{O}}^{+2}_{3})({\mathcal{L}}_0^{-1}{\mathcal{O}}^{+2}_{2})^{m_1}P^N]^*I_{2j}.\eqno(4.219)$$
$${\rm IV}^2_b=\sum_{0\leq m_1+m_2\leq {j-\frac{q}{2}}-1}I_{2j}({\mathcal{L}}_0^{-1}{\mathcal{O}}^{+2}_{2})^{{j-\frac{q}{2}}}P^N
[({\mathcal{L}}_0^{-1}{\mathcal{O}}^{+2}_{2})^{{j-\frac{q}{2}}-1-m_1-m_2}$$
$$\cdot({\mathcal{L}}_0^{-1}A_5)({\mathcal{L}}_0^{-1}{\mathcal{O}}^{+2}_{2})^{m_2}
({\mathcal{L}}_0^{-1}{\mathcal{O}}^{+2}_{3})({\mathcal{L}}_0^{-1}{\mathcal{O}}^{+2}_{2})^{m_1}P^N]^*I_{2j}.\eqno(4.220)$$
$${\rm IV}^3_b=\sum_{0\leq m_1+m_2\leq {j-\frac{q}{2}}-1}I_{2j}({\mathcal{L}}_0^{-1}{\mathcal{O}}^{+2}_{2})^{{j-\frac{q}{2}}}P^N
[({\mathcal{L}}_0^{-1}{\mathcal{O}}^{+2}_{2})^{{j-\frac{q}{2}}-1-m_1-m_2}$$
$$\cdot({\mathcal{L}}_0^{-1}{\mathcal{O}}''_{1})({\mathcal{L}}_0^{-1}{\mathcal{O}}^{+2}_{2})^{m_2}
({\mathcal{L}}_0^{-1}{\mathcal{O}}^{+2}_{3})({\mathcal{L}}_0^{-1}{\mathcal{O}}^{+2}_{2})^{m_1}P^N]^*I_{2j}.\eqno(4.221)$$
Then $$IV_b={\rm IV}^1_b+{\rm IV}^2_b+{\rm IV}^3_b.\eqno(4.222)$$
By (4.218), we get

\begin{eqnarray*}
{\rm IV}^1_b&=&\frac{1}{2^{j-\frac{q}{2}}({j-\frac{q}{2}})!(4\pi)^{2j-q+1}\times 3\pi}\sum_{0\leq m_1+m_2\leq {j-\frac{q}{2}}-1}
I_{2j}({\mathcal{R}}^\bot)(x)^{{j-\frac{q}{2}}}\\
&&\cdot I_{{\rm det}(\overline{W}^*)\otimes E}\left[(-2B^{1,m_1}_0B^{m_1+1,m_1+m_2+1}_1B^{m_1+m_2+1,{j-\frac{q}{2}}}_0\right.\\
&&+2B^{1,m_1}_0B^{m_1+1,m_1+m_2+1}_1B^{m_1+m_2+1,{j-\frac{q}{2}}}_1\\
&&
-B^{1,m_1}_0B^{m_1+1,m_1+m_2+1}_1B^{m_1+m_2+1,{j-\frac{q}{2}}}_2
+B^{1,m_1}_0B^{m_1+1,m_1+m_2+1}_1B^{m_1+m_2+1,{j-\frac{q}{2}}}_1)\\
&&\cdot({\mathcal{R}}^{\bot,*}(x))^{m_1}(\nabla^{\wedge^{0,\bullet}\otimes E}_{\overline{w_l}}({\mathcal{R}}_.))^{\bot,*}
({\mathcal{R}}^{\bot,*}(x))^{{j-\frac{q}{2}}-1-m_1}A_{\overline{r}\overline{l}r}\\
&&-(B^{1,m_1}_0B^{m_1+1,m_1+m_2+1}_1B^{m_1+m_2+1,{j-\frac{q}{2}}}_0
-B^{1,m_1}_0B^{m_1+1,m_1+m_2+1}_1B^{m_1+m_2+1,{j-\frac{q}{2}}}_1\\
&&+B^{1,m_1}_0B^{m_1+1,m_1+m_2+1}_1B^{m_1+m_2+1,{j-\frac{q}{2}}}_2
-B^{1,m_1}_0B^{m_1+1,m_1+m_2+1}_1B^{m_1+m_2+1,{j-\frac{q}{2}}}_1)\\
&&\cdot({\mathcal{R}}^{\bot,*}(x))^{m_1}(\nabla^{\wedge^{0,\bullet}\otimes E}_{\overline{w_l}}({\mathcal{R}}_.))^{\bot,*}
({\mathcal{R}}^{\bot,*}(x))^{{j-\frac{q}{2}}-1-m_1}A_{\overline{l}\overline{r}r}\\
&&+(B^{1,m_1+m_2+1}_0B^{m_1+m_2+1,{j-\frac{q}{2}}}_1-B^{1,m_1+m_2+1}_0B^{m_1+m_2+1,{j-\frac{q}{2}}}_2)\\
&&\cdot({\mathcal{R}}^{\bot,*}(x))^{m_1}(\nabla^{\wedge^{0,\bullet}\otimes E}_{{w_l}}({\mathcal{R}}_.))^{\bot,*}
({\mathcal{R}}^{\bot,*}(x))^{{j-\frac{q}{2}}-1-m_1}(A_{r\overline{r}l}+A_{l\overline{r}r})\\
&&+\frac{1}{2}B^{1,m_1}_0B^{m_1+1,m_1+m_2+1}_1(B^{m_1+m_2+1,j}_0-B^{m_1+m_2+1,{j-\frac{q}{2}}}_1)\\
&&\cdot({\mathcal{R}}^{\bot,*}(x))^{m_1}(\nabla^{\wedge^{0,\bullet}\otimes E}_{\overline{w_l}}({\mathcal{R}}_.))^{\bot,*}
({\mathcal{R}}^{\bot,*}(x))^{{j-\frac{q}{2}}-1-m_1}A_{e_r\overline{l}e_r}\\
&&-\frac{1}{2}B^{1,m_1+m_2+1}_0B^{m_1+m_2+1,{j-\frac{q}{2}}}_1
({\mathcal{R}}^{\bot,*}(x))^{m_1}\\
&&\left.\cdot(\nabla^{\wedge^{0,\bullet}\otimes E}_{{w_l}}({\mathcal{R}}_.))^*
({\mathcal{R}}^{\bot,*}(x))^{{j-\frac{q}{2}}-1-m_1}A_{e_rle_r}\right],~~~~~~~~~~~~~~~~~~~~~~~~~~~~~~~(4.223)
\end{eqnarray*}
By (4.220), when $n=4$, $q=2$ and $j=2$
$${\rm IV}^2_b=I_{4}({\mathcal{L}}_0^{-1}{\mathcal{O}}^{+2}_{2})P^N[
({\mathcal{L}}_0^{-1}A_5)
({\mathcal{L}}_0^{-1}{\mathcal{O}}^{+2}_{3})P^N]^*I_{4}.\eqno(4.224)$$
Then (4.102) and (4.202), we obtain
$$IV^2_b=
\frac{1}{8\pi^4}{\mathcal{R}^\bot}I_{{\rm det}(\overline{W}^*)\otimes E}[\frac{-1}{144}(\partial_{{z}_r }R^L)_{x_0}
(\partial_{\overline{z}_\alpha},\partial_{z_\alpha})
\frac{\partial {\mathcal{R}}^{\bot,*}}{\partial{\overline{z}_r}}$$
$$+\frac{1}{192}(\partial_{j'}R^L)_{x_0}(\partial_{{z}_r},e_{j'})\frac{\partial {\mathcal{R}}^{\bot,*}}{\partial{\overline{z}_r}}]I_4.\eqno(4.225)$$
When $n=4$, $q=2$ and $j=2$,
$$IV^3_b=0.\eqno(4.226)$$
Set
$${\rm IV}^1_c=\sum_{0\leq m_1\leq {j-\frac{q}{2}}-1}\sum_{0\leq m_2\leq {j-\frac{q}{2}}}I_{2j}
({\mathcal{L}}_0^{-1}{\mathcal{O}}^{+2}_{2})^{m_1}
({\mathcal{L}}_0^{-1}{\mathcal{O}}^{+2}_{3})({\mathcal{L}}_0^{-1}{\mathcal{O}}^{+2}_{2})^{{j-\frac{q}{2}}-1-m_1}P^N$$
$$\cdot({\mathcal{O}}^{-2}_{2}{\mathcal{L}}_0^{-1})^{m_2}
({\mathcal{O}}'_{1}{\mathcal{L}}_0^{-1})({\mathcal{O}}^{-2}_{2}{\mathcal{L}}_0^{-1})^{{j-\frac{q}{2}}-m_2}
I_{2j}.\eqno(4.227)$$
$${\rm IV}^2_c=\sum_{0\leq m_1\leq {j-\frac{q}{2}}-1}\sum_{0\leq m_2\leq {j-\frac{q}{2}}}I_{2j}
({\mathcal{L}}_0^{-1}{\mathcal{O}}^{+2}_{2})^{m_1}
({\mathcal{L}}_0^{-1}{\mathcal{O}}^{+2}_{3})({\mathcal{L}}_0^{-1}{\mathcal{O}}^{+2}_{2})^{{j-\frac{q}{2}}-1-m_1}P^N$$
$$\cdot({\mathcal{O}}^{-2}_{2}{\mathcal{L}}_0^{-1})^{m_2}
({\mathcal{O}}''_{1}{\mathcal{L}}_0^{-1})({\mathcal{O}}^{-2}_{2}{\mathcal{L}}_0^{-1})^{{j-\frac{q}{2}}-m_2}
I_{2j}.\eqno(4.228)$$
Then $$IV_c=IV^1_c+IV^2_c.\eqno(4.229)$$
\indent By (4.24), (4.187) and (4.217), we get
\begin{eqnarray*}
{\rm IV}^1_c&=&\sum_{0\leq m_1\leq {j-\frac{q}{2}}-1}\sum_{0\leq m_2\leq {j-\frac{q}{2}}}\frac{1}{3\pi\times (4\pi)^{2j-q+1}}\\
&&\cdot(B^{1,{j-\frac{q}{2}}}_0-B^{1,{j-\frac{q}{2}}-m_1-1}_0B^{{j-\frac{q}{2}}-m_1,{j-\frac{q}{2}}}_1)
({\mathcal{R}}^\bot)(x)^{m_1}\\
&&\cdot(\nabla^{\wedge^{0,\bullet}\otimes E}_{\overline{w_l}}({\mathcal{R}}_.))^\bot
({\mathcal{R}}^\bot)(x)^{{j-\frac{q}{2}}-1-m_1} I_{{\rm det}(\overline{W}^*)\otimes E}({\mathcal{R}}^{\bot,*}(x))^{{j-\frac{q}{2}}}\\
&&\cdot\left[(B^{1,{j-\frac{q}{2}}-m_2}_0B^{{j-\frac{q}{2}}-m_2,{j-\frac{q}{2}}}_0-2
B^{1,{j-\frac{q}{2}}-m_2}_0B^{{j-\frac{q}{2}}-m_2,{j-\frac{q}{2}}}_1\right.\\
&&+B^{1,{j-\frac{q}{2}}-m_2}_0
B^{{j-\frac{q}{2}}-m_2,{j-\frac{q}{2}}}_2)A_{\overline{l}r\overline{r}}\\
&&-(B^{1,{j-\frac{q}{2}}-m_2}_0B^{{j-\frac{q}{2}}-m_2,{j-\frac{q}{2}}}_1-B^{1,{j-\frac{q}{2}}-m_2}_0
B^{{j-\frac{q}{2}}-m_2,{j-\frac{q}{2}}}_2)A_{\overline{r}r\overline{l}}\\
&&\left.-\frac{1}{2}(B^{1,{j-\frac{q}{2}}-m_2}_0B^{{j-\frac{q}{2}}-m_2,{j-\frac{q}{2}}}_0-B^{1,{j-\frac{q}{2}}-m_2}_0
B^{j-m_2,{j-\frac{q}{2}}}_1)A_{e_r\overline{l}e_r}\right],~~~~~~~~(4.230)
\end{eqnarray*}
By (4.228), when $n=4$, $j=2$ and $q=2$, one has
$${\rm IV}^2_c=I_4
({\mathcal{L}}_0^{-1}{\mathcal{O}}^{+2}_{3})P^N
({\mathcal{L}}_0^{-1}{\mathcal{O}}^{+2}_{2}
{\mathcal{L}}_0^{-1}{\mathcal{O}}''_{1}P^N)^*I_4.
\eqno(4.231)$$
Similar to $IV^2_b$, we obtain
$$IV^2_c=
\frac{5\sqrt{-1}}{24\times 144\pi^3}I_4
\frac{\partial {\mathcal{R}}^{\bot}}{\partial{\overline{z}_r}}I_{{\rm det}(\overline{W}^*)\otimes E}
\widehat{\Phi}^*(\partial_{z_r}) {\mathcal{R}}^{0,*}I_4.
\eqno(4.232)$$\\
As in [PZ, p.20], we may write the above formulas in the intrinsic way. When $n=4$, $j=2$ and $q=2$, by (4.1), (4.3), (4.4)-(4.8), (4.11), (4.34), (4.37), (4.76), (4.79)-(4.83), (4.86), (4.151)-(4.153),
(4.159), (4.162)-(4.172), (4.182), (4.183), (4.189), (4.194)-(4.198), (4.212), (4.214), (4.215), (4.222), (4.223), (4.225), (4.226),
(4.229), (4.230), (4.232), we prove Theorem 1.4.~~$\Box$\\

 \noindent {\bf Acknowledgements.} This work
was supported by NSFC No.11271062 and NCET-13-0721. The author is indebted to referees for their careful reading and helpful comments.\\

\noindent{\large \bf References}\\

\noindent[Bi]Bismut, J.: A local index theorem for non K${\rm \ddot{a}}$hler manifolds. Math. Ann. 284, 681-699 (1989)

\noindent[Ca]Catlin, D.: The Bergman kernel and a theorem of Tian. In: Analysis and Geometry in Several Complex
Variables, Katata, 1997. Trends Math., pp. 1-23. Birkhauser, Boston (1999)

\noindent[CS]Charbonneau B., Stern M.: Asymptotic Hodge theory of vector bundles. arXiv:\\
1111.0591 (2011)

\noindent[DLM]Dai, X., Liu, K., Ma, X.: On the asymptotic expansion of Bergman kernel. J. Differ. Geom. 72, 1-41
(2006). Announced in C. R. Math. Acad. Sci. Paris 339, 193-198 (2004)

\noindent[LuW1]Lu, W.: The second coefficient of the asymptotic expansion
of the Bergman kernel of the Hodge-Dolbeault operator. J. Geom. Anal., DOI 10.1007/s12220-013-9412-y.

\noindent[LuW2]Lu, W.: Morse Inequalities and Bergman Kernels. Doctor thesis, (2013)

\noindent[LuZ]Lu, Z.: On the lower order terms of the asymptotic expansion of Tian-Yau-Zelditch. Am. J.Math. 122,
235-273 (2000)

\noindent[MM06]Ma, X., Marinescu, G.: The first coefficients of the asymptotic expansion of the Bergman kernel of
the Spinc Dirac operator. Int. J. Math. 17, 737-759 (2006)

\noindent[MM07]Ma, X., Marinescu, G.: Holomorphic Morse Inequalities and Bergman Kernels. Progr.Math., vol. 254.
Birkhauser, Basel (2007)

\noindent[MM08]Ma, X., Marinescu, G.: Generalized Bergman kernels on symplectic manifolds. Adv. Math. 217,
1756-1815 (2008)

\noindent[MM12]Ma, X., Marinescu, G.: Berezin-Toeplitz quantization on K${\rm \ddot{a}}$hler manifolds. J. Reine Angew. Math.
662, 1-56 (2012)

\noindent[PZ]Puchol, M., Zhu, J.: The first terms in the expansion of the Bergman kernel in higher degrees. arXiv:1210.1717 (2012)

\noindent[Ru]Ruan,W.: Canonical coordinates and Bergman metrics. Commun. Anal. Geom. 6, 589-631 (1998)

\noindent[Ti]Tian, G.: On a set of polarized K${\rm \ddot{a}}$hler metrics on algebraic manifolds. J. Differ. Geom. 32, 99-130
(1990)

\noindent[Wa]Wang, X.: Canonical metrics on stable vector bundles. Commun. Anal. Geom. 13, 253-285 (2005)

\noindent[Ze]Zelditch, S.: Szeg${\rm \ddot{o}}$ kernels and a theorem of Tian. Int. Math. Res. Not. 6, 317-331 (1998)

 \indent{\it School of Mathematics and Statistics, Northeast Normal University, Changchun Jilin, 130024, China }\\
 \indent E-mail: {\it wangy581@nenu.edu.cn}\\

\end{document}